%% file: main.tex
\documentclass{siamart220329}
\newsiamremark{example}{Example}

\usepackage{import}
\subimport{settings/}{_all.tex}

\ifpdf{}
\hypersetup{
  pdftitle={Approximating Higher-Order Derivative Tensors Using Secant Updates},
  pdfauthor={Karl Welzel and Raphael A. Hauser}
}
\fi

\headers{Approximating Higher-Order Derivative Tensors}{Karl Welzel and Raphael A. Hauser}

\title{Approximating Higher-Order Derivative Tensors Using Secant Updates\thanks{Received by the editors DATE.
\funding{This work is supported by the Hong Kong Innovation and Technology Commission (InnoHK Project CIMDA)}}}

\author{Karl Welzel\thanks{Mathematical Institute, University of Oxford, Woodstock Road, Oxford OX2 6GG, United Kingdom
  (\email{welzel@maths.ox.ac.uk}, \email{hauser@maths.ox.ac.uk}).}
\and Raphael A. Hauser\footnotemark[2]}

\begin{document}

\maketitle

\begin{abstract}
    Quasi-Newton methods employ an update rule that gradually improves the Hessian approximation using the already available gradient evaluations.
    We propose higher-order secant updates which generalize this idea to higher-order derivatives, approximating for example third derivatives (which are tensors) from given Hessian evaluations.
    Our generalization is based on the observation that quasi-Newton updates are least-change updates satisfying the secant equation, with different methods using different norms to measure the size of the change. We present a full characterization for least-change updates in weighted Frobenius norms (satisfying an analogue of the secant equation) for derivatives of arbitrary order.
    Moreover, we establish convergence of the approximations to the true derivative under standard assumptions and explore the quality of the generated approximations in numerical experiments.
\end{abstract}



\subimport{}{work.tex}

\bibliographystyle{siamplain}
\bibliography{citations}

\end{document}

%% file: settings/_all.tex
\subimport{.}{math}

\subimport{.}{appearance}
\subimport{.}{plots}

%% file: settings/math.tex

\usepackage{mathtools}
\usepackage{amsfonts}
\usepackage{bm}  
\usepackage{mleftright} \mleftright{}  
\usepackage{braket}  
\usepackage{siunitx}  


\let\e\varepsilon{}               

\let\N\Natural{}
{}
\let\R\Real{}
{}

\newcommand*\diff{\mathop{}\!\mathrm{d}}	


\DeclareMathOperator{\argmin}{arg\,min}


\DeclareMathOperator{\rank}{rank}


\newcommand{\sym}{\mathrm{sym}} 



\newcommand{\deq}{\coloneqq}	

\let\brack\undefined

\DeclarePairedDelimiter{\abs}{\lvert}{\rvert}
\DeclarePairedDelimiter{\norm}{\lVert}{\rVert}
\DeclarePairedDelimiter{\paren}{\lparen}{\rparen}
\DeclarePairedDelimiter{\brack}{\lbrack}{\rbrack}
\DeclarePairedDelimiterX{\innerprod}[2]{\langle}{\rangle}{#1,#2}


\newcommand{\Abs}[1]{\abs*{#1}}
\newcommand{\Norm}[1]{\norm*{#1}}
\newcommand{\Paren}[1]{\paren*{#1}}
\newcommand{\Brack}[1]{\brack*{#1}}
\newcommand{\Innerprod}[2]{\innerprod*{#1}{#2}}


\newcommand{\vek}[1]{\bm{#1}}
\newcommand{\mat}[1]{\bm{#1}}
\newcommand{\tns}[1]{\textsf{\textbf{\textsl{#1}}}}

%% file: settings/appearance.tex
\usepackage{booktabs}  
\usepackage{multirow}
\usepackage{array}  
\usepackage{makecell} 
\usepackage{siunitx} 

\usepackage[shortlabels]{enumitem}

%% file: settings/plots.tex
\usepackage{pgf}

%% file: work.tex
\section{Introduction}\label{sec:introduction}

The incredible success of quasi-Newton methods is largely based on the fact that the rules for updating the Hessian approximations are able to extract crucial second-order information from gradient evaluations.
Unlike finite difference methods they do so without direct control over the evaluation points.
Rather, quasi-Newton rules are designed to handle evaluations of the first derivative at points that are generated by an extraneous process and produce best-effort approximations of the second derivative.
We propose generalizations of these rules that mimic the quasi-Newton approach but approximate \(p\)th derivatives from given evaluations of \((p-1)\)st derivatives for any \(p \geq 2\).

A key motivation for our work is the recent theoretical advances in higher-order optimization methods for unconstrained problems.
Birgin et al.\ \cite{birgin_worst-case_2017} showed that if the objective function \(f\) is \(p\) times continuously differentiable with a Lipschitz continuous \(p\)th derivative and an oracle to compute the first \(p\) derivatives at any point is provided, then an algorithm exists that finds a point with \(\norm{\nabla f(\vek{x})} \leq \e\) in at most \(O(\e^{-(p+1)/p})\) oracle calls.
This result generalized the known cases for \(p=1\) \cite{nesterov_introductory_2004} and \(p=2\) \cite{nesterov_cubic_2006} and was later extended by Cartis, Gould and Toint \cite{cartis_sharp_2020}, who also proved that this bound is sharp for algorithms that minimize regularized Taylor models in each step such as the one used in \cite{birgin_worst-case_2017}.
Simply put, access to more derivatives improves the performance of optimization algorithms.
Approximate higher-order derivatives might provide a way to achieve this improved performance without additional derivative evaluations.
In this paper however, we focus on properties of the updates and on results on the accuracy of the approximations that can be derived without detailed knowledge of how the evaluation points are generated.

This paper will be structured as follows: After introducing the tensor notation we use in \cref{sec:notation}, we will derive the tensor analogues of quasi-Newton updates in \cref{sec:derivation} and give a full characterization of these updates in \cref{sec:characterization}.
The characterization includes an explicit formula and a useful recursive relationship between successive approximations.
Moreover, we will show that these updates exhibit a certain low-rank structure.
\Cref{sec:convergence} contains results on the convergence of the approximations to the exact derivative in the limit under certain conditions on the steps and comparisons of these results to the ones found in the literature on convergence of quasi-Newton matrices.
Lastly, in \cref{sec:numerical-experiments} we present limited numerical experiments to verify the behaviour predicted by the theory and discuss numerical limitations of this method.

\section{Notation}\label{sec:notation}

Along with the notation for tensors and higher-order derivatives that we will use, which is based on the one in \cite{cartis_sharp_2020}, this section also introduces some standard definitions and properties of tensors.
For more information please refer to \cite{kolda_tensor_2009} for an introduction to tensors from an applied perspective and to \cite{hackbusch_tensor_2012} for an in-depth discussion of abstract tensor spaces.

A \emph{\(p\)-tensor} \(\tns{T}\) of dimensions \(n_1 \times \dots \times n_p\) is a multilinear map \(\R^{n_1} \times \dots \times \R^{n_p}\to \R\), so that its evaluation
\begin{equation}
    \tns{T}[\vek{s}_1, \dots, \vek{s}_p], \quad \vek{s}_i \in \R^{n_i}
\end{equation}
is linear in each component.
We refer to these \(p\) components as the \emph{modes} of the tensor and denote the space of such \(p\)-tensors by \(\R^{n_1 \otimes \dots \otimes n_p}\).
If \(n_1 = \dots = n_p = n\) the space is denoted \(\R^{\otimes^p n}\) and \(\tns{T}[\vek{s}, \dots, \vek{s}]\) is abbreviated as \(\tns{T}[\vek{s}]^p\).
The notation also allows to apply the tensor to \(q < p\) vectors, which then results in a \((p-q)\)-tensor.
Moreover, we define the application of matrices \(\mat{W}_1, \dots, \mat{W}_p\) of appropriate dimensions to a tensor by
\begin{equation}
    \big(\tns{T}[\mat{W}_1, \dots, \mat{W}_p]\big)[\vek{s}_1, \dots, \vek{s}_p] = \tns{T}[\mat{W}_1 \vek{s}_1, \dots, \mat{W}_p \vek{s}_p].
\end{equation}

The outer product of a \(p_1\)-tensor \(\tns{T}_1\) with a \(p_2\)-tensor \(\tns{T}_2\) is defined as
\begin{equation}
    (\tns{T}_1 \otimes \tns{T}_2)[\vek{s}_1, \dots, \vek{s}_{p_1 + p_2}] = \tns{T}_1[\vek{s}_1, \dots, \vek{s}_{p_1}] \cdot \tns{T}_2[\vek{s}_{p_1 + 1}, \dots, \vek{s}_{p_1 + p_2}].
\end{equation}
In particular, tensors of the form \(\tns{T} = \vek{v}_1 \otimes \dots \otimes \vek{v}_p\) for vectors \(\vek{v}_1, \dots, \vek{v}_p \in \R^n\), i.e. those where \(\tns{T}[\vek{s}_1, \dots, \vek{s}_p] = \prod_{i = 1}^n \vek{v}_i^T \vek{s}_i\), are called \emph{elementary} or \emph{rank-one} tensors.
If all vectors are the same (\(\vek{v}_1 = \dots = \vek{v}_p = \vek{v}\)), we abbreviate the notation above to \(\otimes^p \vek{v}\).
(Note that this notation is slightly inconsistent since 1-tensors should be row vectors, but are represented by standard column vectors to simplify the notation.
This is why \(\vek{v}[\mat{W}] = \mat{W}^T \vek{v}\).)

For any \(\tns{T} \in \R^{\otimes^p n}\) and any permutation \(\sigma \in S_p\), let \(\sigma(\tns{T}) \in \R^{\otimes^p n}\) be defined by
\begin{equation}
    \sigma(\tns{T})[\vek{s}_1, \dots, \vek{s}_p] = \tns{T}[\vek{s}_{\sigma(1)}, \dots, \vek{s}_{\sigma(p)}].
\end{equation}
If \(\sigma(\tns{T}) = \tns{T}\) for all \(\sigma \in S_p\), then \(\tns{T}\) is called \emph{symmetric}.
The space of all symmetric \(p\)-tensors is denoted \(\R^{\otimes^p n}_{\sym}\).
The projection of \(\R^{\otimes^p n}\) onto \(\R^{\otimes^p n}_{\sym}\) is given by
\begin{equation}
    P_{\sym}(\tns{T}) = \frac{1}{p!} \sum_{\sigma \in S_p} \sigma(\tns{T}),
\end{equation}
see \cite[Proposition 3.76]{hackbusch_tensor_2012}.

Just like matrices, tensors are fully characterized by their actions on basis vectors.
This can be used to represent a \(p\)-tensor \(\tns{T} \in \R^{\otimes^p n}\) as a \(p\)-dimensional array \((t_{i_1, \dots, i_p})_{1 \leq i_j \leq n, 1 \leq j \leq p}\) where
\begin{equation}
    t_{i_1, \dots, i_p} = \tns{T}[\vek{e}_{i_1}, \dots, \vek{e}_{i_p}] \quad \text{and} \quad \tns{T} = \sum_{i_1, \dots, i_p = 1}^n t_{i_1, \dots, i_p} \vek{e}_{i_1} \otimes \dots \otimes \vek{e}_{i_p}.
\end{equation}
There is a Frobenius inner product and a corresponding norm on these \(p\)-dimensional arrays and by extension on \(\R^{\otimes^p n}\), which we will denote by \(\innerprod{\tns{T}_1}{\tns{T}_2}_F\) and \(\norm{\tns{T}}_F\).
Note that \(\innerprod{\tns{T}}{\vek{s}_1 \otimes \dots \otimes \vek{s}_p}_F = \tns{T}[\vek{s}_1, \dots, \vek{s}_p]\).
Another norm is the one induced by the 2-norm on \(\R^n\) (Hackbusch \cite{hackbusch_tensor_2012} calls it the injective norm) which is defined by
\begin{equation}
    \norm{\tns{T}}_2 = \max_{\norm{\vek{s}_i}_2 = 1,\, 1 \leq i \leq p} \abs{\tns{T}[\vek{s}_1, \dots, \vek{s}_p]}.
\end{equation}
Both norms are invariant under orthogonal transformations, so that if \(\mat{Q}_1, \dots, \mat{Q}_p \in \R^{n \times n}\) are orthogonal matrices, then \(\tns{T}[\mat{Q}_1, \dots, \mat{Q}_p]\) has the same Frobenius- and 2-norm as \(\tns{T}\).
As for matrices, the 2-norm is bounded by the Frobenius norm:
\begin{equation}
    \norm{\tns{T}}_2
    = \max_{\norm{\vek{s}_i}_2 = 1,\, 1 \leq i \leq p} \abs{\innerprod{\tns{T}}{\vek{s}_1 \otimes \dots \otimes \vek{s}_p}_F}
    \leq \norm{\tns{T}}_F.
\end{equation}
In the last inequality we used Cauchy-Schwarz and the fact that \(\norm{\vek{s}_1 \otimes \dots \otimes \vek{s}_p}_F = 1\) if all \(\vek{s}_i\) have unit 2-norm.

The \emph{rank} (or \emph{CP rank}) of a tensor is defined as the minimum number \(r\) such that the tensor can be represented as a sum of \(r\) rank-one tensors and denoted as \(\rank(\tns{T})\).
This notion of rank generalizes the familiar notion of the rank of a matrix.
It is however not the only generalization of matrix rank.
Where for matrices the row and column rank always coincide, this is no longer true for tensors.
For each mode \(i\) let \(r_i\) be the dimension of the subspace spanned by the fibers (the analogue of matrix rows and columns) of mode \(i\),\footnote{Equivalently, \(r_i\) is the rank of matrix unfolding of the tensor with dimensions \((\prod_{k \neq i} n_k) \times n_i\).} then the \emph{multilinear rank} (or \emph{Tucker rank}) of \(\tns{T}\) is the tuple \((r_1, \dots, r_p)\).
For example, the multilinear rank of a rank-one tensor is \((1, \dots, 1)\) and the multilinear rank of a generic \(\R^{\otimes^p n}\) tensor is \((n, \dots, n)\).

Using this setup we can now introduce the notation for higher order derivatives.
Let \(f \colon \R^n \to \R\) be a smooth function.
The \(p\)th total derivative of \(f\) at \(\vek{x} \in \R^n\) is denoted \(D^p f(\vek{x})\) and is recursively defined as the total derivative of \(D^{p-1} f\) where \(D^0 f = f\).
This gives a chain of linear maps which we can regard as one multilinear map
\begin{equation}
    D^p f(\vek{x}) \colon \R^n \to \Paren{\R^n \to \Paren{\dots \Paren{\R^n \to \R}}} = \underbrace{\R^n \times \dots \times \R^n}_{p \text{ times}} \to \R
\end{equation}
making \(D^p f(\vek{x})\) a \(p\)-tensor of dimensions \(n \times \dots \times n\).
By this definition the evaluation of this \(p\)-tensor \(D^p f(\vek{x})[\vek{s}_1, \dots, \vek{s}_p]\) is equal to the directional derivative of \(f\) at \(\vek{x}\) along directions \(\vek{s}_1, \dots, \vek{s}_p \in \R^n\).
For example, denoting the gradient by \(\nabla f\) and the Hessian by \(\nabla^2 f\) we have
\begin{equation}
    D^1 f(\vek{x})[\vek{s}_1] = \nabla f(\vek{x})^T \vek{s}_1 \quad \text{and} \quad D^2 f(\vek{x})[\vek{s}_1, \vek{s}_2] = \vek{s}_1^T \nabla^2 f(\vek{x}) \vek{s}_2.
\end{equation}
Moreover, \(D^p f(\vek{x})\) is symmetric, because partial derivatives commute (Schwarz's theorem).
The \(p\)th-order Taylor expansion of \(f\) at \(\vek{x}\) evaluated at an offset \(\vek{s} \in \R^n\) can be expressed in this notation as
\begin{equation}
    T_{f, p}(\vek{x}, \vek{s}) = \sum_{k=0}^p \frac{1}{k!} D^k f(\vek{x})[\vek{s}]^k \approx f(\vek{x} + \vek{s}).
\end{equation}

\section{Derivation}\label{sec:derivation}

We now turn to our derivation of higher-order secant updates by first introducing a few important quasi-Newton updates.
The most well-known update rule for quasi-Newton methods is the BFGS method, which is often described as a rank-two update to the current Hessian approximation.
To motivate the generalization in this paper we take a different view and describe BFGS and similar methods as choosing minimal updates that satisfy the secant equation \cite{dennis_numerical_1996,nocedal_numerical_2006}.
Let \(f \in C^2(\R^n)\) be a twice continuously differentiable function with gradient \(\nabla f\) and Hessian \(\nabla^2 f\).
Assume that we are given some sequence of points \(\vek{x}_k \in \R^n\) for \(k \in \N\) (possibly from minimizing \(f\)), the gradients \(\nabla f(\vek{x}_k)\) at each iterate and some symmetric initial Hessian approximation \(\mat{B}_0 \in \R^{n \times n}_{\sym}\).
Quasi-Newton methods then update \(\mat{B}_k\) at each step such that the new approximation correctly predicts the change in gradients of the previous iteration, that is
\begin{equation}\label{eqn:long-secant-equation}
    \mat{B}_{k+1} (\vek{x}_{k+1} - \vek{x}_k) = \nabla f(\vek{x}_{k+1}) - \nabla f(\vek{x}_k).
\end{equation}
This is called the \emph{secant equation}.
We can write \cref{eqn:long-secant-equation} more succinctly if we define \(\vek{s}_k = \vek{x}_{k+1} - \vek{x}_k\) and let \(\widetilde{\mat{B}}_k = \int_0^1 \nabla^2 f(\vek{x}_k + t \vek{s}_k) \diff t\) be the Hessian of \(f\) averaged over all points on the line from \(\vek{x}_k\) to \(\vek{x}_{k+1}\). This gives the equivalent equation
\begin{equation}\label{eqn:short-secant-equation}
    \mat{B}_{k+1} \vek{s}_k = \widetilde{\mat{B}}_k \vek{s}_k.
\end{equation}

Most quasi-Newton methods then prescribe that among all possible choices of \(\mat{B}_{k+1}\) that are symmetric and satisfy the secant equation we take the one that is closest to \(\mat{B}_k\) in some norm.
The simplest such rule is called the Powell-symmetric-Broyden (PSB) update and is defined as
\begin{equation}\label{eqn:psb}
    \mat{B}_{k+1} = \argmin_{\mat{B} \in \R^{n \times n}_{\sym}} \norm{\mat{B} - \mat{B}_k}_F \text{ s.t. } \mat{B}\vek{s}_k = \widetilde{\mat{B}}_k \vek{s}_k. \tag{PSB}
\end{equation}
It is derived in \cite{powell_new_1970} from Broyden's method \cite{broyden_class_1965} by adding the symmetry constraint.

The Davidon-Fletcher-Powell (DFP) method \cite{davidon_variable_1959,fletcher_rapidly_1963}, even though it has been proposed before PSB, can be understood as a way to make the PSB method scale-invariant by choosing a weighted Frobenius norm.
Let \(\mat{W}_k = \widetilde{\mat{B}}_k^{-1/2}\) (or, in fact, any nonsingular matrix with \(\mat{W}_k^{-T} \mat{W}_k^{-1} \vek{s}_k = \widetilde{\mat{B}}_k \vek{s}_k\)) be the weight matrix, then the DFP update is given by
\begin{equation}\label{eqn:dfp}
    \mat{B}_{k+1} = \argmin_{\mat{B} \in \R^{n \times n}_{\sym}} \norm{\mat{W}_k^T \Paren{\mat{B} - \mat{B}_k} \mat{W}_k}_F \text{ s.t. } \mat{B}\vek{s}_k = \widetilde{\mat{B}}_k \vek{s}_k. \tag{DFP}
\end{equation}
If we rescale the input to \(f\) with the nonsingular matrix \(\mat{A}\), so that \(\bar{f}(\vek{x}) = f(\mat{A}\vek{x})\), and also rescale the iterates using \(\bar{\vek{x}}_k = \mat{A}^{-1} \vek{x}_k\), then the corresponding Hessian approximations of \(\bar{f}\) determined by the DFP method satisfy \(\bar{\mat{B}}_k = \mat{A}^T \mat{B}_k \mat{A}\) as long as it holds for the initial choice \(\bar{\mat{B}}_0\).

Finally, the famous BFGS method named after Broyden, Fletcher, Goldfarb and Shanno \cite{broyden_convergence_1970,fletcher_new_1970,goldfarb_family_1970,shanno_conditioning_1970}, is the dual of DFP in the sense that the new approximation minimizes the difference between inverse matrices in a weighted Frobenius norm:
\begin{equation}\label{eqn:bfgs}
    \mat{B}_{k+1} = \argmin_{\mat{B} \in \R^{n \times n}_{\sym}} \norm{\mat{W}_k^{-1} \Paren{\mat{B}^{-1} - \mat{B}_k^{-1}} \mat{W}_k^{-T}}_F \text{ s.t. } \mat{B} \vek{s}_k = \widetilde{\mat{B}}_k \vek{s}_k \tag{BFGS}
\end{equation}
The weight matrices \(\mat{W}_k\) are the same as above and in the same way they make the method scale invariant.
For more on these updating rules, consult the textbook by Dennis and Schnabel \cite[Chapter 9]{dennis_numerical_1996}

Using this characterization as least-change updates allows a straightforward generalization to tensors, except for the BFGS update.
Since tensors lack the concept of an inverse tensor, \cref{eqn:bfgs} cannot be used, and we will focus on \cref{eqn:psb,eqn:dfp}.
We now need to assume that \(f \colon \R^n \to \R\) is \(p\) times continuously differentiable and that as before a sequence of points \(\vek{x}_k\) with a corresponding sequence of steps \(\vek{s}_k = \vek{x}_{k+1} - \vek{x}_k\) is given.
We will denote the approximations to \(D^p f(\vek{x}_k)\) by \(\tns{C}_k \in \R^{\otimes^p n}_{\sym}\) and the true \(p\)th derivative averaged over all points on the line from \(\vek{x}_k\) to \(\vek{x}_{k+1}\) by
\begin{equation}\label{eqn:C-tilde-definition}
    \widetilde{\tns{C}}_k = \int_0^1 D^p f(\vek{x}_k + t \vek{s}_k) \diff t  \in \R^{\otimes^p n}_{\sym}.
\end{equation}
This means that, in particular, \(\widetilde{\tns{C}}_k[\vek{s}_k] = D^{p-1} f(\vek{x}_{k+1}) - D^{p-1} f(\vek{x}_k)\).
The \(p\)th-order analogue of the secant equation \cref{eqn:short-secant-equation} is given by
\begin{equation}\label{eqn:generalized-secant-equation}
    \tns{C}_{k+1} [\vek{s}_k] = D^{p-1} f(\vek{x}_{k+1}) - D^{p-1} f(\vek{x}_k) = \widetilde{\tns{C}}_k [\vek{s}_k]
\end{equation}
and the generalized update formula for \(\tns{C}_k\) reads
\begin{equation}\label{eqn:update-definition}
    \tns{C}_{k+1} = \argmin_{\tns{C} \in \R^{\otimes^p n}_{\sym}} \norm{(\tns{C} - \tns{C}_k)[\mat{W}_k]^p}_F \text{ s.t. } \tns{C}[\vek{s}_k] = \widetilde{\tns{C}}_k[\vek{s}_k]. \tag{HOSU}
\end{equation}
Even though this update is derived from the updates used in quasi-Newton methods, it is not itself associated with any optimization method and any optimization algorithm using approximate third (or higher) derivatives is also clearly different from Newton's method.
To highlight this distinction we will call the update rule the \emph{higher-order secant update} \cref{eqn:update-definition} because of its connection with the (generalized) secant equation \cref{eqn:generalized-secant-equation}.
It provides a sequence of approximations of the \(p\)th derivative of \(f\) based solely on evaluations of the \((p-1)\)st derivative at the iterates \(\vek{x}_k\), given some initial approximation \(\tns{C}_0\).

Note that unlike for the DFP update we will not assume any specific choice of weight matrices, but rather consider them to be a given sequence of nonsingular matrices.
In particular, that covers the higher-order PSB (\(\mat{W}_k = \mat{I}\)) and DFP (\(\mat{W}_k^{-T} \mat{W}_k^{-1} \vek{s}_k = \nabla f(\vek{x}_{k+1}) - \nabla f(\vek{x}_k)\)) updates, which simplify to \cref{eqn:psb,eqn:dfp} for \(p = 2\).

\section{Characterization of higher-order secant updates}\label{sec:characterization}

The following theorem provides a full characterization of one step of the update in \cref{eqn:update-definition}.
Note that despite the intentional notational similarity the quantities in the statement are independent of the definitions in the previous section.

\begin{theorem}\label{thm:characterization}
    Let \(\tns{C}_{\bullet} \in \R^{\otimes^p n}_{\sym}\) (the current approximation), \(\widetilde{\tns{C}} \in \R^{\otimes^p n}_{\sym}\) (the integrated true derivative), a nonsingular matrix \(\mat{W} \in \R^{n \times n}\) (the weight matrix) and a nonzero \(\vek{s} \in \R^n\) (the step) be given.
    The following equations all have the same unique solution \(\tns{C}_+ \in \R^{\otimes^p n}_{\sym}\) (the new approximation):
    \begin{enumerate}[(a)]
        \item\label{item:characterization-min-frobenius} \(\tns{C}_+ = \argmin_{\tns{C} \in \R^{\otimes^p n}_{\sym}} \norm{(\tns{C} - \tns{C}_{\bullet})[\mat{W}]^p}_F \text{ s.t. } \tns{C}[\vek{s}] = \widetilde{\tns{C}}[\vek{s}]\)
        \item\label{item:characterization-explicit} \(\tns{C}_+ = \tns{C}_{\bullet} + \sum_{j=1}^p (-1)^{j+1} \binom{p}{j} (\vek{v}^T \vek{s})^{-j} P_{\sym} \Paren{\Paren{\otimes^j \vek{v}} \otimes (\widetilde{\tns{C}} - \tns{C}_{\bullet})[\vek{s}]^{j}}\)
        \item\label{item:characterization-low-rank} \(\tns{C}_+ = \tns{C}_{\bullet} + P_{\sym}(\tns{A} \otimes \vek{v})\) for the unique \((p-1)\)-tensor \(\tns{A} \in \R^{\otimes^{p-1} n}_{\sym}\) which satisfies \(P_{\sym}(\tns{A} \otimes \vek{v})[\vek{s}] = (\widetilde{\tns{C}} - \tns{C}_{\bullet})[\vek{s}]\)
        \item\label{item:characterization-recursive} \((\tns{C}_+ - \widetilde{\tns{C}})[\mat{W}]^p = (\tns{C}_{\bullet} - \widetilde{\tns{C}})[\mat{W}]^p \Brack{\mat{I} - \frac{\mat{W}^{-1}\vek{s}\vek{s}^T \mat{W}^{-T}}{\vek{s}^T \mat{W}^{-T} \mat{W}^{-1} \vek{s}}}^p\)
    \end{enumerate}
    where \(\vek{v} = \mat{W}^{-T}\mat{W}^{-1}\vek{s}\).
\end{theorem}

\begin{proof}
    We will first prove the result for \(\mat{W} = \mat{I}\) and then see how that implies the full result for any nonsingular weight matrix \(\mat{W}\).

    The first characterization can be rewritten as
    \begin{equation}\label{eqn:general-psb-update}
        \tns{C}_+ = \tns{C}_{\bullet} + \argmin_{\tns{U} \in \R^{\otimes^p n}_{\sym}} \norm{\tns{U}}_F \text{ s.t. } \tns{U}[\vek{s}] = (\widetilde{\tns{C}} - \tns{C}_{\bullet})[\vek{s}].
    \end{equation}
    Let \(\mat{Q} \in \R^{n \times n}\) be an orthogonal matrix that maps \(\vek{e}_1\) to some scalar multiple of \(\vek{s}\).
    This means \(\tns{U}[\mat{Q}]^p\) has the same Frobenius norm as \(\tns{U}\) and \(\tns{U}[\mat{Q}]^p[\vek{e}_{i_1}, \dots, \vek{e}_{i_p}]\) is fully determined by the equality constraint in \cref{eqn:general-psb-update} if \(1 \in \{i_1, \dots, i_p\}\) because of symmetry.
    On the other hand, if \(1 \notin \{i_1, \dots, i_p\}\) there are no constraints on the values of \(\tns{U}[\mat{Q}]^p[\vek{e}_{i_1}, \dots, \vek{e}_{i_p}]\), so the unique choice that minimizes \(\norm{\tns{U}[\mat{Q}]^p}_F = \norm{\tns{U}}_F\) is clearly to set all of these values to zero.
    Therefore, there is a unique solution to the minimization problem in \labelcref{item:characterization-min-frobenius}, and it is fully characterized by the fact that the update tensor \(\tns{U} = \tns{C}_+ - \tns{C}_{\bullet}\) is symmetric and satisfies the following two properties:
    \begin{subequations}
        \begin{align}
            \tns{U}[\vek{s}]                     & = \paren{\widetilde{\tns{C}} - \tns{C}_{\bullet}}[\vek{s}] \label{eqn:update-secant-equation}        \\
            \tns{U}[\vek{u}_1, \dots, \vek{u}_p] & = 0 \text{ if all } \vek{u}_i \text{ are orthogonal to } \vek{s} \label{eqn:update-minimum-property}
        \end{align}
    \end{subequations}
    We will use this as the basis to show the equivalence with all other characterizations.

    In \labelcref{item:characterization-explicit} we claim that the update \(\tns{U}\) has the form
    \begin{equation}
        \sum_{j=1}^p (-1)^{j+1} \tbinom{p}{j} \norm{\vek{s}}_2^{-2j} P_{\sym} \Paren{\Paren{\otimes^j \vek{s}} \otimes \paren{\widetilde{\tns{C}} - \tns{C}_{\bullet}}[\vek{s}]^{j}}
    \end{equation}
    This tensor is clearly symmetric since it is a sum of symmetric tensors.
    To show that it satisfies property \cref{eqn:update-secant-equation} consider
    \begin{subequations}\label{eqn:update-explicit-cancellation}
        \begin{align}
             & \phantom{={}}\sum_{j=1}^p (-1)^{j+1} \tbinom{p}{j} \norm{\vek{s}}_2^{-2j} P_{\sym} \Paren{\Paren{\otimes^j \vek{s}} \otimes \paren{\widetilde{\tns{C}} - \tns{C}_{\bullet}}[\vek{s}]^{j}}[\vek{s}]                                                                                                                                                \\
            \begin{split}
                {} = \sum_{j=1}^p (-1)^{j+1} \norm{\vek{s}}_2^{-2j} \bigg( \tbinom{p-1}{j-1} P_{\sym} \Paren{\Paren{\otimes^{j-1} \vek{s}} \otimes \paren{\widetilde{\tns{C}} - \tns{C}_{\bullet}}[\vek{s}]^{j}} \norm{\vek{s}}_2^2 \phantom{\bigg)}\\
                {} + \tbinom{p-1}{j} P_{\sym} \Paren{\Paren{\otimes^j \vek{s}} \otimes \paren{\widetilde{\tns{C}} - \tns{C}_{\bullet}}[\vek{s}]^{j+1}} \bigg) \phantom{\norm{\vek{s}_2^2}}
            \end{split} \\
             & = P_{\sym} \Paren{ \paren{\widetilde{\tns{C}} - \tns{C}_{\bullet}}[\vek{s}] } = \paren{\widetilde{\tns{C}} - \tns{C}_{\bullet}}[\vek{s}].
        \end{align}
    \end{subequations}
    The first equality uses the fact that each summand is the symmetric projection of an outer product of two symmetric tensors, a \(j\)-tensor and a \((p-j)\)-tensor.
    Therefore, there are \(\binom{p}{j}\) distinct ways to orient this outer product and for \(\binom{p-1}{j-1}\) of them the vector \(\vek{s}\) is applied to the \(j\)-tensor and for \(\binom{p-1}{j}\) to the \((p-j)\)-tensor.
    A close examination of the expression on the second line shows that it is a telescoping sum where all terms except the ones with coefficients \(\binom{p-1}{0}\) and \(\binom{p-1}{p}\) cancel out.
    Because \(\binom{p-1}{p} = 0\) the only remaining term is the symmetric projection of \((\widetilde{\tns{C}} - \tns{C}_{\bullet})[\vek{s}]\).
    Property \cref{eqn:update-minimum-property} follows from a similar consideration to the one above.
    If \(\vek{u}_1, \dots, \vek{u}_p\) are all orthogonal to \(\vek{s}\), then
    \begin{equation}
        P_{\sym} \Paren{\Paren{\otimes^j \vek{s}} \otimes \paren{\widetilde{\tns{C}} - \tns{C}_{\bullet}}[\vek{s}]^{j}}[\vek{u}_1, \dots, \vek{u}_p] = 0
    \end{equation}
    for \(j \geq 1\) because no matter how the outer product is oriented there is always a factor of \(\vek{s}^T \vek{u}_i = 0\) in the result.

    For \labelcref{item:characterization-low-rank} we need to show that we can always write the update in the form \(P_{\sym}(\tns{A} \otimes \vek{s})\) for some symmetric \((p-1)\)-tensor \(\tns{A}\) and that \(\tns{A}\) is unique such that the corresponding update satisfies the secant equation \cref{eqn:update-secant-equation}.
    Note that the expression for the update in \labelcref{item:characterization-explicit} is already of the form \(P_{\sym}(\tns{A} \otimes \vek{s})\):
    \begin{subequations}
        \begin{align}
             & \phantom{={}}\sum_{j=1}^p (-1)^{j+1} \tbinom{p}{j} \norm{\vek{s}}_2^{-2j} P_{\sym} \Paren{\Paren{\otimes^j \vek{s}} \otimes \paren{\widetilde{\tns{C}} - \tns{C}_{\bullet}}[\vek{s}]^{j}}                               \\
             & = P_{\sym} \Paren{ P_{\sym} \Paren{ \sum_{j=1}^p (-1)^{j+1} \tbinom{p}{j} \norm{\vek{s}}_2^{-2j} \Paren{\otimes^{j-1} \vek{s}} \otimes \paren{\widetilde{\tns{C}} - \tns{C}_{\bullet}}[\vek{s}]^{j} } \otimes \vek{s} }
        \end{align}
    \end{subequations}
    It remains to show that among all updates of the form \(P_{\sym}(\tns{A} \otimes \vek{s})\) there is only one choice of \(\tns{A} \in \R^{\otimes^p n}_{\sym}\) such that \cref{eqn:update-secant-equation} holds.
    Consider the linear map
    \begin{align*}
        \phi \colon \R^{\otimes^{p-1} n}_{\sym} & \to \R^{\otimes^{p-1} n}_{\sym}                    \\
        \tns{A}                                 & \mapsto P_{\sym}(\tns{A} \otimes \vek{s})[\vek{s}]
    \end{align*}
    which maps a finite-dimensional vector space to itself.
    Combining what we already showed for \labelcref{item:characterization-explicit} with the observation that the update in \labelcref{item:characterization-explicit} is of the desired form, this map is surjective (we can prescribe any \((\widetilde{\tns{C}} - \tns{C}_{\bullet})[\vek{s}]\)) and so it must be bijective, which shows uniqueness of \(\tns{A}\).

    Lastly, \labelcref{item:characterization-recursive} claims that the update can be written as
    \begin{equation}
        \paren{\widetilde{\tns{C}} - \tns{C}_{\bullet}} - \paren{\widetilde{\tns{C}} - \tns{C}_{\bullet}} \Brack{\mat{I} - \frac{\vek{s} \vek{s}^T}{\norm{\vek{s}}_2^2}}^p.
    \end{equation}
    Clearly, property \cref{eqn:update-secant-equation} is satisfied because applying this tensor to \(\vek{s}\) makes the second term vanish, leaving only the desired result.
    Moreover, for any vector \(\vek{u}\) that is orthogonal to \(\vek{s}\) the matrix \(\mat{I} - \frac{\vek{s} \vek{s}^T}{\norm{\vek{s}}_2^2}\) maps \(\vek{u}\) to itself.
    This means applying the tensor above to \(\vek{u}_1, \dots, \vek{u}_p\), all of which are orthogonal to \(\vek{s}\), will give zero because both terms cancel out.
    This is property \cref{eqn:update-minimum-property}.

    Now that the equivalence has been established for \(\mat{W} = \mat{I}\), we consider the general case where \(\mat{W} \in \R^{n \times n}\) is any nonsingular matrix.
    The minimization in \labelcref{item:characterization-min-frobenius}
    \begin{equation}
        \tns{C}_+ = \argmin_{\tns{C} \in \R^{\otimes^p n}_{\sym}} \norm{(\tns{C} - \tns{C}_{\bullet})[\mat{W}]^p}_F \text{ s.t. } \tns{C}[\vek{s}] = \widetilde{\tns{C}}[\vek{s}]
    \end{equation}
    can be rewritten using \(\tns{D}_{\bullet} = \tns{C}_{\bullet}[\mat{W}]^p\), \(\widetilde{\tns{D}} = \widetilde{\tns{C}}[\mat{W}]^p\), \(\tns{D}_+ = \tns{C}_+[\mat{W}]^p\) and \(\vek{r} = \mat{W}^{-1}\vek{s}\) as
    \begin{equation}
        \tns{D}_+ = \argmin_{\tns{D} \in \R^{\otimes^p n}_{\sym}} \norm{\tns{D} - \tns{D}_{\bullet}}_F \text{ s.t. } \tns{D}[\vek{r}] = \widetilde{\tns{D}}[\vek{r}].
    \end{equation}
    Applying the existing characterizations and the fact that \(\tns{C}_+ = \tns{D}_+[\mat{W}^{-1}]^p\) we get the claim after some algebraic manipulations.
\end{proof}

The different characterizations listed in the theorem highlight different aspects of the update: \labelcref{item:characterization-min-frobenius} is the least-change update characterization that we motivated from quasi-Newton updates, \labelcref{item:characterization-explicit} gives an explicit formula for the computation, \labelcref{item:characterization-low-rank} shows that the update has a low-rank structure (as discussed below) and \labelcref{item:characterization-recursive} gives a recursive relationship between \(\tns{C}_+\) and \(\tns{C}_\bullet\) that we will make use of in the next section.

An important observation about the update is that the only dependence on the weight matrix \(\mat{W}\) comes in the form of a dependence on \(\vek{v}\) and is moreover invariant under rescaling of \(\vek{v}\).
Most of the degrees of freedom in choosing \(\mat{W}\) are therefore irrelevant.
The only restriction on \(\vek{v}\) comes from the fact that \(\mat{W}^{-T}\mat{W}^{-1}\) is positive definite and so only when \(\vek{v}^T \vek{s} > 0\) there is a matrix \(\mat{W}\) that makes the characterization theorem true.
For the explicit update \labelcref{item:characterization-explicit} to be well-defined we only need \(\vek{v}^T \vek{s} \neq 0\) though, since the update is the same when the sign of \(\vek{v}\) is swapped.

The characterization \labelcref{item:characterization-low-rank} in the above theorem shows that the update exhibits a certain low-rank structure in that it can be expressed as the symmetric projection of a tensor with multilinear rank at most \((n, \dots, n, 1)\).
It might seem like this is suboptimal, and we should aim for an update rule that produces updates with small CP rank.
However, this is impossible for \(p > 2\).
By construction any update \(\tns{U}\) must satisfy the secant equation \(\tns{U}[\vek{s}] = (\widetilde{\tns{C}}-\tns{C}_{\bullet})[\vek{s}]\).
Without assuming any structure of the function \(f\) the right-hand side of the secant equation can be any element of \(\R^{\otimes^{p-1} n}_{\sym}\), which is a space of dimension \(O(n^{p-1})\).
At the same time the space of tensors of rank at most \(r\) has less than \(rpn\) degrees of freedom.
The exponents of \(n\) in these two expressions only match up when \(p=2\).
For \(p=2\) we indeed get that by \labelcref{item:characterization-low-rank} our updates (which include PSB and DFP) can always be expressed as rank-two matrices.\footnote{For \(p=2\) the tensor \(\tns{A}\) is actually a vector so that \(\tns{A} \otimes \vek{v}\) is a rank-one matrix and its symmetric projection has at most rank two.}
For \(p>2\) the low-rank result we have is essentially optimal: If we fix \(\vek{v}\), then the space of tensors of the form \(P_{\sym}(\tns{A} \otimes \vek{v})\) has exactly the same dimension as \(\R^{\otimes^{p-1} n}_{\sym}\) as shown by the uniqueness of \(\tns{A}\).

These considerations imply that a generalization of the symmetric rank-one update (SR1) that stays true to its name is impossible.
However, we can try to choose \(\vek{v}\) in a way that resembles the approach taken by SR1.
For matrices the SR1 update fits our general update formula in \labelcref{item:characterization-explicit} by using \(\vek{v} = (\widetilde{\mat{B}} - \mat{B}_{\bullet})\vek{s}\).
This means \(\vek{v}\) is chosen such that it aligns with the difference between the actual and predicted change in gradients.
In that spirit we could choose \(\vek{v}\) such that \(\otimes^{p-1} \vek{v}\) is aligned as far as possible with the difference between the actual and predicted change in derivatives, i.e.
\begin{equation}
    \vek{v}^* = \argmin_{\norm{\vek{v}}_2 = 1} \abs{\innerprod{(\widetilde{\tns{C}} - \tns{C}_{\bullet})[\vek{s}]}{\otimes^{p-1} \vek{v}}_F}.
\end{equation}
This minimization problem is equivalent to finding the best rank-one approximation to a \((p-1)\)-tensor and so it is NP-hard for \(p>3\) \cite{hillar_most_2013}, but still tractable for \(p=3\) in which case efficient algorithms exist for approximating the eigenvector corresponding to the largest absolute eigenvalue of a symmetric matrix.
By choosing \(\vek{v}\) in this manner we might hope to achieve similar numerical properties to the ones mentioned in \cite{conn_convergence_1991}, where the authors found that SR1 matrices produce significantly better derivative approximations than DFP or even BFGS matrices.

\begin{example}\label{ex:simple-update}
    We discussed some properties that can be deduced from \labelcref{item:characterization-explicit,item:characterization-low-rank} in the previous paragraphs, but their general structure might still seem complicated.
    To show their inner workings we consider a simple example update for the case \(p=2\) and \(p=3\).
    In both cases \(\vek{v}\) and \(\vek{s}\) are the first unit vector \(\vek{e}_1\), which helps to highlight the structure of the outer products involved.

    Consider the matrix case (\(p=2\)) with
    \begin{align}
        \tns{C}_{\bullet} = \begin{pmatrix}
                                0 & 0 & 0 \\
                                0 & 0 & 0 \\
                                0 & 0 & 0
                            \end{pmatrix} \quad \text{and} \quad \widetilde{\tns{C}} = \begin{pmatrix}
                                                                                           1 & 1 & 1 \\
                                                                                           1 & 1 & 1 \\
                                                                                           1 & 1 & 1
                                                                                       \end{pmatrix}
    \end{align}
    first.
    We have
    \begin{subequations}\label{eqn:example-explicit-formula-p2}
        \begin{align}
            \tns{C}_+ & = \tns{C}_{\bullet} + \sum_{j=1}^p (-1)^{j+1} \binom{p}{j} (\vek{v}^T \vek{s})^{-j} P_{\sym} \Paren{\Paren{\otimes^j \vek{v}} \otimes (\widetilde{\tns{C}} - \tns{C}_{\bullet})[\vek{s}]^{j}} \\
                      & = 2 P_{\sym} \begin{pmatrix}
                                         1 & 1 & 1 \\
                                         0 & 0 & 0 \\
                                         0 & 0 & 0
                                     \end{pmatrix} - 1 P_{\sym} \begin{pmatrix}
                                                                    1 & 0 & 0 \\
                                                                    0 & 0 & 0 \\
                                                                    0 & 0 & 0
                                                                \end{pmatrix}                                                                                                                                        \\
                      & = \begin{pmatrix}
                              2 & 1 & 1 \\
                              1 & 0 & 0 \\
                              1 & 0 & 0
                          \end{pmatrix} - \begin{pmatrix}
                                              1 & 0 & 0 \\
                                              0 & 0 & 0 \\
                                              0 & 0 & 0
                                          \end{pmatrix} = \begin{pmatrix}
                                                              1 & 1 & 1 \\
                                                              1 & 0 & 0 \\
                                                              1 & 0 & 0
                                                          \end{pmatrix}
        \end{align}
    \end{subequations}
    which is the smallest symmetric matrix (measured in the Frobenius norm) whose first column is all ones.
    It fits into the low-rank form of \labelcref{item:characterization-low-rank} since
    \begin{equation}
        \tns{C}_+ = P_{\sym} \begin{pmatrix}
            1 & 2 & 2 \\
            0 & 0 & 0 \\
            0 & 0 & 0 \\
        \end{pmatrix} = \tns{C}_{\bullet} + P_{\sym} \Paren{ \begin{pmatrix}
                1 \\
                2 \\
                2
            \end{pmatrix} \otimes \vek{v} }.
    \end{equation}
    The recursive formula \labelcref{item:characterization-recursive} simplifies to
    \begin{equation}\label{eqn:example-simplification-recursive}
        (\tns{C}_+ - \widetilde{\tns{C}}) = (\tns{C}_{\bullet} - \widetilde{\tns{C}})\Brack{ \mat{I} - \frac{\vek{s} \vek{v}^T}{\vek{s}^T \vek{v}} }^p = (\tns{C}_{\bullet} - \widetilde{\tns{C}})\Brack{ \mat{I} - \frac{\vek{e}_1 \vek{e}_1^T}{\vek{e}_1^T \vek{e}_1} }^p
    \end{equation}
    using the definition of \(\vek{v}\) and the values of \(\vek{v}\) and \(\vek{s}\) in this example.
    The matrix \(\mat{I} - \vek{e}_1 \vek{e}_1^T\) is an orthogonal projection onto the subspace orthogonal to \(\vek{e}_1\).
    This means the error in \(\tns{C}_+\) compared to \(\widetilde{\tns{C}}\) is zero in the first row and column and the same as before anywhere else.
    This is exactly the shape we see as the final result in \cref{eqn:example-explicit-formula-p2}.

    Now consider the same example in the tensor case (\(p=3\)), although now in \(\R^2\) to save space.
    The current third derivative approximation and the integrated true third derivative are given by
    \begin{align}
        \tns{C}_{\bullet} = \begin{pmatrix}
                                \begin{pmatrix}
                0 & 0 \\
                0 & 0
            \end{pmatrix} \;
                                \begin{pmatrix}
                0 & 0 \\
                0 & 0
            \end{pmatrix}
                            \end{pmatrix} \quad \text{and} \quad \widetilde{\tns{C}} = \begin{pmatrix}
                                                                                           \begin{pmatrix}
                1 & 1 \\
                1 & 1
            \end{pmatrix} \;
                                                                                           \begin{pmatrix}
                1 & 1 \\
                1 & 1
            \end{pmatrix}
                                                                                       \end{pmatrix}.
    \end{align}
    In this notation the tensor is split into two matrix slices along the first mode.
    If the entries of a tensor \(\tns{T}\) are \(t_{ijk}\) then the first matrix contains the entries of the form \(t_{1jk}\) and the second matrix the entries \(t_{2jk}\).
    In this case, we have
    \begin{subequations}\label{eqn:example-explicit-formula-p3}
        \begin{align}
            \tns{C}_+ & = \tns{C}_{\bullet} + \sum_{j=1}^p (-1)^{j+1} \binom{p}{j} (\vek{v}^T \vek{s})^{-j} P_{\sym} \Paren{\Paren{\otimes^j \vek{v}} \otimes (\widetilde{\tns{C}} - \tns{C}_{\bullet})[\vek{s}]^{j}} \\
                      & =  3 P_{\sym}
            \begin{pmatrix}
                \begin{pmatrix}
                    1 & 1 \\
                    1 & 1
                \end{pmatrix} \;
                \begin{pmatrix}
                    0 & 0 \\
                    0 & 0
                \end{pmatrix}
            \end{pmatrix}
            - 3 P_{\sym}
            \begin{pmatrix}
                \begin{pmatrix}
                    1 & 1 \\
                    0 & 0
                \end{pmatrix} \;
                \begin{pmatrix}
                    0 & 0 \\
                    0 & 0
                \end{pmatrix}
            \end{pmatrix}                                                                                                                                                                              \\ \notag
                      & \phantom{{}=} +1 P_{\sym}
            \begin{pmatrix}
                \begin{pmatrix}
                    1 & 0 \\
                    0 & 0
                \end{pmatrix} \;
                \begin{pmatrix}
                    0 & 0 \\
                    0 & 0
                \end{pmatrix}
            \end{pmatrix}                                                                                                                                                                              \\
                      & =  \begin{pmatrix}
                               \begin{pmatrix}
                    3 & 2 \\
                    2 & 1
                \end{pmatrix} \;
                               \begin{pmatrix}
                    2 & 1 \\
                    1 & 1
                \end{pmatrix}
                           \end{pmatrix}
            -
            \begin{pmatrix}
                \begin{pmatrix}
                    3 & 1 \\
                    1 & 0
                \end{pmatrix} \;
                \begin{pmatrix}
                    1 & 0 \\
                    0 & 0
                \end{pmatrix}
            \end{pmatrix}
            +
            \begin{pmatrix}
                \begin{pmatrix}
                    1 & 0 \\
                    0 & 0
                \end{pmatrix} \;
                \begin{pmatrix}
                    0 & 0 \\
                    0 & 0
                \end{pmatrix}
            \end{pmatrix}                                                                                                                                                                              \\
                      & = \begin{pmatrix}
                              \begin{pmatrix}
                    1 & 1 \\
                    1 & 1
                \end{pmatrix} \;
                              \begin{pmatrix}
                    1 & 1 \\
                    1 & 0
                \end{pmatrix}
                          \end{pmatrix}
        \end{align}
    \end{subequations}
    which is the smallest symmetric 3-tensor (measured in the Frobenius norm) whose first slice is a matrix of all ones.
    The cancellation between the different terms above is exactly the one described in the proof in \cref{eqn:update-explicit-cancellation}.
    It is possible to express the update as the symmetric projection of an outer product as follows:
    \begin{equation}
        \tns{C}_+ = P_{\sym} \begin{pmatrix}
            \begin{pmatrix}
                1   & 3/2 \\
                3/2 & 3
            \end{pmatrix} \; \begin{pmatrix}
                                 0 & 0 \\
                                 0 & 0
                             \end{pmatrix}
        \end{pmatrix} = \tns{C}_{\bullet} + P_{\sym} \Paren{ \begin{pmatrix}
                1   & 3/2 \\
                3/2 & 3
            \end{pmatrix} \otimes \vek{v} }.
    \end{equation}
    \Cref{eqn:example-simplification-recursive} holds as before, so \labelcref{item:characterization-recursive} tells us that the new approximation coincides with \(\widetilde{\tns{C}}\) in the first slice along each mode (entries \(t_{1jk}\), \(t_{i1k}\) and \(t_{ij1}\)) and with \(\tns{C}_{\bullet}\) in the remaining entries.
    Again, this is exactly what we see as the final result of \cref{eqn:example-explicit-formula-p3}.
\end{example}

\section{Convergence of approximate derivates}\label{sec:convergence}

Now that we know how the updates look like, we want to show that the approximate \(p\)th derivatives converge to the true derivative under certain assumptions.
For there to be a true derivative, we will assume in this section that \(f\) is \(p\) times continuously differentiable.
The main tool of this section is the characterization \cref{thm:characterization} \labelcref{item:characterization-recursive} which, when applied to \cref{eqn:update-definition}, becomes
\begin{equation}\label{eqn:update-recursive}
    \paren{\tns{C}_{k+1} - \widetilde{\tns{C}}_k}[\mat{W_k}]^p = (\tns{C}_k - \widetilde{\tns{C}}_k)[\mat{W}_k]^p [\mat{P}_k]^p \text{ where } \mat{P}_k = \mat{I} - \frac{\mat{W}_k^{-1}\vek{s}_k\vek{s}_k^T \mat{W}_k^{-T}}{\vek{s}_k^T \mat{W}_k^{-T} \mat{W}_k^{-1} \vek{s}_k}.
\end{equation}
Note that \(\mat{P}_k\) is the orthogonal projection onto the orthogonal complement of \(\mat{W}_k^{-1} \vek{s}_k\).

\subsection{\texorpdfstring{Convergence for \(p\)th-order polynomials}{Convergence for pth-order polynomials}}

To show the usefulness of \cref{eqn:update-recursive} we consider the case when \(\mat{W}_k\) and \(\widetilde{\tns{C}}_k\) are constant first. If we additionally assume that the scaled steps \(\mat{W}_k^{-1}\vek{s}_k\) are orthogonal, convergence is quite straightforward to prove.

\begin{theorem}
    Let \(\tns{C}_0 \in \R^{\otimes^p n}_{\sym}\) be given and update the approximations \(\tns{C}_k\) according to \cref{eqn:update-definition}.
    Assume \(D^p f(\vek{x}) = \tns{C}_*\) everywhere (which makes \(f\) a \(p\)th-order polynomial) and \(\mat{W}_k = \mat{W}_*\) for every \(k \in \N\).
    After \(n\) steps \(\vek{s}_k\) such that \(\mat{W}_*^{-1}\vek{s}_k\) are orthogonal, we have \(\tns{C}_n = \tns{C}_*\).
\end{theorem}
\begin{proof}
    Under the assumptions of the theorem, repeated application of \cref{eqn:update-recursive} gives
    \begin{equation}
        \Paren{\tns{C}_{n} - \tns{C}_*}[\mat{W}_*]^p = \Paren{\tns{C}_0 - \tns{C}_*}[\mat{W}_*]^p \Brack{\prod_{k=0}^{n-1} \mat{P}_k}^p.
    \end{equation}
    Because \(\mat{W}_0^{-1}\vek{s}_0, \dots, \mat{W}_{n-1}^{-1} \vek{s}_{n-1}\) are orthogonal and \(\mat{P}_k\) are orthogonal projections on their orthogonal complements, \(\prod_{k=0}^{n-1} \mat{P}_k = \mat{0}\), so that \((\tns{C}_n - \tns{C}_*)[\mat{W}_*]^p = \tns{0}\).
    This implies \(\tns{C}_n = \tns{C}_*\) because \(\mat{W}_*\) is nonsingular.
\end{proof}

It is clear that orthogonality of the vectors \(\mat{W}_*^{-1}\vek{s}_k\) is quite a strong assumption.
For the \(p=2\) case, the SR1 update satisfies an equivalent statement with the weaker assumption that the \(n\) steps \(\vek{s}_k\) are linearly independent \cite[Theorem 3.2.1]{fletcher_practical_2000}.
Note however, that this is achieved by using a weight matrix \(\mat{W}_*\) with the property \(\mat{W}_*^{-T} \mat{W}_*^{-1} \vek{s}_k = (\mat{B}_k - \widetilde{\mat{B}}_k)\vek{s}_k\).
The convergence proof then boils down to showing that \(\mat{W}_*^{-1} \vek{s}_k\) are orthogonal.
Similarly, Theorem 3.4.1 in \cite{fletcher_practical_2000} proves convergence of the Broyden class update formulas under the alternative assumption of using exact line searches.
This exact line search condition is used to show that the search directions are conjugate with respect to the constant positive definite Hessian \(\mat{B}_*\).
In the case of DFP the constant weight matrix is given by \(\mat{W}_k = \widetilde{\mat{B}}_k^{-1/2} = \mat{B}_*^{-1/2}\) so that conjugacy of the search directions is equivalent to orthogonality of \(\mat{W}_*^{-1} \vek{s}_k\).
In this context, the result above is essentially optimal without any other assumption on the choice of weight matrix \(\mat{W}_*\) or the choice of steps \(\vek{s}_k\).

\subsection{Bounded deterioration}

For convergence results for general functions \(f\) we first need to establish two lemmas that will help us when \(\vek{x}_k \to \vek{x}_*\) and \(\mat{W}_k \to \mat{W}_*\).
The first one shows that \(\widetilde{\tns{C}}_k\) converges to the \(p\)th derivative at \(\vek{x}_*\), which we denote by \(\tns{C}_* = D^p f(\vek{x}_*)\), and the second one gives a bound on the error term that we incur if we replace \(\widetilde{\tns{C}}_k\) by \(\tns{C}_*\) in \cref{eqn:update-recursive}.
Combining the two gives what Dennis and Schnabel \cite{dennis_numerical_1996} call the \emph{bounded deterioration principle} in that \(\tns{C}_{k+1}\) can only be slightly worse than \(\tns{C}_k\) at approximating \(\tns{C}_*\) as long as \(\vek{x}_k\) and \(\vek{x}_{k+1}\) are close enough to \(\vek{x}_*\).

\begin{lemma}\label{thm:C-tilde-convergence}
    Let \(\vek{x}_* \in \R^n\) and \(\tns{C}_* = D^p f(\vek{x}_*)\).
    \begin{enumerate}[(a)]
        \item For \(\vek{x}_k \to \vek{x}_*\) we have \(\widetilde{\tns{C}}_k \to \tns{C}_*\) as \(k \to \infty\) and
        \item if \(D^p f\) is Lipschitz continuous with constant \(L\), then
              \begin{equation*}
                  \norm{\widetilde{\tns{C}}_k - \tns{C}_*}_2 \leq \frac{L}{2} \Paren{ \norm{\vek{x}_k - \vek{x}_*}_2 + \norm{\vek{x}_{k+1} - \vek{x}_*}_2 }
              \end{equation*}
              for all \(k \in \N\).
    \end{enumerate}
\end{lemma}
\begin{proof}
    By definition in \cref{eqn:C-tilde-definition} we have
    \begin{subequations}
        \begin{align}
            \norm{\widetilde{\tns{C}}_k - \tns{C}_*}_2 & = \norm{\int_0^1 D^p f(\vek{x}_k + t \vek{s}_k) \diff t - D^p f(\vek{x}_*)}_2                                        \\
                                                       & \leq \int_0^1 \norm{D^p f(\vek{x}_k + t \vek{s}_k) - D^p f(\vek{x}_*)}_2 \diff t. \label{eqn:C-tilde-bound-integral}
        \end{align}
    \end{subequations}
    Since \(D^p f\) is continuous the right-hand side will become arbitrarily small as \(\vek{x}_k\) and \(\vek{x}_{k+1}\) converge to \(\vek{x}_*\).
    This shows \(\widetilde{\tns{C}}_k \to \tns{C}_*\).

    If we assume Lipschitz continuity of \(D^p f\), we can bound the integrand in \cref{eqn:C-tilde-bound-integral}:
    \begin{subequations}\label{eqn:C-tilde-Lipschitz-bound}
        \begin{align}
            \norm{\widetilde{\tns{C}}_k - \tns{C}_*}_2 & \leq \int_0^1 L \norm{\vek{x}_k + t \vek{s}_k - \vek{x}_*}_2 \diff t                                  \\
                                                       & \leq L  \int_0^1 (1 - t)\norm{\vek{x}_k - \vek{x}_*}_2 + t \norm{\vek{x}_{k+1} - \vek{x}_*}_2 \diff t \\
                                                       & = \frac{L}{2} \Paren{\norm{\vek{x}_k - \vek{x}_*}_2 + \norm{\vek{x}_{k+1} - \vek{x}_*}_2}
        \end{align}
    \end{subequations}
    This gives the second claim.
\end{proof}

\begin{lemma}\label{thm:error-tensor-bound}
    If we define the error tensor \(\tns{E}_k\) by
    \begin{equation}\label{eqn:error-tensor-recurrence}
        \Paren{\tns{C}_{k+1} - \tns{C}_*}[\mat{W}_k]^p = \Paren{\tns{C}_k - \tns{C}_*}[\mat{W}_k]^p [\mat{P}_k]^p + \tns{E}_k[\mat{W}_k]^p
    \end{equation}
    then \(\norm{\tns{E}_k}_2 \leq (1 + \kappa_2(\mat{W}_k)^p) \norm{\widetilde{\tns{C}}_k - \tns{C}_*}_2\).
\end{lemma}
\begin{proof}
    Subtracting \cref{eqn:update-recursive} from \cref{eqn:error-tensor-recurrence} gives
    \begin{equation}\label{eqn:error-tensor-C-tilde-recurrence}
        \paren{\widetilde{\tns{C}}_k - \tns{C}_*}[\mat{W}_k]^p = \paren{\widetilde{\tns{C}}_k - \tns{C}_*}[\mat{W}_k]^p \Brack{\mat{P}_k}^p + \tns{E}_k[\mat{W}_k]^p.
    \end{equation}
    Multiply both sides by \(\mat{W}_k^{-1}\) from all sides and rearrange to find
    \begin{subequations}
        \begin{align}
            \norm{\tns{E}_k}_2 & = \norm{\paren{\widetilde{\tns{C}}_k - \tns{C}_*} - \paren{\widetilde{\tns{C}}_k - \tns{C}_*}[\mat{W}_k]^p \Brack{\mat{P}_k}^p[\mat{W}_k^{-1}]^p}_2 \\
                               & \leq \Paren{1 + \norm{\mat{W}_k \mat{P}_k \mat{W}_k^{-1}}_2^p} \norm{\widetilde{\tns{C}}_k - \tns{C}_*}_2                                           \\
                               & \leq \Paren{1 + \kappa_2(\mat{W}_k)^p} \norm{\widetilde{\tns{C}}_k - \tns{C}_*}_2
        \end{align}
    \end{subequations}
    as required. In the last step we used \(\norm{\mat{P}_k}_2 \leq 1\).
\end{proof}

\subsection{\texorpdfstring{Convergence for strongly \(\R^n\)-spanning steps}{Convergence for strongly Rn-spanning steps}}\label{sec:convergence-weakly-orthogonal}

Since the updates only have access to \(\widetilde{\tns{C}}_k[\vek{s}_k]\) in each step, we can never hope to recover the true \(p\)th derivative if from some point onward all steps lie in a low-dimensional subspace of \(\R^n\), so we must assume that the steps repeatedly span \(\R^n\).
Indeed, we will assume something slightly stronger, namely
\begin{equation}\label{eqn:orthogonality-condition}
    \Norm{\prod_{k=k_0}^{k_0+m-1} \Paren{\mat{I} - \frac{\mat{W}_*^{-1}\vek{s}_k \vek{s}_k^T \mat{W}_*^{-T}}{\vek{s}_k^T \mat{W}_*^{-T} \mat{W}_*^{-1} \vek{s}_k}}}_2 \leq c < 1
\end{equation}
for some fixed \(m \in \N\), \(c \in \R_{\geq 0}\) and all \(k_0 \in \N\) large enough.
Here, \(\mat{W}_* = \lim_{k \to \infty} \mat{W}_k\) which we assume to be nonsingular as well.
As Moré and Trangenstein \cite{more_global_1976} showed, this assumption is equivalent to uniform linear independence of the scaled steps \(\mat{W}_*^{-1}\vek{s}_k\) which is in turn equivalent to the standard assumption of uniform linear independence of the steps \(\vek{s}_k\) themselves.
Intuitively, the steps being uniformly linearly independent means that every \(m\) consecutive steps span \(\R^n\) and they do so in a way that does not get arbitrarily degenerate.

Under this assumption it is possible to show that the approximations generated by \cref{eqn:update-definition} will converge to the true derivative without assuming Lipschitz continuity of the function or its derivatives.

\begin{theorem}\label{thm:convergence-assuming-orthogonality}
    Let \(\tns{C}_0 \in \R^{\otimes^p n}_{\sym}\) be given and update the approximations \(\tns{C}_k\) according to \cref{eqn:update-definition}.
    Assume \(\vek{x}_k\) converge to \(\vek{x}_* \in \R^n\), \(\mat{W}_k\) converge to some nonsingular matrix \(\mat{W}_* \in \R^{n \times n}\) and the steps are uniformly linearly independent.
    Then \(\tns{C}_k\) converges to \(\tns{C}_* \deq D^p f(\vek{x}_*)\).
\end{theorem}

For this and the other results to follow it suffices to consider the case when \(\mat{W}_* = \mat{I}\).
Otherwise, let \(\bar{f}(\vek{x}) = f(\mat{W}_* \vek{x})\), \(\bar{\vek{x}}_k = \mat{W}_*^{-1}\vek{x}_k\), \(\bar{\vek{x}}_* = \mat{W}_*^{-1}\vek{x}_*\), \(\bar{\mat{W}}_k = \mat{W}_*^{-1}\mat{W}_k\) and \(\bar{\tns{C}}_0 = \tns{C}_0[\mat{W}_*]^p\) and update \(\bar{\tns{C}}_k\) according to the adapted \cref{eqn:update-definition}.
Clearly, the assumptions of \cref{thm:convergence-assuming-orthogonality} are still satisfied for \(\bar{\vek{x}}_k\) and \(\bar{\mat{W}}_k\) and additionally \(\bar{\mat{W}}_k \to \mat{I}\).
By construction, we then have \(\tns{C}_k = \bar{\tns{C}}_k[\mat{W}_*^{-1}]^p\) so that \(\bar{\tns{C}}_k \to D^p \bar{f}(\bar{\vek{x}}_*)\) implies \(\tns{C}_k \to D^p f(\vek{x}_*)\).
Scaling by \(\mat{W}_*^{-1}\) transforms the sequence of approximations into one that gets arbitrarily close to employing the analogue of the PSB update and enables us to use its convergence properties.
Note that the assumption that \(\mat{W}_*\) is nonsingular is crucial for this transformation.

Inside the proof there appear three different matrices that are related to the projection matrices \(\mat{P}_k\):
\begin{subequations}
    \begin{align}
        \mat{P}_k   & = \mat{I} - \frac{\mat{W}_k^{-1}\vek{s}_k \vek{s}_k^T \mat{W}_k^{-T}}{\vek{s}_k^T \mat{W}_k^{-T} \mat{W}_k^{-1} \vek{s}_k}                                       \\
        \mat{P}_k^* & = \mat{I} - \frac{\mat{W}_*^{-1}\vek{s}_k \vek{s}_k^T \mat{W}_*^{-T}}{\vek{s}_k^T \mat{W}_*^{-T} \mat{W}_*^{-1} \vek{s}_k}                                       \\
        \mat{P}'_k  & = \mat{W}_k \mat{P}_k \mat{W}_k^{-1} = \mat{I} - \frac{\vek{s}_k \vek{s}_k^T \mat{W}_k^{-T} \mat{W}_k^{-1}}{\vek{s}_k^T \mat{W}_k^{-T} \mat{W}_k^{-1} \vek{s}_k}
    \end{align}
\end{subequations}
The next lemma shows that as \(\mat{W}_k \to \mat{W}_* = \mat{I}\) the distance between these matrices gets arbitrarily small.

\begin{lemma}\label{thm:projection-matrices-convergence}
    Let the nonsingular matrices \(\mat{W}_k \in \R^{n \times n}\) converge to \(\mat{W}_* = \mat{I}\) then
    \begin{equation}\label{eqn:projection-matrices-convergence}
        \norm{\mat{P}_k - \mat{P}_k^*}_2 \to 0 \quad \text{and} \quad \norm{\mat{P}'_k - \mat{P}_k^*}_2 \to 0
    \end{equation}
    holds for any sequence of steps \((\vek{s}_k)_{k \in \N}\) as \(k \to \infty\).
\end{lemma}
\begin{proof}
    Without loss of generality we can assume that the steps are scaled such that \(\norm{\vek{s}_k}_2 = 1\) for all \(k \in \N\).
    That means
    \begin{equation}
        \norm{\mat{W}_k^{-1}\vek{s}_k - \vek{s}_k}_2 \leq \underbrace{\norm{\mat{W}_k^{-1} - \mat{I}}_2}_{\to 0} \underbrace{\norm{\vek{s}_k}_2}_{=1} \to 0 \quad \text{as } k \to \infty.
    \end{equation}
    Since \(\mat{P}_k\) projects onto the subspace that is orthogonal to \(\mat{W}_k^{-1}\vek{s}_k\) and \(\mat{P}_k^*\) projects onto the subspace that is orthogonal to \(\vek{s}_k\), the difference between the two projection matrices converges to zero. This is the first claim.

    For the second one we find that
    \begin{equation}
        \mat{P}'_k - \mat{P}_k^* = \underbrace{(\mat{W}_k - \mat{I})}_{\to \mat{0}}\mat{P}_k\underbrace{\mat{W}_k^{-1}}_{\to \mat{I}} + \underbrace{(\mat{P}_k - \mat{P}_k^*)}_{\to \mat{0}} \underbrace{\mat{W}_k^{-1}}_{\to \mat{I}} + \mat{P}_k^* \underbrace{(\mat{W}_k^{-1} - \mat{I})}_{\to \mat{0}}
    \end{equation}
    converges to \(\mat{0}\) since the norms of \(\mat{P}_k\) and \(\mat{P}_k^*\) are at most one.
\end{proof}

\begin{proof}[Proof of \cref{thm:convergence-assuming-orthogonality}]
    As argued above, we will only consider the case \(\mat{W}_* = \mat{I}\).
    Let \(\e > 0\) be arbitrary. To prove convergence, we will show \(\norm{\tns{C}_{k} - \tns{C}_*}_2 \leq \e\) for every \(k\) large enough.
    Consider the recurrence relation for \(\tns{C}_k - \tns{C}_*\) established in \cref{thm:error-tensor-bound}.
    We can rearrange it to read
    \begin{equation}
        \Paren{\tns{C}_{k+1} - \tns{C}_*} = \Paren{\tns{C}_k - \tns{C}_*}[\mat{W}_k \mat{P}_k \mat{W}_k^{-1}]^p + \tns{E}_k = \Paren{\tns{C}_k - \tns{C}_*}[\mat{P}'_k]^p + \tns{E}_k.
    \end{equation}
    Note that \(\mat{P}'_k\) is not an orthogonal projection but by \cref{thm:projection-matrices-convergence} it approaches one as \(k \to \infty\).
    In particular, its norm converges to one since \(\norm{\mat{P}'_k}_2 \leq \kappa_2(\mat{W}_k) \to 1\).
    To use the assumption \cref{eqn:orthogonality-condition} later on we apply the previous equality \(m\) times and get
    \begin{equation}\label{eqn:m-step-recurrence}
        \Paren{\tns{C}_{k_0+m} - \tns{C}_*} = \Paren{\tns{C}_{k_0} - \tns{C}_*}\Brack{\prod_{k=k_0}^{k_0+m-1} \mat{P}'_k}^p + \sum_{k=k_0}^{k_0+m-1} \tns{E}_k\Brack{\prod_{l=k+1}^{k_0+m-1} \mat{P}'_l}^p.
    \end{equation}

    Let \(K_1 \in \N\) be such that for all \(k \geq K_1\) we have \(\kappa_2(\mat{W}_k) \leq 2\).
    This means using \cref{thm:error-tensor-bound} we can bound the second term on the right-hand side of the previous equation by
    \begin{equation}
        \Norm{\sum_{k=k_0}^{k_0+m-1} \tns{E}_k\Brack{\prod_{l=k+1}^{k_0+m-1} \mat{P}'_l}}_2 \leq \sum_{k=k_0}^{k_0+m-1} 2^p (1+2^p) \norm{\widetilde{\tns{C}}_k - \tns{C}_*}
    \end{equation}
    for all \(k_0 \geq K_1\).
    Since \cref{thm:C-tilde-convergence} showed that \(\widetilde{\tns{C}}_k \to \tns{C}_*\) the bound above is smaller than \(\e/2 \cdot (1 - (\frac{1+c}{2})^p)\) for \(k_0\) large enough, say \(k_0 \geq K_2\).

    Next, consider first term on the right-hand side of \cref{eqn:m-step-recurrence}.
    It features a product of \(\mat{P}'_k\) whereas our assumption of uniform linear independence \cref{eqn:orthogonality-condition} features a product of \(\mat{P}_k^*\).
    We find that
    \begin{subequations}\label{eqn:projection-product-distance}
        \begin{align}
            \Norm{\prod_{k=k_0}^{k_0+m-1} \mat{P}'_k - \prod_{k=k_0}^{k_0+m-1} \mat{P}_k^*}_2 & = \Norm{\sum_{k=k_0}^{k_0+m-1} \mat{P}_{k_0}^* \cdots \mat{P}_{k-1}^* (\mat{P}'_k - \mat{P}_k^*) \mat{P}'_{k+1} \cdots \mat{P}'_{k_0+m-1}}_2 \\
                                                                                              & \leq \sum_{k=k_0}^{k_0+m-1} 2^{m-1} \norm{\mat{P}'_k - \mat{P}_k^*}_2
        \end{align}
    \end{subequations}
    for \(k_0 \geq K_1\) using \(\norm{\mat{P}'_k} \leq 2\) and \(\norm{\mat{P}_k^*} \leq 1 \leq 2\).
    For \(k_0\) large enough, say \(k_0 \geq K_3\), the right-hand side of \cref{eqn:projection-product-distance} is smaller than \((1-c)/2\) by \cref{thm:projection-matrices-convergence}.
    Similarly, for \(k_0\) large enough, say \(k_0 \geq K_4\), we have \(\norm{\prod_{k=k_0}^{k_0 + m - 1} \mat{P}_k^*} \leq c\) by assumption.
    Therefore, for \(k_0 \geq \max\{K_3, K_4\}\)
    \begin{equation}
        \Norm{\prod_{k=k_0}^{k_0 + m - 1} \mat{P}'_k} \leq \Norm{\prod_{k=k_0}^{k_0+m-1} \mat{P}'_k - \prod_{k=k_0}^{k_0+m-1} \mat{P}_k^*}_2 + \Norm{\prod_{k=k_0}^{k_0 + m - 1} \mat{P}_k^*} \leq \frac{1 + c}{2}.
    \end{equation}

    In the previous two paragraphs we have established that asymptotically for every \(m\) steps the norm of \(\tns{C}_k - \tns{C}_*\) is first multiplied by a factor that is at most slightly larger than \(c\) and then increased by an arbitrarily small error term.
    In particular for \(k_0 \geq K = \max\{K_2, K_3, K_4\}\) we have
    \begin{equation}
        \norm{\tns{C}_{k_0+m} - \tns{C}_*}_2 \leq \norm{\tns{C}_{k_0} - \tns{C}_*}_2 \Paren{\frac{1+c}{2}}^p + \e/2 \cdot \Paren{1 - \Paren{\frac{1+c}{2}}^p}.
    \end{equation}
    Repeatedly applying this inequality shows that \(\norm{\tns{C}_{k_0+im} - \tns{C}_*}_2\) is bounded by a sequence \((a_i)_{i \in \N}\) which converges to \(\e/2\).
    Use this observation for all \(k_0 \in \{K, K+1, \dots, K+m-1\}\) to see that \(\norm{\tns{C}_k - \tns{C}_*}_2 \leq \e\) for all \(k\) large enough, as claimed.
\end{proof}

\Cref{thm:convergence-assuming-orthogonality} can be seen as a global convergence result for the approximations \(\tns{C}_k\).
The main assumptions are that \(\mat{W}_k\) converges to \(\mat{W}_*\), that \(\mat{W}_*\) is nonsingular, and that the steps \(\vek{s}_k\) are uniformly linearly independent.
The first one might already look nonsensical given that for any update step one can replace the weight matrix \(\mat{W}_k\) by another one from an infinite family of matrices without changing the sequence of approximations.
It really should be understood as the requirement that \emph{there exists} a sequence of weight matrices compatible with the updates which converges to a nonsingular matrix \(\mat{W}_*\) or, in other words, we require the relationship between \(\vek{s}_k\) and \(\vek{v}_k\) to be become linear as \(k \to \infty\), namely \(\vek{v}_k \approx \mat{W}_*^{-T} \mat{W}_*^{-1} \vek{s}_k\).
For PSB and DFP methods this assumption is satisfied.
For PSB the update is compatible with choosing \(\mat{W}_k = \mat{I}\) as the weight matrix, so clearly this sequence converges to a nonsingular matrix.
The DFP update is compatible with \(\mat{W}_k = \widetilde{\mat{B}}_k^{-1/2}\) where \(\widetilde{\mat{B}}_k\) is defined as the averaged true Hessian on the line from \(\vek{x}_k\) to \(\vek{x}_{k+1}\).
Clearly, \(\widetilde{\mat{B}}_k \to \nabla^2 f(\vek{x}_*)\) as \(\vek{x}_k \to \vek{x}_*\) and so assuming positive definiteness of the Hessian at the limit point we get that \(\mat{W}_*\) is positive definite as well.

For the assumption of uniform linear independence of the steps we already argued that it is necessary to assume that the steps repeatedly span \(\R^n\) to have any hope of convergence.
The assumption, however, is indeed stronger than that, so we might ask whether this is warranted.
Ge and Powell \cite{ge_convergence_1983} give an example where the DFP method fails to produce a sequence of matrices that converges to the true Hessian even in the case where the function is a quadratic function of two variables and the Hessian is the identity.
The steps are chosen in such a way that the angle between consecutive steps converges to zero as \(1/k\).
This shows that repeatedly spanning \(\R^n\) does not suffice to ensure convergence as any two consecutive steps do span \(\R^2\) in the example.
Uniform linear independence is sufficiently strong to rule out these steps as it would require that there is an \(m \in \N\) such that the maximum angle between two of \(m\) consecutive steps is uniformly bounded below throughout the sequence.

To further emphasize the relevance of these assumptions, we take a look at the existing literature of convergence of classical quasi-Newton approximations.
Powell \cite[Theorem 5]{powell_new_1970} showed in the same paper that introduced the PSB update that, assuming boundedness and Lipschitz continuity of the second derivative, as well as uniform linear independence of the steps, the approximations converge to the true Hessian.
The proof has many similarities to the one shown here, except for the fact that neither boundedness nor Lipschitz continuity are necessary for our result.
It is interesting that Powell suggests making every third iteration of his optimization method a ``special iteration'' in order to enforce the uniform linear independence assumption in practice.\footnote{
    In special iterations the last \(2n\) steps are checked to see whether they satisfy a uniform linear independence condition and if necessary a step pointing towards the missing direction is introduced.
    The required bookkeeping to identify such directions is described in \cite[Section 7]{powell_fortran_1968}.
    One could also revise this scheme to introduce the missing direction not as an additional step, but as a perturbation of a computed step, potentially saving one function and derivative evaluation whenever such a correction is needed.
}
This idea is not found in modern implementations of quasi-Newton methods and uniform linear independence remains a theoretical device that cannot be ensured in practice.

Conn, Gould and Toint \cite{conn_convergence_1991} provide a global convergence theorem for matrices generated by the SR1 update.
Just as Powell \cite{powell_new_1970} they aim to cover the trust-region case and therefore do not impose any specific structure of the steps.
Instead, they again assume Lipschitz continuity of the Hessian and uniform linear independence of the steps as well as a bound on the near-orthogonality of \(\vek{s}_k\) and \((\mat{B}_k - \widetilde{\mat{B}}_k)\vek{s}_k\).
Note that our theorem above does not cover the SR1 case since the mapping from \(\vek{s}_k\) to \(\vek{v}_k = (\mat{B}_k - \widetilde{\mat{B}}_k)\vek{s}_k\) does not approach a constant nonsingular linear map.

For DFP and BFGS updates Ge and Powell \cite{ge_convergence_1983} were able to show that assuming Lipschitz continuity of the Hessian, positive definiteness of the Hessian at the limit point and steps of the form \(\vek{s}_k = -\mat{B}_k^{-1} \nabla f(\vek{x}_k)\) the sequence of quasi-Newton matrices starting sufficiently close to the true Hessian will converge, although not necessarily to the true Hessian.
They drop the uniform linear independence assumption but achieve a weaker result in return.

The case of DFP and BFGS convergence for unstructured steps was covered in a very technical paper by Boggs and Tolle \cite{boggs_convergence_1994}.
In the case of a quadratic function \(f\) with nonsingular Hessian they are able to show global convergence of the DFP approximations using a notion that is slightly weaker than uniform linear independence. 
The same is true for BFGS with the added assumption that the initial approximation is nonsingular. 
For general functions \(f\) they are able to generalize convergence of DFP but not of BFGS updates. 
Their main contribution however is that they cover the case where the steps asymptotically fall into a subspace and we only have uniform linear independence for the projected steps.
Say that subspace is spanned by the columns of \(\mat{V}\), then we get convergence of \(\mat{B}_k \mat{V}\) to \(\nabla^2 f(\vek{x}_*) \mat{V}\) in all previously discussed cases, which is the best we can hope for.

Even though \cref{thm:convergence-assuming-orthogonality} considers a different family of updates, if we specialize to PSB or DFP updates and the case \(p=2\) we recover a slightly more general statement than the one given by Powell \cite{powell_new_1970} and a slightly weaker statement than the one given by Boggs and Tolle \cite{boggs_convergence_1994} respectively.
Fundamentally though, the results in the literature also require (a variant of) uniform linear independence and in the case of DFP positive definiteness of the Hessian at \(\vek{x}_*\) and so \cref{thm:convergence-assuming-orthogonality} subsumes and generalizes these convergence theorems.

\subsection{Generalized Dennis--Mor\'e condition}\label{sec:dennis-more}

Dennis and Moré \cite{dennis_characterization_1974} showed that if optimization methods choose their iterates using \(\vek{x}_{k+1} = \vek{x}_k - \mat{B}_k^{-1} \nabla f(\vek{x}_k)\) and converge to \(\vek{x}_*\) then this convergence is Q-superlinear if and only if
\begin{equation}\label{eqn:dennis-more}
    \lim_{k \to \infty} \frac{\norm{(\mat{B}_k - \nabla^2 f(\vek{x}_*))\vek{s}_k}_F}{\norm{\vek{s}_k}_2} = 0.
\end{equation}
This equation is called the \emph{Dennis--Mor\'e condition} and is weaker than convergence of \(\mat{B}_k\) to \(\nabla^2 f(\vek{x}_*)\).
Instead, it suffices when \(\mat{B}_k\) converges along the step directions.

We are not concerned with convergence rates for any particular optimization method in this paper, but it seems natural to ask whether the approximations \(\tns{C}_k\) satisfy a generalized Dennis--Mor\'e condition
\begin{equation}\label{eqn:generalized-dennis-more}
    \lim_{k \to \infty} \frac{\norm{(\tns{C}_k - D^p f(\vek{x}_*))[\vek{s}_k]}_F}{\norm{\vek{s}_k}_2} = 0
\end{equation}
when they are updated according to \cref{eqn:update-definition}.

\begin{theorem}\label{thm:convergence-dennis-more}
    Let \(\tns{C}_0 \in \R^{\otimes^p n}_{\sym}\) be given and update the approximations \(\tns{C}_k\) according to \cref{eqn:update-definition} where the function has a Lipschitz continuous \(p\)th derivative \(D^p f\).
    Assume \(\vek{x}_k\) converge to \(\vek{x}_* \in \R^n\) and \(\mat{W}_k\) converge to some nonsingular matrix \(\mat{W}_* \in \R^{n \times n}\) fast enough such that
    \begin{equation}
        \sum_{k \geq 0} \norm{\vek{x}_k - \vek{x}_*}_2 < \infty \quad \text{and} \quad \sum_{k \geq 0} \norm{\mat{W}_k - \mat{W}_*}_2 < \infty.
    \end{equation}
    Then the generalized Dennis--Mor\'e condition \cref{eqn:generalized-dennis-more} holds.
\end{theorem}
\begin{proof}
    As in the proof of \cref{thm:convergence-assuming-orthogonality} we assume \(\mat{W}_* = \mat{I}\) without loss of generality.
    Since the condition number is locally Lipschitz continuous around \(\mat{I}\), the assumption \(\sum_{k \geq 0} \norm{\mat{W}_k - \mat{I}}_2 < \infty\) also implies \(C_{\kappa} \deq \sum_{k \geq 0} \Paren{\kappa_2(\mat{W}_k) - 1} < \infty\).

    As a first step we need to use the bounded deterioration principle to show that \(\norm{\tns{C}_k - \tns{C}_*}_2\) stays bounded.
    For this we again consider the \(m\)-step recursive formula established in \cref{eqn:m-step-recurrence}:
    \begin{equation}
        \Paren{\tns{C}_m - \tns{C}_*} = \Paren{\tns{C}_0 - \tns{C}_*}\Brack{\prod_{k=0}^{m-1} \mat{P}'_k}^p + \sum_{k=0}^{m-1} \tns{E}_k \Brack{\prod_{l=k+1}^{m-1} \mat{P}'_l}^p
    \end{equation}
    where \(\mat{P}'_k = \mat{W}_k \mat{P}_k \mat{W}_k^{-1}\) and \(\tns{E}_k\) is defined in \cref{thm:error-tensor-bound}.
    Clearly,
    \begin{equation}
        \ln\Paren{\Norm{\prod_{k=0}^{m-1} \mat{P}'_k}_2} \leq \sum_{k=0}^{m-1} \ln\Paren{\norm{\mat{P}'_k}_2} \leq \sum_{k=0}^{m-1} \ln\Paren{\kappa_2(\mat{W}_k)} \leq \sum_{k=0}^{m-1} \Paren{\kappa_2(\mat{W}_k) - 1} \leq C_{\kappa}
    \end{equation}
    is uniformly bounded for all \(m\), which means the same is true for \(\norm{\prod_{k=0}^{m-1} \mat{P}'_k}_2\) itself.
    Therefore,
    \begin{equation}
        \Norm{\tns{C}_m - \tns{C}_*}_2 \leq \exp(C_{\kappa})^p \Paren{ \Norm{\tns{C}_0 - \tns{C}_*}_2 + \sum_{k=0}^{m-1} \norm{\tns{E}_k}_2 }.
    \end{equation}
    Because \(D^p f\) is Lipschitz continuous and \(\kappa_2(\mat{W}_k)\) stays bounded, \cref{thm:C-tilde-convergence,thm:error-tensor-bound} show that there is a constant \(C_{\tns{E}} < \infty\) such that
    \begin{equation}\label{eqn:dennis-more-error-Lipschitz-bound}
        \norm{\tns{E}_k}_2 \leq C_{\tns{E}}/2 \Paren{\norm{\vek{x}_k - \vek{x}_*} + \norm{\vek{x}_{k+1} - \vek{x}_*}}.
    \end{equation}
    This implies \(\sum_{k=0}^{m-1} \norm{\tns{E}_k}_2 \leq C_{\tns{E}} \sum_{k=0}^m \norm{\vek{x}_k - \vek{x}_*}_2\) is uniformly bounded for all \(m\) and therefore \(\norm{\tns{C}_m - \tns{C}_*}_2\) is as well.

    Now we are ready to tackle the main claim.
    To introduce the quantity of interest \(\norm{(\tns{C}_k - \tns{C}_*)[\vek{s}_k]}_F / \norm{\vek{s}_k}_2\) we use a trick similar to the one used in the proof of Theorem 8.2.2 in \cite{dennis_numerical_1996}.
    Note that in the Frobenius norm for any tensor \(\tns{T} \in \R^{\otimes^p n}\) and any nonzero vector \(\vek{w} \in \R^n\) we have
    \begin{equation}
        \norm{\tns{T}}_F^2 = \Norm{\tns{T}\Brack{\frac{\vek{w} \vek{w}^T}{\vek{w}^T \vek{w}}}}_F^2 + \Norm{\tns{T}\Brack{\mat{I} - \frac{\vek{w} \vek{w}^T}{\vek{w}^T \vek{w}}}}_F^2 = \frac{\Norm{\tns{T}\Brack{\vek{w}}}_F^2}{\norm{\vek{w}}_2^2} + \Norm{\tns{T}\Brack{\mat{I} - \frac{\vek{w} \vek{w}^T}{\vek{w}^T \vek{w}}}}_F^2
    \end{equation}
    because the matrices in brackets are orthogonal projections.
    We can apply this to the case where \(\mat{I} - \vek{w}\vek{w}^T/(\vek{w}^T \vek{w})\) is \(\mat{P}_k\) to get
    \begin{subequations}\label{eqn:frobenius-norm-bound-dennis-more}
        \begin{align}
             & \phantom{{}={}} \Norm{\Paren{\tns{C}_{k} - \tns{C}_*}\Brack{\mat{W}_k}^p \Brack{\mat{P}_k}^p}_F^2                                                                                                          \\
             & \leq \Norm{\Paren{\tns{C}_{k} - \tns{C}_*}\Brack{\mat{W}_k}^p \Brack{\mat{P}_k}}_F^2                                                                                                                       \\
             & = \Norm{\Paren{\tns{C}_{k} - \tns{C}_*}\Brack{\mat{W}_k}^p}_F^2 - \frac{\Norm{\Paren{\tns{C}_{k} - \tns{C}_*}\Brack{\mat{W}_k}^p \Brack{\mat{W}_k^{-1}\vek{s}_k}}_F^2}{\norm{\mat{W}_k^{-1}\vek{s}_k}_2^2} \\
             & = \Norm{\Paren{\tns{C}_{k} - \tns{C}_*}\Brack{\mat{W}_k}^p}_F^2 - \frac{\Norm{\Paren{\tns{C}_{k} - \tns{C}_*}[\vek{s}_k]\Brack{\mat{W}_k}^{p-1}}_F^2}{\norm{\mat{W}_k^{-1}\vek{s}_k}_2^2}                  \\
             & \leq \Norm{\Paren{\tns{C}_{k} - \tns{C}_*}\Brack{\mat{W}_k}^p}_F^2 - \frac{\Norm{\Paren{\tns{C}_{k} - \tns{C}_*}[\vek{s}_k]}_F^2}{\norm{\vek{s}_k}_2^2} \norm{\mat{W}_k^{-1}}_2^{-2p}.
        \end{align}
    \end{subequations}
    Adding \(\tns{E}_k[\mat{W}_k]^p\) into the Frobenius norm on the left-hand side gives
    \begin{subequations}\label{eqn:frobenius-norm-bound-recurrence}
        \begin{align}
            \norm{\Paren{\tns{C}_{k+1} - \tns{C}_*}[\mat{W}_k]^p}_F^2
             & = \norm{\Paren{\tns{C}_{k} - \tns{C}_*}\Brack{\mat{W}_k}^p \Brack{\mat{P}_k}^p + \tns{E}_k\Brack{\mat{W}_k}}_F^2                                                                \\
             & = \norm{\Paren{\tns{C}_{k} - \tns{C}_*}\Brack{\mat{W}_k}^p \Brack{\mat{P}_k}^p}_F^2 + \norm{\tns{E}_k\Brack{\mat{W}_k}^p}_F^2. \label{eqn:frobenius-norm-bound-recurrence-last}
        \end{align}
    \end{subequations}

    The inner product term \(\Innerprod{\Paren{\tns{C}_{k} - \tns{C}_*}\Brack{\mat{W}_k}^p \Brack{\mat{P}_k}^p}{\tns{E}_k\Brack{\mat{W}_k}^p}_F\) is missing in \cref{eqn:frobenius-norm-bound-recurrence-last} because it is zero.
    To show that, note that we can rewrite \cref{eqn:error-tensor-C-tilde-recurrence} from the proof of \cref{thm:error-tensor-bound} as
    \begin{equation}
        \paren{\widetilde{\tns{C}}_k - \tns{E}_k - \tns{C}_*}[\mat{W}_k]^p = \paren{\widetilde{\tns{C}}_k - \tns{C}_*}[\mat{W}_k]^p \Brack{\mat{P}_k}^p.
    \end{equation}
    Using the equivalence between \cref{thm:characterization} \labelcref{item:characterization-recursive} and \labelcref{item:characterization-low-rank} the error tensor can be written explicitly as \(-\tns{E}_k = P_{\sym}(\tns{A} \otimes \mat{W}_k^{-T}\mat{W}_k^{-1}\vek{s}_k)\) for some \((p-1)\)-tensor \(\tns{A}\) and
    \begin{equation}
        \tns{E}_k\Brack{\mat{W}_k}^p = -P_{\sym}(\tns{A}\Brack{\mat{W}_k}^{p-1} \otimes \mat{W}_k^{-1}\vek{s}_k).
    \end{equation}
    In other words, \(\tns{E}_k \Brack{\mat{W}_k}^p\) can be expressed as a sum of outer products between \(\mat{W}_k^{-1}\vek{s}_k\) and some \((p-1)\)-tensor.
    Any inner product of such a tensor with \(\Paren{\tns{C}_{k} - \tns{C}_*}\Brack{\mat{W}_k}^p \Brack{\mat{P}_k}^p\) must be zero as \(\mat{P}_k\) maps \(\mat{W}_k^{-1}\vek{s}_k\) to zero.

    Combining \cref{eqn:frobenius-norm-bound-dennis-more,eqn:frobenius-norm-bound-recurrence} and rearranging gives
    \begin{subequations}\label{eqn:dennis-more-bound-internal}
        \begin{align}
            \frac{\Norm{\Paren{\tns{C}_{k} - \tns{C}_*}[\vek{s}_k]}_F^2}{\norm{\vek{s}_k}_2^2} & \leq \norm{\mat{W}_k^{-1}}_2^{2p} \big(\norm{\Paren{\tns{C}_k - \tns{C}_*}\Brack{\mat{W}_k}^p}_F^2 \notag                                                             \\
                                                                                               & \phantom{\leq \norm{\mat{W}_k^{-1}}_2^{2p} \big(} - \norm{\Paren{\tns{C}_{k+1} - \tns{C}_*}\Brack{\mat{W}_k}^p}_F^2 + \norm{\tns{E}_k\Brack{\mat{W}_k}^p}_F^2 \big)   \\
                                                                                               & \leq \norm{\Paren{\tns{C}_k - \tns{C}_*}}_F^2 \kappa_2(\mat{W}_k)^{2p} - \norm{\Paren{\tns{C}_{k+1} - \tns{C}_*}}_F^2 + \norm{\tns{E}_k}_F^2 \kappa_2(\mat{W}_k)^{2p}
        \end{align}
    \end{subequations}
    Lastly, we wish to sum up both sides of the previous inequality over all \(k \geq 0\) and show that the right-hand side stays bounded.
    This immediately gives the claim.
    A few technicalities are needed.
    We already showed that \(\tns{C}_k - \tns{C}_*\) stays bounded, so let \(C_{\Delta} < \infty\) be a constant such that \(\norm{\tns{C}_k - \tns{C}_*}_F \leq C_{\Delta}\) for all \(k \in \N\).
    As established above, \(\sum_{k\geq 0} \Paren{\kappa_2(\mat{W}_k) - 1} < \infty\).
    This also implies that \(\sum_{k\geq 0} \Paren{\kappa_2(\mat{W}_k)^{2p} - 1} = C_{\kappa, 2p} < \infty\) for some constant \(C_{\kappa, 2p}\) since \(\kappa_2(\mat{W}_k)^{2p} - 1 \leq 4p \Paren{\kappa_2(\mat{W}_k) - 1}\) for \(\kappa_2(\mat{W}_k)\) small enough.

    Consider the \(\tns{C}_k - \tns{C}_*\) terms on the right-hand side of \cref{eqn:dennis-more-bound-internal} first:
    \begin{subequations}
        \begin{align}
             & \phantom{={}} \sum_{k=0}^{m-1} \Paren{\norm{\Paren{\tns{C}_k - \tns{C}_*}}_F^2 \kappa_2(\mat{W}_k)^{2p} - \norm{\Paren{\tns{C}_{k+1} - \tns{C}_*}}_F^2}                                                                                                                                                         \\
             & = \underbrace{\norm{\tns{C}_0 - \tns{C}_*}_F^2 \kappa_2(\mat{W}_0)^{2p}}_{\text{constant}} + \underbrace{\sum_{k=1}^{m-1} \Paren{\kappa_2(\mat{W}_k)^{2p} - 1}}_{= C_{\kappa, 2p}} \underbrace{\norm{\tns{C}_k - \tns{C}_*}_F^2}_{\leq C_{\Delta}^2} - \underbrace{\norm{\tns{C}_{m} - \tns{C}_*}_F^2}_{\geq 0}
        \end{align}
    \end{subequations}
    is uniformly bounded for all \(m\).
    The same is true for the \(\tns{E}_k\) term.
    Because the Frobenius norm and the 2-norm for \(p\)-tensors are both norms on finite-dimensional vector spaces they are equivalent.
    Moreover, \(\kappa_2(\mat{W}_k)\) stays bounded, so for some constant \(C\)
    \begin{subequations}
        \begin{align}
            \sum_{k=0}^{m-1} \norm{\tns{E}_k}_F^2 \kappa_2(\mat{W}_k)^{2p} \leq C \sum_{k=0}^{m-1} \norm{\tns{E}_k}_2^2 \leq C \Paren{\sum_{k=0}^{m-1} \norm{\tns{E}_k}_2}^2
        \end{align}
    \end{subequations}
    holds, and the term is uniformly bounded.
    Therefore,
    \begin{equation}
        \sum_{k \geq 0} \frac{\Norm{\Paren{\tns{C}_{k} - \tns{C}_*}[\vek{s}_k]}_F^2}{\norm{\vek{s}_k}_2^2} < \infty \quad \text{and} \quad \frac{\Norm{\Paren{\tns{C}_{k} - \tns{C}_*}[\vek{s}_k]}_F}{\norm{\vek{s}_k}_2} \to 0 \text{ as } k \to \infty
    \end{equation}
    as claimed.
\end{proof}

Let us compare this result to the one obtained by Dennis and Mor\'e in the paper that introduced the Dennis--Mor\'e condition and showed its relevance for superlinear convergence \cite{dennis_characterization_1974}.
To explain superlinear convergence of existing quasi-Newton methods in this new framework they needed to establish that these update formulas indeed satisfy condition \cref{eqn:dennis-more}.
For the DFP update they do so by assuming Lipschitz continuity of the Hessian, positive definiteness of the Hessian at \(\vek{x}_*\) and boundedness of \(\sum_{k \geq 0} \norm{\vek{x}_k - \vek{x}_*}\).\footnote{
    Strictly speaking, they only assume existence of constants \(L, \alpha > 0\) such that \(\norm{\nabla^2 f(\vek{x}) - \nabla^2 f(\vek{x}_*)} \leq L \norm{\vek{x} - \vek{x}_*}^{\alpha}\) for any \(\vek{x}\) and correspondingly \(\sum_{k \geq 0} \norm{\vek{x}_k - \vek{x}_*}^{\alpha}\) on top of positive definiteness of \(\nabla^2 f(\vek{x}_*)\).
    By adapting \cref{eqn:C-tilde-Lipschitz-bound,eqn:dennis-more-error-Lipschitz-bound} appropriately, our proof also covers this case.
}
Under these same assumptions \cref{thm:convergence-dennis-more} also implies convergence of the DFP matrices.
As mentioned before the DFP method is compatible with \(\mat{W}_k = \widetilde{\mat{B}}_k^{-1/2} \to \mat{W}_* = \nabla^2 f(\vek{x}_*)^{-1/2}\).
Since \(\nabla^2 f(\vek{x}_*)\) is positive definite, \(\mat{W}_*\) is well-defined and nonsingular.
Moreover, \(\mat{A} \mapsto \mat{A}^{-1/2}\) is differentiable at any positive definite matrix, so the map is also locally Lipschitz around \(\nabla^2 f(\vek{x}_*)\) and \(\sum_{k \geq 0} \norm{\vek{x}_k - \vek{x}_*}_2 < \infty\) implies \(\sum_{k \geq 0} \norm{\mat{W}_k - \mat{W}_*}_2 < \infty\).
Therefore, for DFP all the assumptions of \cref{thm:convergence-dennis-more} are covered by the assumptions in \cite{dennis_characterization_1974} and vice versa.
The same is true for the PSB update, where Dennis and Mor\'e mention that the positive definiteness assumption can be dropped.
This agrees perfectly with our theorem since for \(\mat{W}_k = \mat{I}\) both existence of \(\mat{W}_*\) and \(\sum_{k \geq 0} \norm{\mat{W}_k - \mat{W}_*}_2 < \infty\) are obvious.
Again, this result extends the known cases to higher-order updates and clarifies the relevant convergence conditions for all updates that can be expressed as least-change updates in weighted Frobenius norms.

\section{Numerical experiments}\label{sec:numerical-experiments}

While the previous sections investigated the theoretical properties of the higher-order secant updates, we now turn to numerical experiments to understand how quickly convergence of the approximations sets in and how the algorithm behaves for different kinds of iterates.
This is not supposed to be a comprehensive treatment of the numerical performance of the algorithm, but rather give an idea of its general behaviour by considering a small toy problem.

The implementation was done in Python using NumPy \cite{harris_array_2020} and uses the explicit formula in \cref{thm:characterization} \labelcref{item:characterization-explicit} at its core:
\begin{equation}
    \tns{C}_{k+1} = \tns{C}_k + \sum_{j=1}^p (-1)^{j+1} \tbinom{p}{j} \Paren{\vek{v}_k^T \vek{s}_k}^{-j} P_{\sym} \Paren{\Paren{\otimes^j \vek{v}_k} \otimes \Paren{\tns{D}_k - \tns{C}_k[\vek{s}_k]}[\vek{s}_k]^{j-1}}
\end{equation}
where \(\vek{v}_k = \mat{W}_k^{-T} \mat{W}_k^{-1} \vek{s}_k\) and \(\tns{D}_k = D^{p-1} f(\vek{x}_{k+1}) - D^{p-1} f(\vek{x}_k)\).
This makes it particularly easy to implement the analogues of PSB (where \(\vek{v}_k = \vek{s}_k\)) and DFP (where \(\vek{v}_k = (\nabla f(\vek{x}_{k+1}) - \nabla f(\vek{x}_k))/\norm{\vek{s}_k}_2\)).
To increase legibility and since these two choices produce roughly similar approximations, only the PSB variant will be shown below.

We chose to use the two-dimensional Rosenbrock function \(f(x, y) = (1 - x)^2 + 100(y - x^2)^2\) to test our algorithm on because it is a simple, yet widely used test function with nonconstant third derivative
\begin{equation}
    D^3 f(x, y) = \begin{pmatrix}
        \begin{pmatrix}
            -2400x & -400 \\
            -400   & 0
        \end{pmatrix}
        \; \begin{pmatrix}
               -400 & 0 \\
               0    & 0
           \end{pmatrix}
    \end{pmatrix}.\footnote{We use the notation for 3-tensors introduced in \cref{ex:simple-update} here.}
\end{equation}
The iterates \(\vek{x}_k\) are generated by different minimization algorithms starting at \((0, 0)\) and converging to the global minimum of \(f\) at \(\vek{x}_* = (1, 1)\).
It is important to point out that the iterates are computed without using the approximated third derivatives.
In fact any number of sequences \(\vek{x}_k\) could have been chosen to explore the behaviour of the algorithm.
We chose sequences generated by optimization methods since they seem to be particularly relevant to our intended application, but it should be clear any algorithm that depends on the approximated third derivatives will produce different iterates to the ones we used.
Lastly, the initial approximation \(\tns{C}_0\) is just the zero 3-tensor in the following.

\subsection{Numerical limitations}

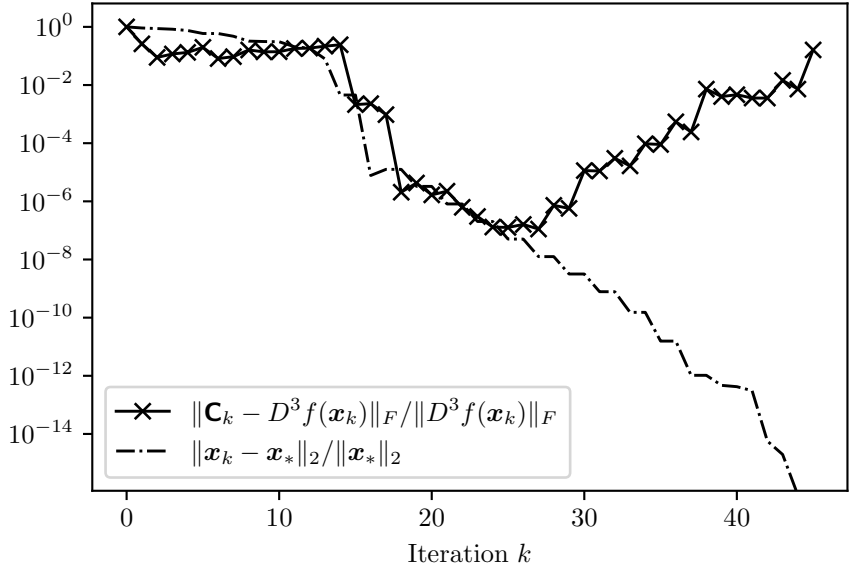
\begin{figure}
    \centering
    \subimport{plots/}{convergence_linearly_converging_min_simple.pgf}
    \caption{Convergence of iterates and approximations for nonlinear CG}
    \label{fig:convergence_linearly_converging_min_simple}
\end{figure}

In the first experiment, a nonlinear CG method\footnote{\texttt{scipy.optimize.minimize(method="CG")} in SciPy version 1.9.3, implementing the Polak--Ribi\`ere variant of nonlinear CG \cite[p. 122]{nocedal_numerical_2006}} was used to minimize the Rosenbrock function.
Although nonlinear CG methods with restarts can achieve superlinear local convergence, this implementation does not include restarts, and we observe roughly linear convergence of \(\vek{x}_k\) to \(\vek{x}_*\) in \cref{fig:convergence_linearly_converging_min_simple}.
The relative error of each \(\tns{C}_k\) is not measured with respect to \(\tns{C}_* = D^3 f(1,1)\) here but instead with respect to the true third derivative at each iterate \(D^3 f(\vek{x}_k)\), since this is the more relevant metric in practice.
From the convergence theorems in the previous section we expect that this quantity will also converge to zero, since both \(\tns{C}_k\) and \(D^3 f(\vek{x}_k)\) converge to \(\tns{C}_*\).

Unfortunately, although at first this seems to be true and the two error curves roughly coincide, from iteration 27 onwards the error in \(\tns{C}_k\) increases quite considerably.
The issue, as it turns out, stems from rounding errors in finite precision arithmetic, which we did not consider in the theory.
Specifically the computation of \(\tns{D}_k\) becomes more ill-conditioned the smaller the step \(\vek{s}_k\) is.

In each iteration, the current approximation \(\tns{C}_k\) moves closer to the integrated derivative \(\widetilde{\tns{C}}_k\), so we cannot expect \(\tns{C}_k\) to approximate \(D^3 f(\vek{x}_k)\) better than \(\widetilde{\tns{C}}_k\).
Of course, the only part of \(\widetilde{\tns{C}}_k\) which is used is its component in the direction of \(\vek{s}_k\), i.e. \(\tns{D}_k / \norm{\vek{s}_k}_2\).
In exact arithmetic the proof of \cref{thm:C-tilde-convergence} also shows that
\begin{equation}
    \frac{\tns{D}_k}{\norm{\vek{s}_k}_2} = \widetilde{\tns{C}}_k\Brack{\vek{s}_k^{\rightarrow}} = D^3 f(\vek{x}_k)\Brack{\vek{s}_k^{\rightarrow}} + \Delta \tns{D}_k \quad \text{with} \quad \norm{\Delta \tns{D}_k}_2 \leq \frac{L}{2} \norm{\vek{s}_k}_2
\end{equation}
where \(\vek{s}_k^{\rightarrow}\) is the normed step \(\vek{s}_k/\norm{\vek{s}_k}_2\), i.e.\ the unit norm vector pointing in the same direction as \(\vek{s}_k\), and \(L\) is the (local) Lipschitz constant of \(D^p f\).

Now, let \(\widehat{\tns{D}}_k\) be the computed \(\tns{D}_k = D^{p-1} f(\vek{x}_{k+1}) - D^{p-1} f(\vek{x}_k)\) under the influence of rounding errors.
As Higham \cite[p. 9]{higham_accuracy_2002} explains, if we subtract two numbers \(\hat{a} = a(1 + \Delta a)\) and \(\hat{b} = b(1 + \Delta b)\) from each other, and assume the relative errors \(\Delta a\) and \(\Delta b\) are bounded by \(\delta\), the absolute error in the result is bounded by
\begin{equation}\label{eqn:cancellation}
    \Abs{-a \Delta a + b \Delta b} \leq \delta \Paren{\abs{a} + \abs{b}}.
\end{equation}
The \(a\) and \(b\) in our case are the entries of \(D^{p-1} f(\vek{x}_{k+1})\) and \(D^{p-1} f(\vek{x}_k)\).
They are computed with the exact formulas, but stored in finite precision, so the best possible error bound \(\delta\) is the machine precision \(\e_{\text{mach}} \approx 10^{-16}\).
This gives
\begin{equation}\label{eqn:two-error-terms}
    \frac{\widehat{\tns{D}}_k}{\norm{\vek{s}_k}_2} = \frac{\tns{D}_k}{\norm{\vek{s}_k}_2} + \Delta \widehat{\tns{D}}_k = D^3 f(\vek{x}_k)\Brack{\vek{s}_k^{\rightarrow}} + \Delta \tns{D}_k + \Delta \widehat{\tns{D}}_k
\end{equation}
where \(\norm{\Delta \widehat{\tns{D}}_k}_F \leq \sqrt{2} \e_{\text{mach}} \Paren{\norm{D^{p-1} f(\vek{x}_{k+1})}_F + \norm{D^{p-1} f(\vek{x}_k)}_F} / \norm{\vek{s}_k}_2\).

\Cref{eqn:two-error-terms} shows that there are two sources of error, one from using the secant equation and one from calculating \(\tns{D}_k\).
The former is proportional to \(\norm{\vek{s}_k}_2\) whereas the latter is proportional to \(1/\norm{\vek{s}_k}_2\).
This leads to the V-shaped graph in \cref{fig:convergence_linearly_converging_min_simple}.
If we assume that \(D^{p-1} f\), \(D^p f\) and \(L\) have roughly the same scale, we can estimate that the lowest possible relative error in \(\widehat{\tns{D}}_k/\norm{\vek{s}_k}_2\) (compared to \(D^3 f(\vek{x}_k)\Brack{\vek{s}_k^{\rightarrow}}\)) is \(\sqrt{\e_{\text{mach}}}\) and is achieved when \(\norm{\vek{s}_k}_2 \approx \sqrt{\e_{\text{mach}}}\).
This analysis is analogous to one for numerical differentiation schemes where the same lower bound is derived, see for example \cite[Section 5.7]{press_numerical_2007}.
One notable exception to the rule occurs when \(D^{p-1} f(\vek{x}_*) = \tns{0}\).
In that case \(\norm{\Delta \widehat{\tns{D}}_k}_F\) stays close to machine precision and the approximations get better and better as \(\vek{s}_k\) converges to \(\vek{0}\).
This is one of the reasons that quasi-Newton methods (\(p = 2\)) work very well for optimization algorithms as they converge to stationary points.

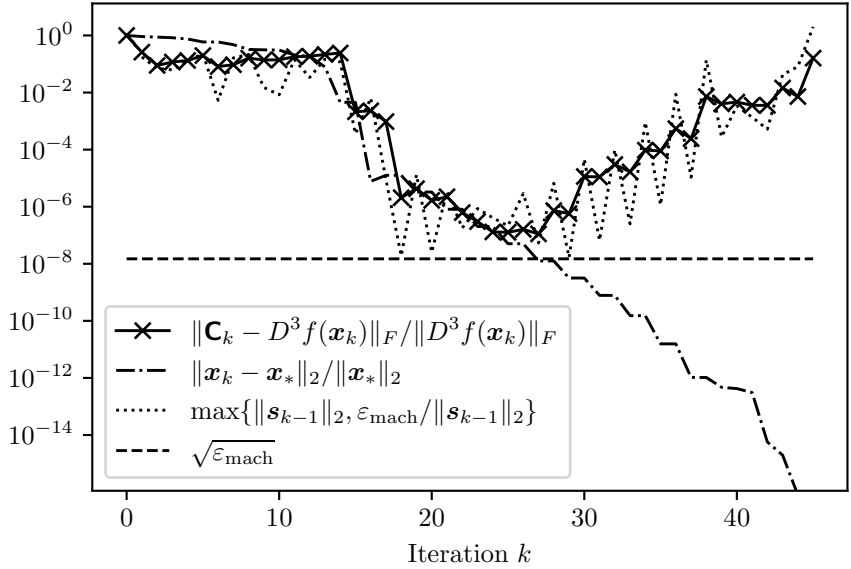
\begin{figure}
    \centering
    \subimport{plots/}{convergence_linearly_converging_min.pgf}
    \caption{Convergence of iterates and approximations for nonlinear CG (extended)}
    \label{fig:convergence_linearly_converging_min}
\end{figure}

In \cref{fig:convergence_linearly_converging_min} one can see that indeed the best relative error achieved is roughly \(\sqrt{\e_{\text{mach}}}\) and that the maximum of \(\norm{\vek{s}_{k-1}}_2\) and \(\e_{\text{mach}}/\norm{\vek{s}_{k-1}}_2\) is a pretty good proxy for a lower bound on the error.

Note that these numerical issues cannot be overcome with a different implementation but are inherent in this approach of extracting third-order information from successive evaluations of second-order derivatives, since computing \(\tns{D}_k\) is an ill-conditioned problem.
A practical way to avoid losing accuracy in the last few iterations would be to employ a heuristic that skips updating \(\tns{C}_k\) when the expected size of errors \(\Delta \widehat{\tns{D}}_k\) exceeds the size of the update or simply when \(\norm{\vek{s}_k}_2\) becomes too small.

\subsection{Convergence in a subspace}

\begin{figure}
    \centering
    \subimport{plots/}{convergence_quadratically_converging_min.pgf}
    \caption{Convergence of iterates and approximations for exact trust region}
    \label{fig:convergence_quadratically_converging_min}
\end{figure}
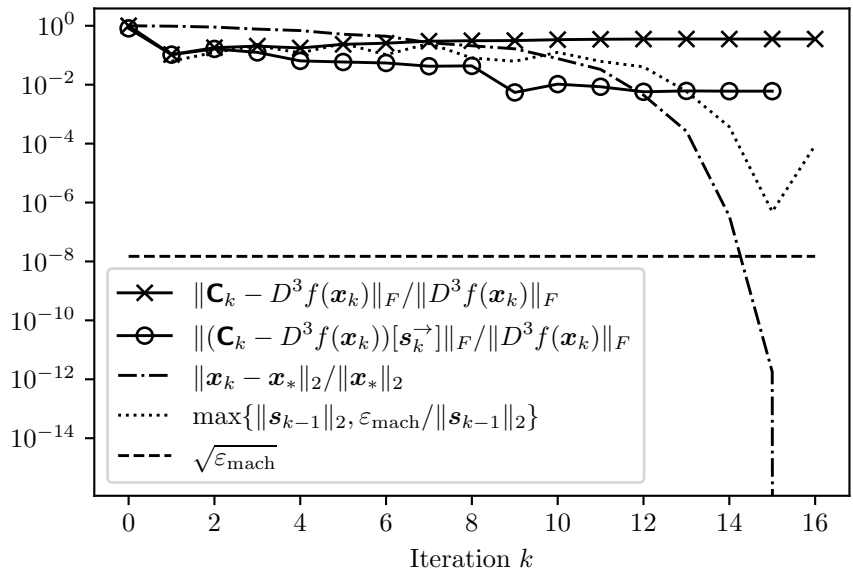

In the second experiment, we used a trust region approach with exact Hessian evaluations to generate the iterates.\footnote{\texttt{scipy.optimize.minimize(method="trust-exact")} in SciPy version 1.9.3, see \cite[pp. 169--200]{conn_trust_2000} for more details. Only the successful iterations were used.}
As predicted by the theory for these methods, we observe local quadratic convergence to the minimizer (and much fewer iterations in general).
The dotted line shows that there are very few iterations in which accurate information about the third derivatives can be obtained, but even then the relative error in \(\tns{C}_k\) is several orders of magnitude larger than the lower bound.

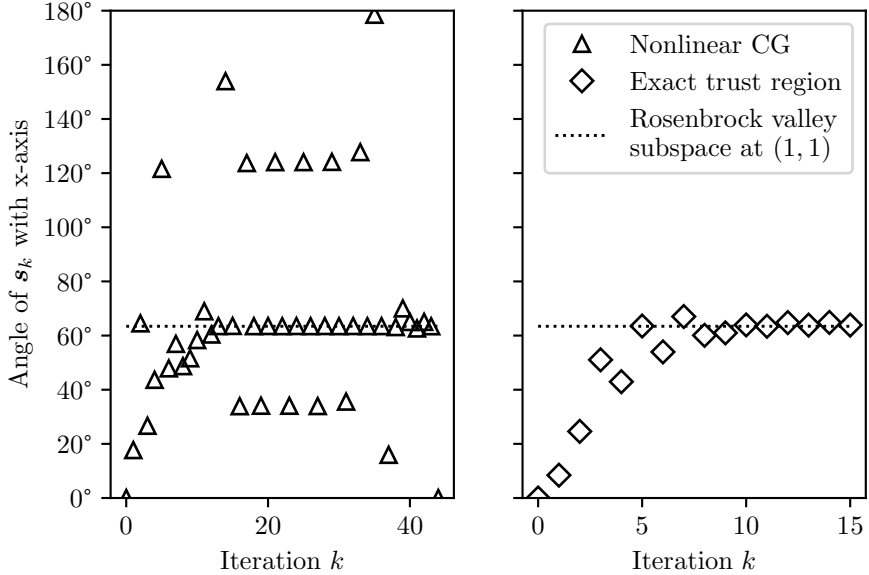
\begin{figure}
    \centering
    \subimport{plots/}{angles_comparison.pgf}
    \caption{Subspaces of the steps for nonlinear CG and exact trust region}
    \label{fig:angles}
\end{figure}

A key difference in the two experiments is how the directions of the steps are distributed.
\Cref{fig:angles} plots the angle of each \(\vek{s}_k\) with the x-axis, normalized between 0\textdegree{} and 180\textdegree{} so that opposite directions coincide.
Whereas the steps generated by the nonlinear CG algorithm cover multiple well-separated directions during the main part of the algorithm, the steps generated by the trust region method tend to fall into a one-dimensional subspace, especially towards the end when convergence happens.
This directly explains why the relative Frobenius error stays high and even increases towards the end in the second experiment:
All the information we can extract from the (averaged) true derivative \(\widetilde{\tns{C}}_k\) is its evaluation in the direction \(\vek{s}_k\) and since most of the steps point in the same direction at the end, the information about the other directions gets more and more outdated.

In addition to the relative Frobenius norm, we also included (a relative version of) the Dennis--Mor\'e measure from \cref{sec:dennis-more} in \cref{fig:convergence_quadratically_converging_min}.
This one measures the error in \(\tns{C}_k\) only in the direction of the step \(\vek{s}_k\).
Now, one might expect that when the approximation is updated in one specific direction and the Dennis--Mor\'e error only measures how good the approximation is in this one direction, the error must track the lower bound derived in the previous subsection quite well.
Indeed, this is the case when the steps all lie \emph{exactly} in one subspace as we could verify with manually generated iterates.
For the exact trust region method however the step directions vary by about 1\textdegree{} in successive iterations at the end, which can be seen in \cref{fig:angles}.
Let \(\vek{s}_k^{\rightarrow} = \lambda_1 \vek{s}_{k-1}^{\rightarrow} + \lambda_2 \vek{u}_k\) where \(\vek{u}_k\) is chosen such that \(\vek{s}_{k-1}^{\rightarrow}\) and \(\vek{u}_k\) form an orthonormal basis of \(\R^2\), then \(\lambda_1^2 + \lambda_2^2 = 1\) and by multilinearity of \(\tns{C}_k\),
\begin{equation}
    \tns{C}_k[\vek{s}_k^{\rightarrow}] = \lambda_1 \tns{C}_k[\vek{s}_{k-1}^{\rightarrow}] + \lambda_2 \tns{C}_k[\vek{u}_k].
\end{equation}
Therefore, in each of the last few iterations the approximation in direction \(\vek{s}_k^{\rightarrow}\) is a linear combination of the very accurate information in direction \(\vek{s}_{k-1}^{\rightarrow}\) and very inaccurate information in direction \(\vek{u}_k\).
In particular, if the angle with the x-axis varies by about 1\textdegree{} we get that \(\lambda_2 \approx \sin(1^\circ) \approx 10^{-2}\).
Combining this with the knowledge that the relative overall error in \(\tns{C}_k\) is on the order of \(1\), we expect that the Dennis--Mor\'e measure will hover around \(10^{-2}\).
This agrees very well with the graph in \cref{fig:convergence_quadratically_converging_min} and shows that it is important for this method to gather accurate derivative information in all directions.

\section{Conclusion}

We have seen in this paper that quasi-Newton updates described as least-change updates admit fairly straightforward generalizations to higher-order derivatives, which we call higher-order secant updates.
These updates have a closed form solution with a certain low-rank structure to it, generalizing the rank-two characterization of regular quasi-Newton updates.
The theoretical results suggest that, as long as the directions of the steps span the space and stay well separated, the generated approximations converge to the true derivative in the limit and under suitably fast convergence of the iterates they even converge (in a subspace) if these assumptions are violated.
This is however not the behaviour we see in experiments, since the problem of computing the difference between Hessian evaluations becomes more and more ill-conditioned as the distance between consecutive iterates becomes smaller.
These numerical limitations lead to a loss in accuracy:
If the Hessians are computed with relative error \(\delta\) we cannot expect the errors in the generated approximations to go below \(\sqrt{\delta}\).
Our experiments show that, as long as the directions of the steps stay well separated, the method indeed generates accurate approximations up to the numerical limit.

Considering these preliminary experiments and the similarities between the convergence results for conventional quasi-Newton updates and our updates, we hope that the generated approximations will be similarly useful for optimization methods.
For example, they could be used inside a third-order tensor method, such as the one investigated by Birgin et al.~\cite{birgin_use_2020},
in order to achieve their favourable complexity results without requiring access to exact third derivatives.
We aim to investigate such methods in a future paper.

\section*{Acknowledgments}

We would like to thank Coralia Cartis, Yuji Nakatsukasa and the anonymous reviewers for their comments on an earlier draft of this paper.
They helped to improve the presentation and readability of our work.

%% file: plots/convergence_linearly_converging_min_simple.pgf
\begingroup%
\makeatletter%
\begin{pgfpicture}%
\pgfpathrectangle{\pgfpointorigin}{\pgfqpoint{4.593087in}{3.152241in}}%
\pgfusepath{use as bounding box, clip}%
\begin{pgfscope}%
\pgfsetbuttcap%
\pgfsetmiterjoin%
\definecolor{currentfill}{rgb}{1.000000,1.000000,1.000000}%
\pgfsetfillcolor{currentfill}%
\pgfsetlinewidth{0.000000pt}%
\definecolor{currentstroke}{rgb}{1.000000,1.000000,1.000000}%
\pgfsetstrokecolor{currentstroke}%
\pgfsetdash{}{0pt}%
\pgfpathmoveto{\pgfqpoint{0.000000in}{0.000000in}}%
\pgfpathlineto{\pgfqpoint{4.593087in}{0.000000in}}%
\pgfpathlineto{\pgfqpoint{4.593087in}{3.152241in}}%
\pgfpathlineto{\pgfqpoint{0.000000in}{3.152241in}}%
\pgfpathlineto{\pgfqpoint{0.000000in}{0.000000in}}%
\pgfpathclose%
\pgfusepath{fill}%
\end{pgfscope}%
\begin{pgfscope}%
\pgfsetbuttcap%
\pgfsetmiterjoin%
\definecolor{currentfill}{rgb}{1.000000,1.000000,1.000000}%
\pgfsetfillcolor{currentfill}%
\pgfsetlinewidth{0.000000pt}%
\definecolor{currentstroke}{rgb}{0.000000,0.000000,0.000000}%
\pgfsetstrokecolor{currentstroke}%
\pgfsetstrokeopacity{0.000000}%
\pgfsetdash{}{0pt}%
\pgfpathmoveto{\pgfqpoint{0.540587in}{0.499691in}}%
\pgfpathlineto{\pgfqpoint{4.493088in}{0.499691in}}%
\pgfpathlineto{\pgfqpoint{4.493088in}{3.052241in}}%
\pgfpathlineto{\pgfqpoint{0.540587in}{3.052241in}}%
\pgfpathlineto{\pgfqpoint{0.540587in}{0.499691in}}%
\pgfpathclose%
\pgfusepath{fill}%
\end{pgfscope}%
\begin{pgfscope}%
\pgfsetbuttcap%
\pgfsetroundjoin%
\definecolor{currentfill}{rgb}{0.000000,0.000000,0.000000}%
\pgfsetfillcolor{currentfill}%
\pgfsetfillopacity{0.000000}%
\pgfsetlinewidth{0.803000pt}%
\definecolor{currentstroke}{rgb}{0.000000,0.000000,0.000000}%
\pgfsetstrokecolor{currentstroke}%
\pgfsetdash{}{0pt}%
\pgfsys@defobject{currentmarker}{\pgfqpoint{0.000000in}{-0.048611in}}{\pgfqpoint{0.000000in}{0.000000in}}{%
\pgfpathmoveto{\pgfqpoint{0.000000in}{0.000000in}}%
\pgfpathlineto{\pgfqpoint{0.000000in}{-0.048611in}}%
\pgfusepath{stroke,fill}%
}%
\begin{pgfscope}%
\pgfsys@transformshift{0.720247in}{0.499691in}%
\pgfsys@useobject{currentmarker}{}%
\end{pgfscope}%
\end{pgfscope}%
\begin{pgfscope}%
\definecolor{textcolor}{rgb}{0.000000,0.000000,0.000000}%
\pgfsetstrokecolor{textcolor}%
\pgfsetfillcolor{textcolor}%
\pgftext[x=0.720247in,y=0.402469in,,top]{\color{textcolor}\rmfamily\fontsize{10.000000}{12.000000}\selectfont \(\displaystyle {0}\)}%
\end{pgfscope}%
\begin{pgfscope}%
\pgfsetbuttcap%
\pgfsetroundjoin%
\definecolor{currentfill}{rgb}{0.000000,0.000000,0.000000}%
\pgfsetfillcolor{currentfill}%
\pgfsetfillopacity{0.000000}%
\pgfsetlinewidth{0.803000pt}%
\definecolor{currentstroke}{rgb}{0.000000,0.000000,0.000000}%
\pgfsetstrokecolor{currentstroke}%
\pgfsetdash{}{0pt}%
\pgfsys@defobject{currentmarker}{\pgfqpoint{0.000000in}{-0.048611in}}{\pgfqpoint{0.000000in}{0.000000in}}{%
\pgfpathmoveto{\pgfqpoint{0.000000in}{0.000000in}}%
\pgfpathlineto{\pgfqpoint{0.000000in}{-0.048611in}}%
\pgfusepath{stroke,fill}%
}%
\begin{pgfscope}%
\pgfsys@transformshift{1.518731in}{0.499691in}%
\pgfsys@useobject{currentmarker}{}%
\end{pgfscope}%
\end{pgfscope}%
\begin{pgfscope}%
\definecolor{textcolor}{rgb}{0.000000,0.000000,0.000000}%
\pgfsetstrokecolor{textcolor}%
\pgfsetfillcolor{textcolor}%
\pgftext[x=1.518731in,y=0.402469in,,top]{\color{textcolor}\rmfamily\fontsize{10.000000}{12.000000}\selectfont \(\displaystyle {10}\)}%
\end{pgfscope}%
\begin{pgfscope}%
\pgfsetbuttcap%
\pgfsetroundjoin%
\definecolor{currentfill}{rgb}{0.000000,0.000000,0.000000}%
\pgfsetfillcolor{currentfill}%
\pgfsetfillopacity{0.000000}%
\pgfsetlinewidth{0.803000pt}%
\definecolor{currentstroke}{rgb}{0.000000,0.000000,0.000000}%
\pgfsetstrokecolor{currentstroke}%
\pgfsetdash{}{0pt}%
\pgfsys@defobject{currentmarker}{\pgfqpoint{0.000000in}{-0.048611in}}{\pgfqpoint{0.000000in}{0.000000in}}{%
\pgfpathmoveto{\pgfqpoint{0.000000in}{0.000000in}}%
\pgfpathlineto{\pgfqpoint{0.000000in}{-0.048611in}}%
\pgfusepath{stroke,fill}%
}%
\begin{pgfscope}%
\pgfsys@transformshift{2.317216in}{0.499691in}%
\pgfsys@useobject{currentmarker}{}%
\end{pgfscope}%
\end{pgfscope}%
\begin{pgfscope}%
\definecolor{textcolor}{rgb}{0.000000,0.000000,0.000000}%
\pgfsetstrokecolor{textcolor}%
\pgfsetfillcolor{textcolor}%
\pgftext[x=2.317216in,y=0.402469in,,top]{\color{textcolor}\rmfamily\fontsize{10.000000}{12.000000}\selectfont \(\displaystyle {20}\)}%
\end{pgfscope}%
\begin{pgfscope}%
\pgfsetbuttcap%
\pgfsetroundjoin%
\definecolor{currentfill}{rgb}{0.000000,0.000000,0.000000}%
\pgfsetfillcolor{currentfill}%
\pgfsetfillopacity{0.000000}%
\pgfsetlinewidth{0.803000pt}%
\definecolor{currentstroke}{rgb}{0.000000,0.000000,0.000000}%
\pgfsetstrokecolor{currentstroke}%
\pgfsetdash{}{0pt}%
\pgfsys@defobject{currentmarker}{\pgfqpoint{0.000000in}{-0.048611in}}{\pgfqpoint{0.000000in}{0.000000in}}{%
\pgfpathmoveto{\pgfqpoint{0.000000in}{0.000000in}}%
\pgfpathlineto{\pgfqpoint{0.000000in}{-0.048611in}}%
\pgfusepath{stroke,fill}%
}%
\begin{pgfscope}%
\pgfsys@transformshift{3.115701in}{0.499691in}%
\pgfsys@useobject{currentmarker}{}%
\end{pgfscope}%
\end{pgfscope}%
\begin{pgfscope}%
\definecolor{textcolor}{rgb}{0.000000,0.000000,0.000000}%
\pgfsetstrokecolor{textcolor}%
\pgfsetfillcolor{textcolor}%
\pgftext[x=3.115701in,y=0.402469in,,top]{\color{textcolor}\rmfamily\fontsize{10.000000}{12.000000}\selectfont \(\displaystyle {30}\)}%
\end{pgfscope}%
\begin{pgfscope}%
\pgfsetbuttcap%
\pgfsetroundjoin%
\definecolor{currentfill}{rgb}{0.000000,0.000000,0.000000}%
\pgfsetfillcolor{currentfill}%
\pgfsetfillopacity{0.000000}%
\pgfsetlinewidth{0.803000pt}%
\definecolor{currentstroke}{rgb}{0.000000,0.000000,0.000000}%
\pgfsetstrokecolor{currentstroke}%
\pgfsetdash{}{0pt}%
\pgfsys@defobject{currentmarker}{\pgfqpoint{0.000000in}{-0.048611in}}{\pgfqpoint{0.000000in}{0.000000in}}{%
\pgfpathmoveto{\pgfqpoint{0.000000in}{0.000000in}}%
\pgfpathlineto{\pgfqpoint{0.000000in}{-0.048611in}}%
\pgfusepath{stroke,fill}%
}%
\begin{pgfscope}%
\pgfsys@transformshift{3.914186in}{0.499691in}%
\pgfsys@useobject{currentmarker}{}%
\end{pgfscope}%
\end{pgfscope}%
\begin{pgfscope}%
\definecolor{textcolor}{rgb}{0.000000,0.000000,0.000000}%
\pgfsetstrokecolor{textcolor}%
\pgfsetfillcolor{textcolor}%
\pgftext[x=3.914186in,y=0.402469in,,top]{\color{textcolor}\rmfamily\fontsize{10.000000}{12.000000}\selectfont \(\displaystyle {40}\)}%
\end{pgfscope}%
\begin{pgfscope}%
\definecolor{textcolor}{rgb}{0.000000,0.000000,0.000000}%
\pgfsetstrokecolor{textcolor}%
\pgfsetfillcolor{textcolor}%
\pgftext[x=2.516837in,y=0.223457in,,top]{\color{textcolor}\rmfamily\fontsize{10.000000}{12.000000}\selectfont Iteration \(\displaystyle k\)}%
\end{pgfscope}%
\begin{pgfscope}%
\pgfsetbuttcap%
\pgfsetroundjoin%
\definecolor{currentfill}{rgb}{0.000000,0.000000,0.000000}%
\pgfsetfillcolor{currentfill}%
\pgfsetfillopacity{0.000000}%
\pgfsetlinewidth{0.803000pt}%
\definecolor{currentstroke}{rgb}{0.000000,0.000000,0.000000}%
\pgfsetstrokecolor{currentstroke}%
\pgfsetdash{}{0pt}%
\pgfsys@defobject{currentmarker}{\pgfqpoint{-0.048611in}{0.000000in}}{\pgfqpoint{-0.000000in}{0.000000in}}{%
\pgfpathmoveto{\pgfqpoint{-0.000000in}{0.000000in}}%
\pgfpathlineto{\pgfqpoint{-0.048611in}{0.000000in}}%
\pgfusepath{stroke,fill}%
}%
\begin{pgfscope}%
\pgfsys@transformshift{0.540587in}{0.797378in}%
\pgfsys@useobject{currentmarker}{}%
\end{pgfscope}%
\end{pgfscope}%
\begin{pgfscope}%
\definecolor{textcolor}{rgb}{0.000000,0.000000,0.000000}%
\pgfsetstrokecolor{textcolor}%
\pgfsetfillcolor{textcolor}%
\pgftext[x=0.100000in, y=0.749153in, left, base]{\color{textcolor}\rmfamily\fontsize{10.000000}{12.000000}\selectfont \(\displaystyle {10^{-14}}\)}%
\end{pgfscope}%
\begin{pgfscope}%
\pgfsetbuttcap%
\pgfsetroundjoin%
\definecolor{currentfill}{rgb}{0.000000,0.000000,0.000000}%
\pgfsetfillcolor{currentfill}%
\pgfsetfillopacity{0.000000}%
\pgfsetlinewidth{0.803000pt}%
\definecolor{currentstroke}{rgb}{0.000000,0.000000,0.000000}%
\pgfsetstrokecolor{currentstroke}%
\pgfsetdash{}{0pt}%
\pgfsys@defobject{currentmarker}{\pgfqpoint{-0.048611in}{0.000000in}}{\pgfqpoint{-0.000000in}{0.000000in}}{%
\pgfpathmoveto{\pgfqpoint{-0.000000in}{0.000000in}}%
\pgfpathlineto{\pgfqpoint{-0.048611in}{0.000000in}}%
\pgfusepath{stroke,fill}%
}%
\begin{pgfscope}%
\pgfsys@transformshift{0.540587in}{1.101981in}%
\pgfsys@useobject{currentmarker}{}%
\end{pgfscope}%
\end{pgfscope}%
\begin{pgfscope}%
\definecolor{textcolor}{rgb}{0.000000,0.000000,0.000000}%
\pgfsetstrokecolor{textcolor}%
\pgfsetfillcolor{textcolor}%
\pgftext[x=0.100000in, y=1.053756in, left, base]{\color{textcolor}\rmfamily\fontsize{10.000000}{12.000000}\selectfont \(\displaystyle {10^{-12}}\)}%
\end{pgfscope}%
\begin{pgfscope}%
\pgfsetbuttcap%
\pgfsetroundjoin%
\definecolor{currentfill}{rgb}{0.000000,0.000000,0.000000}%
\pgfsetfillcolor{currentfill}%
\pgfsetfillopacity{0.000000}%
\pgfsetlinewidth{0.803000pt}%
\definecolor{currentstroke}{rgb}{0.000000,0.000000,0.000000}%
\pgfsetstrokecolor{currentstroke}%
\pgfsetdash{}{0pt}%
\pgfsys@defobject{currentmarker}{\pgfqpoint{-0.048611in}{0.000000in}}{\pgfqpoint{-0.000000in}{0.000000in}}{%
\pgfpathmoveto{\pgfqpoint{-0.000000in}{0.000000in}}%
\pgfpathlineto{\pgfqpoint{-0.048611in}{0.000000in}}%
\pgfusepath{stroke,fill}%
}%
\begin{pgfscope}%
\pgfsys@transformshift{0.540587in}{1.406584in}%
\pgfsys@useobject{currentmarker}{}%
\end{pgfscope}%
\end{pgfscope}%
\begin{pgfscope}%
\definecolor{textcolor}{rgb}{0.000000,0.000000,0.000000}%
\pgfsetstrokecolor{textcolor}%
\pgfsetfillcolor{textcolor}%
\pgftext[x=0.100000in, y=1.358359in, left, base]{\color{textcolor}\rmfamily\fontsize{10.000000}{12.000000}\selectfont \(\displaystyle {10^{-10}}\)}%
\end{pgfscope}%
\begin{pgfscope}%
\pgfsetbuttcap%
\pgfsetroundjoin%
\definecolor{currentfill}{rgb}{0.000000,0.000000,0.000000}%
\pgfsetfillcolor{currentfill}%
\pgfsetfillopacity{0.000000}%
\pgfsetlinewidth{0.803000pt}%
\definecolor{currentstroke}{rgb}{0.000000,0.000000,0.000000}%
\pgfsetstrokecolor{currentstroke}%
\pgfsetdash{}{0pt}%
\pgfsys@defobject{currentmarker}{\pgfqpoint{-0.048611in}{0.000000in}}{\pgfqpoint{-0.000000in}{0.000000in}}{%
\pgfpathmoveto{\pgfqpoint{-0.000000in}{0.000000in}}%
\pgfpathlineto{\pgfqpoint{-0.048611in}{0.000000in}}%
\pgfusepath{stroke,fill}%
}%
\begin{pgfscope}%
\pgfsys@transformshift{0.540587in}{1.711187in}%
\pgfsys@useobject{currentmarker}{}%
\end{pgfscope}%
\end{pgfscope}%
\begin{pgfscope}%
\definecolor{textcolor}{rgb}{0.000000,0.000000,0.000000}%
\pgfsetstrokecolor{textcolor}%
\pgfsetfillcolor{textcolor}%
\pgftext[x=0.155363in, y=1.662962in, left, base]{\color{textcolor}\rmfamily\fontsize{10.000000}{12.000000}\selectfont \(\displaystyle {10^{-8}}\)}%
\end{pgfscope}%
\begin{pgfscope}%
\pgfsetbuttcap%
\pgfsetroundjoin%
\definecolor{currentfill}{rgb}{0.000000,0.000000,0.000000}%
\pgfsetfillcolor{currentfill}%
\pgfsetfillopacity{0.000000}%
\pgfsetlinewidth{0.803000pt}%
\definecolor{currentstroke}{rgb}{0.000000,0.000000,0.000000}%
\pgfsetstrokecolor{currentstroke}%
\pgfsetdash{}{0pt}%
\pgfsys@defobject{currentmarker}{\pgfqpoint{-0.048611in}{0.000000in}}{\pgfqpoint{-0.000000in}{0.000000in}}{%
\pgfpathmoveto{\pgfqpoint{-0.000000in}{0.000000in}}%
\pgfpathlineto{\pgfqpoint{-0.048611in}{0.000000in}}%
\pgfusepath{stroke,fill}%
}%
\begin{pgfscope}%
\pgfsys@transformshift{0.540587in}{2.015790in}%
\pgfsys@useobject{currentmarker}{}%
\end{pgfscope}%
\end{pgfscope}%
\begin{pgfscope}%
\definecolor{textcolor}{rgb}{0.000000,0.000000,0.000000}%
\pgfsetstrokecolor{textcolor}%
\pgfsetfillcolor{textcolor}%
\pgftext[x=0.155363in, y=1.967565in, left, base]{\color{textcolor}\rmfamily\fontsize{10.000000}{12.000000}\selectfont \(\displaystyle {10^{-6}}\)}%
\end{pgfscope}%
\begin{pgfscope}%
\pgfsetbuttcap%
\pgfsetroundjoin%
\definecolor{currentfill}{rgb}{0.000000,0.000000,0.000000}%
\pgfsetfillcolor{currentfill}%
\pgfsetfillopacity{0.000000}%
\pgfsetlinewidth{0.803000pt}%
\definecolor{currentstroke}{rgb}{0.000000,0.000000,0.000000}%
\pgfsetstrokecolor{currentstroke}%
\pgfsetdash{}{0pt}%
\pgfsys@defobject{currentmarker}{\pgfqpoint{-0.048611in}{0.000000in}}{\pgfqpoint{-0.000000in}{0.000000in}}{%
\pgfpathmoveto{\pgfqpoint{-0.000000in}{0.000000in}}%
\pgfpathlineto{\pgfqpoint{-0.048611in}{0.000000in}}%
\pgfusepath{stroke,fill}%
}%
\begin{pgfscope}%
\pgfsys@transformshift{0.540587in}{2.320393in}%
\pgfsys@useobject{currentmarker}{}%
\end{pgfscope}%
\end{pgfscope}%
\begin{pgfscope}%
\definecolor{textcolor}{rgb}{0.000000,0.000000,0.000000}%
\pgfsetstrokecolor{textcolor}%
\pgfsetfillcolor{textcolor}%
\pgftext[x=0.155363in, y=2.272168in, left, base]{\color{textcolor}\rmfamily\fontsize{10.000000}{12.000000}\selectfont \(\displaystyle {10^{-4}}\)}%
\end{pgfscope}%
\begin{pgfscope}%
\pgfsetbuttcap%
\pgfsetroundjoin%
\definecolor{currentfill}{rgb}{0.000000,0.000000,0.000000}%
\pgfsetfillcolor{currentfill}%
\pgfsetfillopacity{0.000000}%
\pgfsetlinewidth{0.803000pt}%
\definecolor{currentstroke}{rgb}{0.000000,0.000000,0.000000}%
\pgfsetstrokecolor{currentstroke}%
\pgfsetdash{}{0pt}%
\pgfsys@defobject{currentmarker}{\pgfqpoint{-0.048611in}{0.000000in}}{\pgfqpoint{-0.000000in}{0.000000in}}{%
\pgfpathmoveto{\pgfqpoint{-0.000000in}{0.000000in}}%
\pgfpathlineto{\pgfqpoint{-0.048611in}{0.000000in}}%
\pgfusepath{stroke,fill}%
}%
\begin{pgfscope}%
\pgfsys@transformshift{0.540587in}{2.624996in}%
\pgfsys@useobject{currentmarker}{}%
\end{pgfscope}%
\end{pgfscope}%
\begin{pgfscope}%
\definecolor{textcolor}{rgb}{0.000000,0.000000,0.000000}%
\pgfsetstrokecolor{textcolor}%
\pgfsetfillcolor{textcolor}%
\pgftext[x=0.155363in, y=2.576771in, left, base]{\color{textcolor}\rmfamily\fontsize{10.000000}{12.000000}\selectfont \(\displaystyle {10^{-2}}\)}%
\end{pgfscope}%
\begin{pgfscope}%
\pgfsetbuttcap%
\pgfsetroundjoin%
\definecolor{currentfill}{rgb}{0.000000,0.000000,0.000000}%
\pgfsetfillcolor{currentfill}%
\pgfsetfillopacity{0.000000}%
\pgfsetlinewidth{0.803000pt}%
\definecolor{currentstroke}{rgb}{0.000000,0.000000,0.000000}%
\pgfsetstrokecolor{currentstroke}%
\pgfsetdash{}{0pt}%
\pgfsys@defobject{currentmarker}{\pgfqpoint{-0.048611in}{0.000000in}}{\pgfqpoint{-0.000000in}{0.000000in}}{%
\pgfpathmoveto{\pgfqpoint{-0.000000in}{0.000000in}}%
\pgfpathlineto{\pgfqpoint{-0.048611in}{0.000000in}}%
\pgfusepath{stroke,fill}%
}%
\begin{pgfscope}%
\pgfsys@transformshift{0.540587in}{2.929600in}%
\pgfsys@useobject{currentmarker}{}%
\end{pgfscope}%
\end{pgfscope}%
\begin{pgfscope}%
\definecolor{textcolor}{rgb}{0.000000,0.000000,0.000000}%
\pgfsetstrokecolor{textcolor}%
\pgfsetfillcolor{textcolor}%
\pgftext[x=0.242169in, y=2.881374in, left, base]{\color{textcolor}\rmfamily\fontsize{10.000000}{12.000000}\selectfont \(\displaystyle {10^{0}}\)}%
\end{pgfscope}%
\begin{pgfscope}%
\pgfpathrectangle{\pgfqpoint{0.540587in}{0.499691in}}{\pgfqpoint{3.952500in}{2.552550in}}%
\pgfusepath{clip}%
\pgfsetrectcap%
\pgfsetroundjoin%
\pgfsetlinewidth{1.104125pt}%
\definecolor{currentstroke}{rgb}{0.000000,0.000000,0.000000}%
\pgfsetstrokecolor{currentstroke}%
\pgfsetdash{}{0pt}%
\pgfpathmoveto{\pgfqpoint{0.720247in}{2.929600in}}%
\pgfpathlineto{\pgfqpoint{0.800095in}{2.840666in}}%
\pgfpathlineto{\pgfqpoint{0.879944in}{2.769655in}}%
\pgfpathlineto{\pgfqpoint{0.959792in}{2.787699in}}%
\pgfpathlineto{\pgfqpoint{1.039641in}{2.796051in}}%
\pgfpathlineto{\pgfqpoint{1.119489in}{2.822332in}}%
\pgfpathlineto{\pgfqpoint{1.199337in}{2.763636in}}%
\pgfpathlineto{\pgfqpoint{1.279186in}{2.773643in}}%
\pgfpathlineto{\pgfqpoint{1.359034in}{2.808530in}}%
\pgfpathlineto{\pgfqpoint{1.438883in}{2.798782in}}%
\pgfpathlineto{\pgfqpoint{1.518731in}{2.799688in}}%
\pgfpathlineto{\pgfqpoint{1.598580in}{2.816330in}}%
\pgfpathlineto{\pgfqpoint{1.678428in}{2.818468in}}%
\pgfpathlineto{\pgfqpoint{1.758277in}{2.828509in}}%
\pgfpathlineto{\pgfqpoint{1.838125in}{2.836001in}}%
\pgfpathlineto{\pgfqpoint{1.917974in}{2.522547in}}%
\pgfpathlineto{\pgfqpoint{1.997822in}{2.527549in}}%
\pgfpathlineto{\pgfqpoint{2.077671in}{2.469980in}}%
\pgfpathlineto{\pgfqpoint{2.157519in}{2.063711in}}%
\pgfpathlineto{\pgfqpoint{2.237368in}{2.111829in}}%
\pgfpathlineto{\pgfqpoint{2.317216in}{2.049573in}}%
\pgfpathlineto{\pgfqpoint{2.397065in}{2.068900in}}%
\pgfpathlineto{\pgfqpoint{2.476913in}{1.984499in}}%
\pgfpathlineto{\pgfqpoint{2.556762in}{1.936322in}}%
\pgfpathlineto{\pgfqpoint{2.636610in}{1.881096in}}%
\pgfpathlineto{\pgfqpoint{2.716459in}{1.879050in}}%
\pgfpathlineto{\pgfqpoint{2.796307in}{1.894727in}}%
\pgfpathlineto{\pgfqpoint{2.876156in}{1.870654in}}%
\pgfpathlineto{\pgfqpoint{2.956004in}{1.994749in}}%
\pgfpathlineto{\pgfqpoint{3.035853in}{1.978936in}}%
\pgfpathlineto{\pgfqpoint{3.115701in}{2.176765in}}%
\pgfpathlineto{\pgfqpoint{3.195550in}{2.174821in}}%
\pgfpathlineto{\pgfqpoint{3.275398in}{2.241178in}}%
\pgfpathlineto{\pgfqpoint{3.355247in}{2.201395in}}%
\pgfpathlineto{\pgfqpoint{3.435095in}{2.317407in}}%
\pgfpathlineto{\pgfqpoint{3.514944in}{2.313621in}}%
\pgfpathlineto{\pgfqpoint{3.594792in}{2.432503in}}%
\pgfpathlineto{\pgfqpoint{3.674641in}{2.379267in}}%
\pgfpathlineto{\pgfqpoint{3.754489in}{2.603585in}}%
\pgfpathlineto{\pgfqpoint{3.834338in}{2.564686in}}%
\pgfpathlineto{\pgfqpoint{3.914186in}{2.574364in}}%
\pgfpathlineto{\pgfqpoint{3.994034in}{2.556808in}}%
\pgfpathlineto{\pgfqpoint{4.073883in}{2.556632in}}%
\pgfpathlineto{\pgfqpoint{4.153731in}{2.649274in}}%
\pgfpathlineto{\pgfqpoint{4.233580in}{2.603720in}}%
\pgfpathlineto{\pgfqpoint{4.313428in}{2.808758in}}%
\pgfusepath{stroke}%
\end{pgfscope}%
\begin{pgfscope}%
\pgfpathrectangle{\pgfqpoint{0.540587in}{0.499691in}}{\pgfqpoint{3.952500in}{2.552550in}}%
\pgfusepath{clip}%
\pgfsetbuttcap%
\pgfsetroundjoin%
\definecolor{currentfill}{rgb}{0.000000,0.000000,0.000000}%
\pgfsetfillcolor{currentfill}%
\pgfsetfillopacity{0.000000}%
\pgfsetlinewidth{1.003750pt}%
\definecolor{currentstroke}{rgb}{0.000000,0.000000,0.000000}%
\pgfsetstrokecolor{currentstroke}%
\pgfsetdash{}{0pt}%
\pgfsys@defobject{currentmarker}{\pgfqpoint{-0.041667in}{-0.041667in}}{\pgfqpoint{0.041667in}{0.041667in}}{%
\pgfpathmoveto{\pgfqpoint{-0.041667in}{-0.041667in}}%
\pgfpathlineto{\pgfqpoint{0.041667in}{0.041667in}}%
\pgfpathmoveto{\pgfqpoint{-0.041667in}{0.041667in}}%
\pgfpathlineto{\pgfqpoint{0.041667in}{-0.041667in}}%
\pgfusepath{stroke,fill}%
}%
\begin{pgfscope}%
\pgfsys@transformshift{0.720247in}{2.929600in}%
\pgfsys@useobject{currentmarker}{}%
\end{pgfscope}%
\begin{pgfscope}%
\pgfsys@transformshift{0.800095in}{2.840666in}%
\pgfsys@useobject{currentmarker}{}%
\end{pgfscope}%
\begin{pgfscope}%
\pgfsys@transformshift{0.879944in}{2.769655in}%
\pgfsys@useobject{currentmarker}{}%
\end{pgfscope}%
\begin{pgfscope}%
\pgfsys@transformshift{0.959792in}{2.787699in}%
\pgfsys@useobject{currentmarker}{}%
\end{pgfscope}%
\begin{pgfscope}%
\pgfsys@transformshift{1.039641in}{2.796051in}%
\pgfsys@useobject{currentmarker}{}%
\end{pgfscope}%
\begin{pgfscope}%
\pgfsys@transformshift{1.119489in}{2.822332in}%
\pgfsys@useobject{currentmarker}{}%
\end{pgfscope}%
\begin{pgfscope}%
\pgfsys@transformshift{1.199337in}{2.763636in}%
\pgfsys@useobject{currentmarker}{}%
\end{pgfscope}%
\begin{pgfscope}%
\pgfsys@transformshift{1.279186in}{2.773643in}%
\pgfsys@useobject{currentmarker}{}%
\end{pgfscope}%
\begin{pgfscope}%
\pgfsys@transformshift{1.359034in}{2.808530in}%
\pgfsys@useobject{currentmarker}{}%
\end{pgfscope}%
\begin{pgfscope}%
\pgfsys@transformshift{1.438883in}{2.798782in}%
\pgfsys@useobject{currentmarker}{}%
\end{pgfscope}%
\begin{pgfscope}%
\pgfsys@transformshift{1.518731in}{2.799688in}%
\pgfsys@useobject{currentmarker}{}%
\end{pgfscope}%
\begin{pgfscope}%
\pgfsys@transformshift{1.598580in}{2.816330in}%
\pgfsys@useobject{currentmarker}{}%
\end{pgfscope}%
\begin{pgfscope}%
\pgfsys@transformshift{1.678428in}{2.818468in}%
\pgfsys@useobject{currentmarker}{}%
\end{pgfscope}%
\begin{pgfscope}%
\pgfsys@transformshift{1.758277in}{2.828509in}%
\pgfsys@useobject{currentmarker}{}%
\end{pgfscope}%
\begin{pgfscope}%
\pgfsys@transformshift{1.838125in}{2.836001in}%
\pgfsys@useobject{currentmarker}{}%
\end{pgfscope}%
\begin{pgfscope}%
\pgfsys@transformshift{1.917974in}{2.522547in}%
\pgfsys@useobject{currentmarker}{}%
\end{pgfscope}%
\begin{pgfscope}%
\pgfsys@transformshift{1.997822in}{2.527549in}%
\pgfsys@useobject{currentmarker}{}%
\end{pgfscope}%
\begin{pgfscope}%
\pgfsys@transformshift{2.077671in}{2.469980in}%
\pgfsys@useobject{currentmarker}{}%
\end{pgfscope}%
\begin{pgfscope}%
\pgfsys@transformshift{2.157519in}{2.063711in}%
\pgfsys@useobject{currentmarker}{}%
\end{pgfscope}%
\begin{pgfscope}%
\pgfsys@transformshift{2.237368in}{2.111829in}%
\pgfsys@useobject{currentmarker}{}%
\end{pgfscope}%
\begin{pgfscope}%
\pgfsys@transformshift{2.317216in}{2.049573in}%
\pgfsys@useobject{currentmarker}{}%
\end{pgfscope}%
\begin{pgfscope}%
\pgfsys@transformshift{2.397065in}{2.068900in}%
\pgfsys@useobject{currentmarker}{}%
\end{pgfscope}%
\begin{pgfscope}%
\pgfsys@transformshift{2.476913in}{1.984499in}%
\pgfsys@useobject{currentmarker}{}%
\end{pgfscope}%
\begin{pgfscope}%
\pgfsys@transformshift{2.556762in}{1.936322in}%
\pgfsys@useobject{currentmarker}{}%
\end{pgfscope}%
\begin{pgfscope}%
\pgfsys@transformshift{2.636610in}{1.881096in}%
\pgfsys@useobject{currentmarker}{}%
\end{pgfscope}%
\begin{pgfscope}%
\pgfsys@transformshift{2.716459in}{1.879050in}%
\pgfsys@useobject{currentmarker}{}%
\end{pgfscope}%
\begin{pgfscope}%
\pgfsys@transformshift{2.796307in}{1.894727in}%
\pgfsys@useobject{currentmarker}{}%
\end{pgfscope}%
\begin{pgfscope}%
\pgfsys@transformshift{2.876156in}{1.870654in}%
\pgfsys@useobject{currentmarker}{}%
\end{pgfscope}%
\begin{pgfscope}%
\pgfsys@transformshift{2.956004in}{1.994749in}%
\pgfsys@useobject{currentmarker}{}%
\end{pgfscope}%
\begin{pgfscope}%
\pgfsys@transformshift{3.035853in}{1.978936in}%
\pgfsys@useobject{currentmarker}{}%
\end{pgfscope}%
\begin{pgfscope}%
\pgfsys@transformshift{3.115701in}{2.176765in}%
\pgfsys@useobject{currentmarker}{}%
\end{pgfscope}%
\begin{pgfscope}%
\pgfsys@transformshift{3.195550in}{2.174821in}%
\pgfsys@useobject{currentmarker}{}%
\end{pgfscope}%
\begin{pgfscope}%
\pgfsys@transformshift{3.275398in}{2.241178in}%
\pgfsys@useobject{currentmarker}{}%
\end{pgfscope}%
\begin{pgfscope}%
\pgfsys@transformshift{3.355247in}{2.201395in}%
\pgfsys@useobject{currentmarker}{}%
\end{pgfscope}%
\begin{pgfscope}%
\pgfsys@transformshift{3.435095in}{2.317407in}%
\pgfsys@useobject{currentmarker}{}%
\end{pgfscope}%
\begin{pgfscope}%
\pgfsys@transformshift{3.514944in}{2.313621in}%
\pgfsys@useobject{currentmarker}{}%
\end{pgfscope}%
\begin{pgfscope}%
\pgfsys@transformshift{3.594792in}{2.432503in}%
\pgfsys@useobject{currentmarker}{}%
\end{pgfscope}%
\begin{pgfscope}%
\pgfsys@transformshift{3.674641in}{2.379267in}%
\pgfsys@useobject{currentmarker}{}%
\end{pgfscope}%
\begin{pgfscope}%
\pgfsys@transformshift{3.754489in}{2.603585in}%
\pgfsys@useobject{currentmarker}{}%
\end{pgfscope}%
\begin{pgfscope}%
\pgfsys@transformshift{3.834338in}{2.564686in}%
\pgfsys@useobject{currentmarker}{}%
\end{pgfscope}%
\begin{pgfscope}%
\pgfsys@transformshift{3.914186in}{2.574364in}%
\pgfsys@useobject{currentmarker}{}%
\end{pgfscope}%
\begin{pgfscope}%
\pgfsys@transformshift{3.994034in}{2.556808in}%
\pgfsys@useobject{currentmarker}{}%
\end{pgfscope}%
\begin{pgfscope}%
\pgfsys@transformshift{4.073883in}{2.556632in}%
\pgfsys@useobject{currentmarker}{}%
\end{pgfscope}%
\begin{pgfscope}%
\pgfsys@transformshift{4.153731in}{2.649274in}%
\pgfsys@useobject{currentmarker}{}%
\end{pgfscope}%
\begin{pgfscope}%
\pgfsys@transformshift{4.233580in}{2.603720in}%
\pgfsys@useobject{currentmarker}{}%
\end{pgfscope}%
\begin{pgfscope}%
\pgfsys@transformshift{4.313428in}{2.808758in}%
\pgfsys@useobject{currentmarker}{}%
\end{pgfscope}%
\end{pgfscope}%
\begin{pgfscope}%
\pgfpathrectangle{\pgfqpoint{0.540587in}{0.499691in}}{\pgfqpoint{3.952500in}{2.552550in}}%
\pgfusepath{clip}%
\pgfsetbuttcap%
\pgfsetroundjoin%
\pgfsetlinewidth{1.104125pt}%
\definecolor{currentstroke}{rgb}{0.000000,0.000000,0.000000}%
\pgfsetstrokecolor{currentstroke}%
\pgfsetdash{{7.040000pt}{1.760000pt}{1.100000pt}{1.760000pt}}{0.000000pt}%
\pgfpathmoveto{\pgfqpoint{0.720247in}{2.929600in}}%
\pgfpathlineto{\pgfqpoint{0.800095in}{2.923802in}}%
\pgfpathlineto{\pgfqpoint{0.879944in}{2.920759in}}%
\pgfpathlineto{\pgfqpoint{0.959792in}{2.917424in}}%
\pgfpathlineto{\pgfqpoint{1.039641in}{2.910661in}}%
\pgfpathlineto{\pgfqpoint{1.119489in}{2.894944in}}%
\pgfpathlineto{\pgfqpoint{1.199337in}{2.894761in}}%
\pgfpathlineto{\pgfqpoint{1.279186in}{2.880505in}}%
\pgfpathlineto{\pgfqpoint{1.359034in}{2.855053in}}%
\pgfpathlineto{\pgfqpoint{1.438883in}{2.852896in}}%
\pgfpathlineto{\pgfqpoint{1.518731in}{2.851656in}}%
\pgfpathlineto{\pgfqpoint{1.598580in}{2.825642in}}%
\pgfpathlineto{\pgfqpoint{1.678428in}{2.817621in}}%
\pgfpathlineto{\pgfqpoint{1.758277in}{2.761954in}}%
\pgfpathlineto{\pgfqpoint{1.838125in}{2.572776in}}%
\pgfpathlineto{\pgfqpoint{1.917974in}{2.572643in}}%
\pgfpathlineto{\pgfqpoint{1.997822in}{2.151097in}}%
\pgfpathlineto{\pgfqpoint{2.077671in}{2.183292in}}%
\pgfpathlineto{\pgfqpoint{2.157519in}{2.183257in}}%
\pgfpathlineto{\pgfqpoint{2.237368in}{2.094166in}}%
\pgfpathlineto{\pgfqpoint{2.317216in}{2.094057in}}%
\pgfpathlineto{\pgfqpoint{2.397065in}{2.002056in}}%
\pgfpathlineto{\pgfqpoint{2.476913in}{2.002021in}}%
\pgfpathlineto{\pgfqpoint{2.556762in}{1.910247in}}%
\pgfpathlineto{\pgfqpoint{2.636610in}{1.910138in}}%
\pgfpathlineto{\pgfqpoint{2.716459in}{1.818188in}}%
\pgfpathlineto{\pgfqpoint{2.796307in}{1.818153in}}%
\pgfpathlineto{\pgfqpoint{2.876156in}{1.726299in}}%
\pgfpathlineto{\pgfqpoint{2.956004in}{1.726190in}}%
\pgfpathlineto{\pgfqpoint{3.035853in}{1.635148in}}%
\pgfpathlineto{\pgfqpoint{3.115701in}{1.635113in}}%
\pgfpathlineto{\pgfqpoint{3.195550in}{1.542506in}}%
\pgfpathlineto{\pgfqpoint{3.275398in}{1.542386in}}%
\pgfpathlineto{\pgfqpoint{3.355247in}{1.434389in}}%
\pgfpathlineto{\pgfqpoint{3.435095in}{1.434356in}}%
\pgfpathlineto{\pgfqpoint{3.514944in}{1.283745in}}%
\pgfpathlineto{\pgfqpoint{3.594792in}{1.283713in}}%
\pgfpathlineto{\pgfqpoint{3.674641in}{1.104507in}}%
\pgfpathlineto{\pgfqpoint{3.754489in}{1.104458in}}%
\pgfpathlineto{\pgfqpoint{3.834338in}{1.051860in}}%
\pgfpathlineto{\pgfqpoint{3.914186in}{1.045004in}}%
\pgfpathlineto{\pgfqpoint{3.994034in}{1.024079in}}%
\pgfpathlineto{\pgfqpoint{4.073883in}{0.760054in}}%
\pgfpathlineto{\pgfqpoint{4.153731in}{0.689834in}}%
\pgfpathlineto{\pgfqpoint{4.226238in}{0.496358in}}%
\pgfusepath{stroke}%
\end{pgfscope}%
\begin{pgfscope}%
\pgfsetrectcap%
\pgfsetmiterjoin%
\pgfsetlinewidth{0.803000pt}%
\definecolor{currentstroke}{rgb}{0.000000,0.000000,0.000000}%
\pgfsetstrokecolor{currentstroke}%
\pgfsetdash{}{0pt}%
\pgfpathmoveto{\pgfqpoint{0.540587in}{0.499691in}}%
\pgfpathlineto{\pgfqpoint{0.540587in}{3.052241in}}%
\pgfusepath{stroke}%
\end{pgfscope}%
\begin{pgfscope}%
\pgfsetrectcap%
\pgfsetmiterjoin%
\pgfsetlinewidth{0.803000pt}%
\definecolor{currentstroke}{rgb}{0.000000,0.000000,0.000000}%
\pgfsetstrokecolor{currentstroke}%
\pgfsetdash{}{0pt}%
\pgfpathmoveto{\pgfqpoint{4.493088in}{0.499691in}}%
\pgfpathlineto{\pgfqpoint{4.493088in}{3.052241in}}%
\pgfusepath{stroke}%
\end{pgfscope}%
\begin{pgfscope}%
\pgfsetrectcap%
\pgfsetmiterjoin%
\pgfsetlinewidth{0.803000pt}%
\definecolor{currentstroke}{rgb}{0.000000,0.000000,0.000000}%
\pgfsetstrokecolor{currentstroke}%
\pgfsetdash{}{0pt}%
\pgfpathmoveto{\pgfqpoint{0.540587in}{0.499691in}}%
\pgfpathlineto{\pgfqpoint{4.493088in}{0.499691in}}%
\pgfusepath{stroke}%
\end{pgfscope}%
\begin{pgfscope}%
\pgfsetrectcap%
\pgfsetmiterjoin%
\pgfsetlinewidth{0.803000pt}%
\definecolor{currentstroke}{rgb}{0.000000,0.000000,0.000000}%
\pgfsetstrokecolor{currentstroke}%
\pgfsetdash{}{0pt}%
\pgfpathmoveto{\pgfqpoint{0.540587in}{3.052241in}}%
\pgfpathlineto{\pgfqpoint{4.493088in}{3.052241in}}%
\pgfusepath{stroke}%
\end{pgfscope}%
\begin{pgfscope}%
\pgfsetbuttcap%
\pgfsetmiterjoin%
\definecolor{currentfill}{rgb}{1.000000,1.000000,1.000000}%
\pgfsetfillcolor{currentfill}%
\pgfsetfillopacity{0.800000}%
\pgfsetlinewidth{1.003750pt}%
\definecolor{currentstroke}{rgb}{0.800000,0.800000,0.800000}%
\pgfsetstrokecolor{currentstroke}%
\pgfsetstrokeopacity{0.800000}%
\pgfsetdash{}{0pt}%
\pgfpathmoveto{\pgfqpoint{0.637810in}{0.569136in}}%
\pgfpathlineto{\pgfqpoint{3.013827in}{0.569136in}}%
\pgfpathquadraticcurveto{\pgfqpoint{3.041605in}{0.569136in}}{\pgfqpoint{3.041605in}{0.596913in}}%
\pgfpathlineto{\pgfqpoint{3.041605in}{1.015525in}}%
\pgfpathquadraticcurveto{\pgfqpoint{3.041605in}{1.043303in}}{\pgfqpoint{3.013827in}{1.043303in}}%
\pgfpathlineto{\pgfqpoint{0.637810in}{1.043303in}}%
\pgfpathquadraticcurveto{\pgfqpoint{0.610032in}{1.043303in}}{\pgfqpoint{0.610032in}{1.015525in}}%
\pgfpathlineto{\pgfqpoint{0.610032in}{0.596913in}}%
\pgfpathquadraticcurveto{\pgfqpoint{0.610032in}{0.569136in}}{\pgfqpoint{0.637810in}{0.569136in}}%
\pgfpathlineto{\pgfqpoint{0.637810in}{0.569136in}}%
\pgfpathclose%
\pgfusepath{stroke,fill}%
\end{pgfscope}%
\begin{pgfscope}%
\pgfsetrectcap%
\pgfsetroundjoin%
\pgfsetlinewidth{1.104125pt}%
\definecolor{currentstroke}{rgb}{0.000000,0.000000,0.000000}%
\pgfsetstrokecolor{currentstroke}%
\pgfsetdash{}{0pt}%
\pgfpathmoveto{\pgfqpoint{0.665587in}{0.916358in}}%
\pgfpathlineto{\pgfqpoint{0.804476in}{0.916358in}}%
\pgfpathlineto{\pgfqpoint{0.943365in}{0.916358in}}%
\pgfusepath{stroke}%
\end{pgfscope}%
\begin{pgfscope}%
\pgfsetbuttcap%
\pgfsetroundjoin%
\definecolor{currentfill}{rgb}{0.000000,0.000000,0.000000}%
\pgfsetfillcolor{currentfill}%
\pgfsetfillopacity{0.000000}%
\pgfsetlinewidth{1.003750pt}%
\definecolor{currentstroke}{rgb}{0.000000,0.000000,0.000000}%
\pgfsetstrokecolor{currentstroke}%
\pgfsetdash{}{0pt}%
\pgfsys@defobject{currentmarker}{\pgfqpoint{-0.041667in}{-0.041667in}}{\pgfqpoint{0.041667in}{0.041667in}}{%
\pgfpathmoveto{\pgfqpoint{-0.041667in}{-0.041667in}}%
\pgfpathlineto{\pgfqpoint{0.041667in}{0.041667in}}%
\pgfpathmoveto{\pgfqpoint{-0.041667in}{0.041667in}}%
\pgfpathlineto{\pgfqpoint{0.041667in}{-0.041667in}}%
\pgfusepath{stroke,fill}%
}%
\begin{pgfscope}%
\pgfsys@transformshift{0.804476in}{0.916358in}%
\pgfsys@useobject{currentmarker}{}%
\end{pgfscope}%
\end{pgfscope}%
\begin{pgfscope}%
\definecolor{textcolor}{rgb}{0.000000,0.000000,0.000000}%
\pgfsetstrokecolor{textcolor}%
\pgfsetfillcolor{textcolor}%
\pgftext[x=1.054476in,y=0.867747in,left,base]{\color{textcolor}\rmfamily\fontsize{10.000000}{12.000000}\selectfont \(\displaystyle \norm{\tns{C}_k - D^3 f(\vek{x}_k)}_F / \norm{D^3 f(\vek{x}_k)}_F\)}%
\end{pgfscope}%
\begin{pgfscope}%
\pgfsetbuttcap%
\pgfsetroundjoin%
\pgfsetlinewidth{1.104125pt}%
\definecolor{currentstroke}{rgb}{0.000000,0.000000,0.000000}%
\pgfsetstrokecolor{currentstroke}%
\pgfsetdash{{7.040000pt}{1.760000pt}{1.100000pt}{1.760000pt}}{0.000000pt}%
\pgfpathmoveto{\pgfqpoint{0.665587in}{0.708024in}}%
\pgfpathlineto{\pgfqpoint{0.804476in}{0.708024in}}%
\pgfpathlineto{\pgfqpoint{0.943365in}{0.708024in}}%
\pgfusepath{stroke}%
\end{pgfscope}%
\begin{pgfscope}%
\definecolor{textcolor}{rgb}{0.000000,0.000000,0.000000}%
\pgfsetstrokecolor{textcolor}%
\pgfsetfillcolor{textcolor}%
\pgftext[x=1.054476in,y=0.659413in,left,base]{\color{textcolor}\rmfamily\fontsize{10.000000}{12.000000}\selectfont \(\displaystyle \norm{\vek{x}_k - \vek{x}_*}_2 / \norm{\vek{x}_*}_2\)}%
\end{pgfscope}%
\end{pgfpicture}%
\makeatother%
\endgroup%

%% file: plots/convergence_linearly_converging_min.pgf
\begingroup%
\makeatletter%
\begin{pgfpicture}%
\pgfpathrectangle{\pgfpointorigin}{\pgfqpoint{4.593087in}{3.152241in}}%
\pgfusepath{use as bounding box, clip}%
\begin{pgfscope}%
\pgfsetbuttcap%
\pgfsetmiterjoin%
\definecolor{currentfill}{rgb}{1.000000,1.000000,1.000000}%
\pgfsetfillcolor{currentfill}%
\pgfsetlinewidth{0.000000pt}%
\definecolor{currentstroke}{rgb}{1.000000,1.000000,1.000000}%
\pgfsetstrokecolor{currentstroke}%
\pgfsetdash{}{0pt}%
\pgfpathmoveto{\pgfqpoint{0.000000in}{0.000000in}}%
\pgfpathlineto{\pgfqpoint{4.593087in}{0.000000in}}%
\pgfpathlineto{\pgfqpoint{4.593087in}{3.152241in}}%
\pgfpathlineto{\pgfqpoint{0.000000in}{3.152241in}}%
\pgfpathlineto{\pgfqpoint{0.000000in}{0.000000in}}%
\pgfpathclose%
\pgfusepath{fill}%
\end{pgfscope}%
\begin{pgfscope}%
\pgfsetbuttcap%
\pgfsetmiterjoin%
\definecolor{currentfill}{rgb}{1.000000,1.000000,1.000000}%
\pgfsetfillcolor{currentfill}%
\pgfsetlinewidth{0.000000pt}%
\definecolor{currentstroke}{rgb}{0.000000,0.000000,0.000000}%
\pgfsetstrokecolor{currentstroke}%
\pgfsetstrokeopacity{0.000000}%
\pgfsetdash{}{0pt}%
\pgfpathmoveto{\pgfqpoint{0.540587in}{0.499691in}}%
\pgfpathlineto{\pgfqpoint{4.493088in}{0.499691in}}%
\pgfpathlineto{\pgfqpoint{4.493088in}{3.052241in}}%
\pgfpathlineto{\pgfqpoint{0.540587in}{3.052241in}}%
\pgfpathlineto{\pgfqpoint{0.540587in}{0.499691in}}%
\pgfpathclose%
\pgfusepath{fill}%
\end{pgfscope}%
\begin{pgfscope}%
\pgfsetbuttcap%
\pgfsetroundjoin%
\definecolor{currentfill}{rgb}{0.000000,0.000000,0.000000}%
\pgfsetfillcolor{currentfill}%
\pgfsetfillopacity{0.000000}%
\pgfsetlinewidth{0.803000pt}%
\definecolor{currentstroke}{rgb}{0.000000,0.000000,0.000000}%
\pgfsetstrokecolor{currentstroke}%
\pgfsetdash{}{0pt}%
\pgfsys@defobject{currentmarker}{\pgfqpoint{0.000000in}{-0.048611in}}{\pgfqpoint{0.000000in}{0.000000in}}{%
\pgfpathmoveto{\pgfqpoint{0.000000in}{0.000000in}}%
\pgfpathlineto{\pgfqpoint{0.000000in}{-0.048611in}}%
\pgfusepath{stroke,fill}%
}%
\begin{pgfscope}%
\pgfsys@transformshift{0.720247in}{0.499691in}%
\pgfsys@useobject{currentmarker}{}%
\end{pgfscope}%
\end{pgfscope}%
\begin{pgfscope}%
\definecolor{textcolor}{rgb}{0.000000,0.000000,0.000000}%
\pgfsetstrokecolor{textcolor}%
\pgfsetfillcolor{textcolor}%
\pgftext[x=0.720247in,y=0.402469in,,top]{\color{textcolor}\rmfamily\fontsize{10.000000}{12.000000}\selectfont \(\displaystyle {0}\)}%
\end{pgfscope}%
\begin{pgfscope}%
\pgfsetbuttcap%
\pgfsetroundjoin%
\definecolor{currentfill}{rgb}{0.000000,0.000000,0.000000}%
\pgfsetfillcolor{currentfill}%
\pgfsetfillopacity{0.000000}%
\pgfsetlinewidth{0.803000pt}%
\definecolor{currentstroke}{rgb}{0.000000,0.000000,0.000000}%
\pgfsetstrokecolor{currentstroke}%
\pgfsetdash{}{0pt}%
\pgfsys@defobject{currentmarker}{\pgfqpoint{0.000000in}{-0.048611in}}{\pgfqpoint{0.000000in}{0.000000in}}{%
\pgfpathmoveto{\pgfqpoint{0.000000in}{0.000000in}}%
\pgfpathlineto{\pgfqpoint{0.000000in}{-0.048611in}}%
\pgfusepath{stroke,fill}%
}%
\begin{pgfscope}%
\pgfsys@transformshift{1.518731in}{0.499691in}%
\pgfsys@useobject{currentmarker}{}%
\end{pgfscope}%
\end{pgfscope}%
\begin{pgfscope}%
\definecolor{textcolor}{rgb}{0.000000,0.000000,0.000000}%
\pgfsetstrokecolor{textcolor}%
\pgfsetfillcolor{textcolor}%
\pgftext[x=1.518731in,y=0.402469in,,top]{\color{textcolor}\rmfamily\fontsize{10.000000}{12.000000}\selectfont \(\displaystyle {10}\)}%
\end{pgfscope}%
\begin{pgfscope}%
\pgfsetbuttcap%
\pgfsetroundjoin%
\definecolor{currentfill}{rgb}{0.000000,0.000000,0.000000}%
\pgfsetfillcolor{currentfill}%
\pgfsetfillopacity{0.000000}%
\pgfsetlinewidth{0.803000pt}%
\definecolor{currentstroke}{rgb}{0.000000,0.000000,0.000000}%
\pgfsetstrokecolor{currentstroke}%
\pgfsetdash{}{0pt}%
\pgfsys@defobject{currentmarker}{\pgfqpoint{0.000000in}{-0.048611in}}{\pgfqpoint{0.000000in}{0.000000in}}{%
\pgfpathmoveto{\pgfqpoint{0.000000in}{0.000000in}}%
\pgfpathlineto{\pgfqpoint{0.000000in}{-0.048611in}}%
\pgfusepath{stroke,fill}%
}%
\begin{pgfscope}%
\pgfsys@transformshift{2.317216in}{0.499691in}%
\pgfsys@useobject{currentmarker}{}%
\end{pgfscope}%
\end{pgfscope}%
\begin{pgfscope}%
\definecolor{textcolor}{rgb}{0.000000,0.000000,0.000000}%
\pgfsetstrokecolor{textcolor}%
\pgfsetfillcolor{textcolor}%
\pgftext[x=2.317216in,y=0.402469in,,top]{\color{textcolor}\rmfamily\fontsize{10.000000}{12.000000}\selectfont \(\displaystyle {20}\)}%
\end{pgfscope}%
\begin{pgfscope}%
\pgfsetbuttcap%
\pgfsetroundjoin%
\definecolor{currentfill}{rgb}{0.000000,0.000000,0.000000}%
\pgfsetfillcolor{currentfill}%
\pgfsetfillopacity{0.000000}%
\pgfsetlinewidth{0.803000pt}%
\definecolor{currentstroke}{rgb}{0.000000,0.000000,0.000000}%
\pgfsetstrokecolor{currentstroke}%
\pgfsetdash{}{0pt}%
\pgfsys@defobject{currentmarker}{\pgfqpoint{0.000000in}{-0.048611in}}{\pgfqpoint{0.000000in}{0.000000in}}{%
\pgfpathmoveto{\pgfqpoint{0.000000in}{0.000000in}}%
\pgfpathlineto{\pgfqpoint{0.000000in}{-0.048611in}}%
\pgfusepath{stroke,fill}%
}%
\begin{pgfscope}%
\pgfsys@transformshift{3.115701in}{0.499691in}%
\pgfsys@useobject{currentmarker}{}%
\end{pgfscope}%
\end{pgfscope}%
\begin{pgfscope}%
\definecolor{textcolor}{rgb}{0.000000,0.000000,0.000000}%
\pgfsetstrokecolor{textcolor}%
\pgfsetfillcolor{textcolor}%
\pgftext[x=3.115701in,y=0.402469in,,top]{\color{textcolor}\rmfamily\fontsize{10.000000}{12.000000}\selectfont \(\displaystyle {30}\)}%
\end{pgfscope}%
\begin{pgfscope}%
\pgfsetbuttcap%
\pgfsetroundjoin%
\definecolor{currentfill}{rgb}{0.000000,0.000000,0.000000}%
\pgfsetfillcolor{currentfill}%
\pgfsetfillopacity{0.000000}%
\pgfsetlinewidth{0.803000pt}%
\definecolor{currentstroke}{rgb}{0.000000,0.000000,0.000000}%
\pgfsetstrokecolor{currentstroke}%
\pgfsetdash{}{0pt}%
\pgfsys@defobject{currentmarker}{\pgfqpoint{0.000000in}{-0.048611in}}{\pgfqpoint{0.000000in}{0.000000in}}{%
\pgfpathmoveto{\pgfqpoint{0.000000in}{0.000000in}}%
\pgfpathlineto{\pgfqpoint{0.000000in}{-0.048611in}}%
\pgfusepath{stroke,fill}%
}%
\begin{pgfscope}%
\pgfsys@transformshift{3.914186in}{0.499691in}%
\pgfsys@useobject{currentmarker}{}%
\end{pgfscope}%
\end{pgfscope}%
\begin{pgfscope}%
\definecolor{textcolor}{rgb}{0.000000,0.000000,0.000000}%
\pgfsetstrokecolor{textcolor}%
\pgfsetfillcolor{textcolor}%
\pgftext[x=3.914186in,y=0.402469in,,top]{\color{textcolor}\rmfamily\fontsize{10.000000}{12.000000}\selectfont \(\displaystyle {40}\)}%
\end{pgfscope}%
\begin{pgfscope}%
\definecolor{textcolor}{rgb}{0.000000,0.000000,0.000000}%
\pgfsetstrokecolor{textcolor}%
\pgfsetfillcolor{textcolor}%
\pgftext[x=2.516837in,y=0.223457in,,top]{\color{textcolor}\rmfamily\fontsize{10.000000}{12.000000}\selectfont Iteration \(\displaystyle k\)}%
\end{pgfscope}%
\begin{pgfscope}%
\pgfsetbuttcap%
\pgfsetroundjoin%
\definecolor{currentfill}{rgb}{0.000000,0.000000,0.000000}%
\pgfsetfillcolor{currentfill}%
\pgfsetfillopacity{0.000000}%
\pgfsetlinewidth{0.803000pt}%
\definecolor{currentstroke}{rgb}{0.000000,0.000000,0.000000}%
\pgfsetstrokecolor{currentstroke}%
\pgfsetdash{}{0pt}%
\pgfsys@defobject{currentmarker}{\pgfqpoint{-0.048611in}{0.000000in}}{\pgfqpoint{-0.000000in}{0.000000in}}{%
\pgfpathmoveto{\pgfqpoint{-0.000000in}{0.000000in}}%
\pgfpathlineto{\pgfqpoint{-0.048611in}{0.000000in}}%
\pgfusepath{stroke,fill}%
}%
\begin{pgfscope}%
\pgfsys@transformshift{0.540587in}{0.791868in}%
\pgfsys@useobject{currentmarker}{}%
\end{pgfscope}%
\end{pgfscope}%
\begin{pgfscope}%
\definecolor{textcolor}{rgb}{0.000000,0.000000,0.000000}%
\pgfsetstrokecolor{textcolor}%
\pgfsetfillcolor{textcolor}%
\pgftext[x=0.100000in, y=0.743643in, left, base]{\color{textcolor}\rmfamily\fontsize{10.000000}{12.000000}\selectfont \(\displaystyle {10^{-14}}\)}%
\end{pgfscope}%
\begin{pgfscope}%
\pgfsetbuttcap%
\pgfsetroundjoin%
\definecolor{currentfill}{rgb}{0.000000,0.000000,0.000000}%
\pgfsetfillcolor{currentfill}%
\pgfsetfillopacity{0.000000}%
\pgfsetlinewidth{0.803000pt}%
\definecolor{currentstroke}{rgb}{0.000000,0.000000,0.000000}%
\pgfsetstrokecolor{currentstroke}%
\pgfsetdash{}{0pt}%
\pgfsys@defobject{currentmarker}{\pgfqpoint{-0.048611in}{0.000000in}}{\pgfqpoint{-0.000000in}{0.000000in}}{%
\pgfpathmoveto{\pgfqpoint{-0.000000in}{0.000000in}}%
\pgfpathlineto{\pgfqpoint{-0.048611in}{0.000000in}}%
\pgfusepath{stroke,fill}%
}%
\begin{pgfscope}%
\pgfsys@transformshift{0.540587in}{1.090833in}%
\pgfsys@useobject{currentmarker}{}%
\end{pgfscope}%
\end{pgfscope}%
\begin{pgfscope}%
\definecolor{textcolor}{rgb}{0.000000,0.000000,0.000000}%
\pgfsetstrokecolor{textcolor}%
\pgfsetfillcolor{textcolor}%
\pgftext[x=0.100000in, y=1.042607in, left, base]{\color{textcolor}\rmfamily\fontsize{10.000000}{12.000000}\selectfont \(\displaystyle {10^{-12}}\)}%
\end{pgfscope}%
\begin{pgfscope}%
\pgfsetbuttcap%
\pgfsetroundjoin%
\definecolor{currentfill}{rgb}{0.000000,0.000000,0.000000}%
\pgfsetfillcolor{currentfill}%
\pgfsetfillopacity{0.000000}%
\pgfsetlinewidth{0.803000pt}%
\definecolor{currentstroke}{rgb}{0.000000,0.000000,0.000000}%
\pgfsetstrokecolor{currentstroke}%
\pgfsetdash{}{0pt}%
\pgfsys@defobject{currentmarker}{\pgfqpoint{-0.048611in}{0.000000in}}{\pgfqpoint{-0.000000in}{0.000000in}}{%
\pgfpathmoveto{\pgfqpoint{-0.000000in}{0.000000in}}%
\pgfpathlineto{\pgfqpoint{-0.048611in}{0.000000in}}%
\pgfusepath{stroke,fill}%
}%
\begin{pgfscope}%
\pgfsys@transformshift{0.540587in}{1.389797in}%
\pgfsys@useobject{currentmarker}{}%
\end{pgfscope}%
\end{pgfscope}%
\begin{pgfscope}%
\definecolor{textcolor}{rgb}{0.000000,0.000000,0.000000}%
\pgfsetstrokecolor{textcolor}%
\pgfsetfillcolor{textcolor}%
\pgftext[x=0.100000in, y=1.341572in, left, base]{\color{textcolor}\rmfamily\fontsize{10.000000}{12.000000}\selectfont \(\displaystyle {10^{-10}}\)}%
\end{pgfscope}%
\begin{pgfscope}%
\pgfsetbuttcap%
\pgfsetroundjoin%
\definecolor{currentfill}{rgb}{0.000000,0.000000,0.000000}%
\pgfsetfillcolor{currentfill}%
\pgfsetfillopacity{0.000000}%
\pgfsetlinewidth{0.803000pt}%
\definecolor{currentstroke}{rgb}{0.000000,0.000000,0.000000}%
\pgfsetstrokecolor{currentstroke}%
\pgfsetdash{}{0pt}%
\pgfsys@defobject{currentmarker}{\pgfqpoint{-0.048611in}{0.000000in}}{\pgfqpoint{-0.000000in}{0.000000in}}{%
\pgfpathmoveto{\pgfqpoint{-0.000000in}{0.000000in}}%
\pgfpathlineto{\pgfqpoint{-0.048611in}{0.000000in}}%
\pgfusepath{stroke,fill}%
}%
\begin{pgfscope}%
\pgfsys@transformshift{0.540587in}{1.688762in}%
\pgfsys@useobject{currentmarker}{}%
\end{pgfscope}%
\end{pgfscope}%
\begin{pgfscope}%
\definecolor{textcolor}{rgb}{0.000000,0.000000,0.000000}%
\pgfsetstrokecolor{textcolor}%
\pgfsetfillcolor{textcolor}%
\pgftext[x=0.155363in, y=1.640537in, left, base]{\color{textcolor}\rmfamily\fontsize{10.000000}{12.000000}\selectfont \(\displaystyle {10^{-8}}\)}%
\end{pgfscope}%
\begin{pgfscope}%
\pgfsetbuttcap%
\pgfsetroundjoin%
\definecolor{currentfill}{rgb}{0.000000,0.000000,0.000000}%
\pgfsetfillcolor{currentfill}%
\pgfsetfillopacity{0.000000}%
\pgfsetlinewidth{0.803000pt}%
\definecolor{currentstroke}{rgb}{0.000000,0.000000,0.000000}%
\pgfsetstrokecolor{currentstroke}%
\pgfsetdash{}{0pt}%
\pgfsys@defobject{currentmarker}{\pgfqpoint{-0.048611in}{0.000000in}}{\pgfqpoint{-0.000000in}{0.000000in}}{%
\pgfpathmoveto{\pgfqpoint{-0.000000in}{0.000000in}}%
\pgfpathlineto{\pgfqpoint{-0.048611in}{0.000000in}}%
\pgfusepath{stroke,fill}%
}%
\begin{pgfscope}%
\pgfsys@transformshift{0.540587in}{1.987727in}%
\pgfsys@useobject{currentmarker}{}%
\end{pgfscope}%
\end{pgfscope}%
\begin{pgfscope}%
\definecolor{textcolor}{rgb}{0.000000,0.000000,0.000000}%
\pgfsetstrokecolor{textcolor}%
\pgfsetfillcolor{textcolor}%
\pgftext[x=0.155363in, y=1.939502in, left, base]{\color{textcolor}\rmfamily\fontsize{10.000000}{12.000000}\selectfont \(\displaystyle {10^{-6}}\)}%
\end{pgfscope}%
\begin{pgfscope}%
\pgfsetbuttcap%
\pgfsetroundjoin%
\definecolor{currentfill}{rgb}{0.000000,0.000000,0.000000}%
\pgfsetfillcolor{currentfill}%
\pgfsetfillopacity{0.000000}%
\pgfsetlinewidth{0.803000pt}%
\definecolor{currentstroke}{rgb}{0.000000,0.000000,0.000000}%
\pgfsetstrokecolor{currentstroke}%
\pgfsetdash{}{0pt}%
\pgfsys@defobject{currentmarker}{\pgfqpoint{-0.048611in}{0.000000in}}{\pgfqpoint{-0.000000in}{0.000000in}}{%
\pgfpathmoveto{\pgfqpoint{-0.000000in}{0.000000in}}%
\pgfpathlineto{\pgfqpoint{-0.048611in}{0.000000in}}%
\pgfusepath{stroke,fill}%
}%
\begin{pgfscope}%
\pgfsys@transformshift{0.540587in}{2.286692in}%
\pgfsys@useobject{currentmarker}{}%
\end{pgfscope}%
\end{pgfscope}%
\begin{pgfscope}%
\definecolor{textcolor}{rgb}{0.000000,0.000000,0.000000}%
\pgfsetstrokecolor{textcolor}%
\pgfsetfillcolor{textcolor}%
\pgftext[x=0.155363in, y=2.238466in, left, base]{\color{textcolor}\rmfamily\fontsize{10.000000}{12.000000}\selectfont \(\displaystyle {10^{-4}}\)}%
\end{pgfscope}%
\begin{pgfscope}%
\pgfsetbuttcap%
\pgfsetroundjoin%
\definecolor{currentfill}{rgb}{0.000000,0.000000,0.000000}%
\pgfsetfillcolor{currentfill}%
\pgfsetfillopacity{0.000000}%
\pgfsetlinewidth{0.803000pt}%
\definecolor{currentstroke}{rgb}{0.000000,0.000000,0.000000}%
\pgfsetstrokecolor{currentstroke}%
\pgfsetdash{}{0pt}%
\pgfsys@defobject{currentmarker}{\pgfqpoint{-0.048611in}{0.000000in}}{\pgfqpoint{-0.000000in}{0.000000in}}{%
\pgfpathmoveto{\pgfqpoint{-0.000000in}{0.000000in}}%
\pgfpathlineto{\pgfqpoint{-0.048611in}{0.000000in}}%
\pgfusepath{stroke,fill}%
}%
\begin{pgfscope}%
\pgfsys@transformshift{0.540587in}{2.585656in}%
\pgfsys@useobject{currentmarker}{}%
\end{pgfscope}%
\end{pgfscope}%
\begin{pgfscope}%
\definecolor{textcolor}{rgb}{0.000000,0.000000,0.000000}%
\pgfsetstrokecolor{textcolor}%
\pgfsetfillcolor{textcolor}%
\pgftext[x=0.155363in, y=2.537431in, left, base]{\color{textcolor}\rmfamily\fontsize{10.000000}{12.000000}\selectfont \(\displaystyle {10^{-2}}\)}%
\end{pgfscope}%
\begin{pgfscope}%
\pgfsetbuttcap%
\pgfsetroundjoin%
\definecolor{currentfill}{rgb}{0.000000,0.000000,0.000000}%
\pgfsetfillcolor{currentfill}%
\pgfsetfillopacity{0.000000}%
\pgfsetlinewidth{0.803000pt}%
\definecolor{currentstroke}{rgb}{0.000000,0.000000,0.000000}%
\pgfsetstrokecolor{currentstroke}%
\pgfsetdash{}{0pt}%
\pgfsys@defobject{currentmarker}{\pgfqpoint{-0.048611in}{0.000000in}}{\pgfqpoint{-0.000000in}{0.000000in}}{%
\pgfpathmoveto{\pgfqpoint{-0.000000in}{0.000000in}}%
\pgfpathlineto{\pgfqpoint{-0.048611in}{0.000000in}}%
\pgfusepath{stroke,fill}%
}%
\begin{pgfscope}%
\pgfsys@transformshift{0.540587in}{2.884621in}%
\pgfsys@useobject{currentmarker}{}%
\end{pgfscope}%
\end{pgfscope}%
\begin{pgfscope}%
\definecolor{textcolor}{rgb}{0.000000,0.000000,0.000000}%
\pgfsetstrokecolor{textcolor}%
\pgfsetfillcolor{textcolor}%
\pgftext[x=0.242169in, y=2.836396in, left, base]{\color{textcolor}\rmfamily\fontsize{10.000000}{12.000000}\selectfont \(\displaystyle {10^{0}}\)}%
\end{pgfscope}%
\begin{pgfscope}%
\pgfpathrectangle{\pgfqpoint{0.540587in}{0.499691in}}{\pgfqpoint{3.952500in}{2.552550in}}%
\pgfusepath{clip}%
\pgfsetrectcap%
\pgfsetroundjoin%
\pgfsetlinewidth{1.104125pt}%
\definecolor{currentstroke}{rgb}{0.000000,0.000000,0.000000}%
\pgfsetstrokecolor{currentstroke}%
\pgfsetdash{}{0pt}%
\pgfpathmoveto{\pgfqpoint{0.720247in}{2.884621in}}%
\pgfpathlineto{\pgfqpoint{0.800095in}{2.797333in}}%
\pgfpathlineto{\pgfqpoint{0.879944in}{2.727637in}}%
\pgfpathlineto{\pgfqpoint{0.959792in}{2.745347in}}%
\pgfpathlineto{\pgfqpoint{1.039641in}{2.753544in}}%
\pgfpathlineto{\pgfqpoint{1.119489in}{2.779339in}}%
\pgfpathlineto{\pgfqpoint{1.199337in}{2.721730in}}%
\pgfpathlineto{\pgfqpoint{1.279186in}{2.731552in}}%
\pgfpathlineto{\pgfqpoint{1.359034in}{2.765793in}}%
\pgfpathlineto{\pgfqpoint{1.438883in}{2.756225in}}%
\pgfpathlineto{\pgfqpoint{1.518731in}{2.757114in}}%
\pgfpathlineto{\pgfqpoint{1.598580in}{2.773448in}}%
\pgfpathlineto{\pgfqpoint{1.678428in}{2.775546in}}%
\pgfpathlineto{\pgfqpoint{1.758277in}{2.785402in}}%
\pgfpathlineto{\pgfqpoint{1.838125in}{2.792755in}}%
\pgfpathlineto{\pgfqpoint{1.917974in}{2.485103in}}%
\pgfpathlineto{\pgfqpoint{1.997822in}{2.490012in}}%
\pgfpathlineto{\pgfqpoint{2.077671in}{2.433509in}}%
\pgfpathlineto{\pgfqpoint{2.157519in}{2.034761in}}%
\pgfpathlineto{\pgfqpoint{2.237368in}{2.081988in}}%
\pgfpathlineto{\pgfqpoint{2.317216in}{2.020884in}}%
\pgfpathlineto{\pgfqpoint{2.397065in}{2.039853in}}%
\pgfpathlineto{\pgfqpoint{2.476913in}{1.957014in}}%
\pgfpathlineto{\pgfqpoint{2.556762in}{1.909730in}}%
\pgfpathlineto{\pgfqpoint{2.636610in}{1.855526in}}%
\pgfpathlineto{\pgfqpoint{2.716459in}{1.853518in}}%
\pgfpathlineto{\pgfqpoint{2.796307in}{1.868905in}}%
\pgfpathlineto{\pgfqpoint{2.876156in}{1.845277in}}%
\pgfpathlineto{\pgfqpoint{2.956004in}{1.967075in}}%
\pgfpathlineto{\pgfqpoint{3.035853in}{1.951555in}}%
\pgfpathlineto{\pgfqpoint{3.115701in}{2.145722in}}%
\pgfpathlineto{\pgfqpoint{3.195550in}{2.143814in}}%
\pgfpathlineto{\pgfqpoint{3.275398in}{2.208943in}}%
\pgfpathlineto{\pgfqpoint{3.355247in}{2.169896in}}%
\pgfpathlineto{\pgfqpoint{3.435095in}{2.283760in}}%
\pgfpathlineto{\pgfqpoint{3.514944in}{2.280045in}}%
\pgfpathlineto{\pgfqpoint{3.594792in}{2.396726in}}%
\pgfpathlineto{\pgfqpoint{3.674641in}{2.344475in}}%
\pgfpathlineto{\pgfqpoint{3.754489in}{2.564642in}}%
\pgfpathlineto{\pgfqpoint{3.834338in}{2.526462in}}%
\pgfpathlineto{\pgfqpoint{3.914186in}{2.535962in}}%
\pgfpathlineto{\pgfqpoint{3.994034in}{2.518730in}}%
\pgfpathlineto{\pgfqpoint{4.073883in}{2.518557in}}%
\pgfpathlineto{\pgfqpoint{4.153731in}{2.609485in}}%
\pgfpathlineto{\pgfqpoint{4.233580in}{2.564773in}}%
\pgfpathlineto{\pgfqpoint{4.313428in}{2.766016in}}%
\pgfusepath{stroke}%
\end{pgfscope}%
\begin{pgfscope}%
\pgfpathrectangle{\pgfqpoint{0.540587in}{0.499691in}}{\pgfqpoint{3.952500in}{2.552550in}}%
\pgfusepath{clip}%
\pgfsetbuttcap%
\pgfsetroundjoin%
\definecolor{currentfill}{rgb}{0.000000,0.000000,0.000000}%
\pgfsetfillcolor{currentfill}%
\pgfsetfillopacity{0.000000}%
\pgfsetlinewidth{1.003750pt}%
\definecolor{currentstroke}{rgb}{0.000000,0.000000,0.000000}%
\pgfsetstrokecolor{currentstroke}%
\pgfsetdash{}{0pt}%
\pgfsys@defobject{currentmarker}{\pgfqpoint{-0.041667in}{-0.041667in}}{\pgfqpoint{0.041667in}{0.041667in}}{%
\pgfpathmoveto{\pgfqpoint{-0.041667in}{-0.041667in}}%
\pgfpathlineto{\pgfqpoint{0.041667in}{0.041667in}}%
\pgfpathmoveto{\pgfqpoint{-0.041667in}{0.041667in}}%
\pgfpathlineto{\pgfqpoint{0.041667in}{-0.041667in}}%
\pgfusepath{stroke,fill}%
}%
\begin{pgfscope}%
\pgfsys@transformshift{0.720247in}{2.884621in}%
\pgfsys@useobject{currentmarker}{}%
\end{pgfscope}%
\begin{pgfscope}%
\pgfsys@transformshift{0.800095in}{2.797333in}%
\pgfsys@useobject{currentmarker}{}%
\end{pgfscope}%
\begin{pgfscope}%
\pgfsys@transformshift{0.879944in}{2.727637in}%
\pgfsys@useobject{currentmarker}{}%
\end{pgfscope}%
\begin{pgfscope}%
\pgfsys@transformshift{0.959792in}{2.745347in}%
\pgfsys@useobject{currentmarker}{}%
\end{pgfscope}%
\begin{pgfscope}%
\pgfsys@transformshift{1.039641in}{2.753544in}%
\pgfsys@useobject{currentmarker}{}%
\end{pgfscope}%
\begin{pgfscope}%
\pgfsys@transformshift{1.119489in}{2.779339in}%
\pgfsys@useobject{currentmarker}{}%
\end{pgfscope}%
\begin{pgfscope}%
\pgfsys@transformshift{1.199337in}{2.721730in}%
\pgfsys@useobject{currentmarker}{}%
\end{pgfscope}%
\begin{pgfscope}%
\pgfsys@transformshift{1.279186in}{2.731552in}%
\pgfsys@useobject{currentmarker}{}%
\end{pgfscope}%
\begin{pgfscope}%
\pgfsys@transformshift{1.359034in}{2.765793in}%
\pgfsys@useobject{currentmarker}{}%
\end{pgfscope}%
\begin{pgfscope}%
\pgfsys@transformshift{1.438883in}{2.756225in}%
\pgfsys@useobject{currentmarker}{}%
\end{pgfscope}%
\begin{pgfscope}%
\pgfsys@transformshift{1.518731in}{2.757114in}%
\pgfsys@useobject{currentmarker}{}%
\end{pgfscope}%
\begin{pgfscope}%
\pgfsys@transformshift{1.598580in}{2.773448in}%
\pgfsys@useobject{currentmarker}{}%
\end{pgfscope}%
\begin{pgfscope}%
\pgfsys@transformshift{1.678428in}{2.775546in}%
\pgfsys@useobject{currentmarker}{}%
\end{pgfscope}%
\begin{pgfscope}%
\pgfsys@transformshift{1.758277in}{2.785402in}%
\pgfsys@useobject{currentmarker}{}%
\end{pgfscope}%
\begin{pgfscope}%
\pgfsys@transformshift{1.838125in}{2.792755in}%
\pgfsys@useobject{currentmarker}{}%
\end{pgfscope}%
\begin{pgfscope}%
\pgfsys@transformshift{1.917974in}{2.485103in}%
\pgfsys@useobject{currentmarker}{}%
\end{pgfscope}%
\begin{pgfscope}%
\pgfsys@transformshift{1.997822in}{2.490012in}%
\pgfsys@useobject{currentmarker}{}%
\end{pgfscope}%
\begin{pgfscope}%
\pgfsys@transformshift{2.077671in}{2.433509in}%
\pgfsys@useobject{currentmarker}{}%
\end{pgfscope}%
\begin{pgfscope}%
\pgfsys@transformshift{2.157519in}{2.034761in}%
\pgfsys@useobject{currentmarker}{}%
\end{pgfscope}%
\begin{pgfscope}%
\pgfsys@transformshift{2.237368in}{2.081988in}%
\pgfsys@useobject{currentmarker}{}%
\end{pgfscope}%
\begin{pgfscope}%
\pgfsys@transformshift{2.317216in}{2.020884in}%
\pgfsys@useobject{currentmarker}{}%
\end{pgfscope}%
\begin{pgfscope}%
\pgfsys@transformshift{2.397065in}{2.039853in}%
\pgfsys@useobject{currentmarker}{}%
\end{pgfscope}%
\begin{pgfscope}%
\pgfsys@transformshift{2.476913in}{1.957014in}%
\pgfsys@useobject{currentmarker}{}%
\end{pgfscope}%
\begin{pgfscope}%
\pgfsys@transformshift{2.556762in}{1.909730in}%
\pgfsys@useobject{currentmarker}{}%
\end{pgfscope}%
\begin{pgfscope}%
\pgfsys@transformshift{2.636610in}{1.855526in}%
\pgfsys@useobject{currentmarker}{}%
\end{pgfscope}%
\begin{pgfscope}%
\pgfsys@transformshift{2.716459in}{1.853518in}%
\pgfsys@useobject{currentmarker}{}%
\end{pgfscope}%
\begin{pgfscope}%
\pgfsys@transformshift{2.796307in}{1.868905in}%
\pgfsys@useobject{currentmarker}{}%
\end{pgfscope}%
\begin{pgfscope}%
\pgfsys@transformshift{2.876156in}{1.845277in}%
\pgfsys@useobject{currentmarker}{}%
\end{pgfscope}%
\begin{pgfscope}%
\pgfsys@transformshift{2.956004in}{1.967075in}%
\pgfsys@useobject{currentmarker}{}%
\end{pgfscope}%
\begin{pgfscope}%
\pgfsys@transformshift{3.035853in}{1.951555in}%
\pgfsys@useobject{currentmarker}{}%
\end{pgfscope}%
\begin{pgfscope}%
\pgfsys@transformshift{3.115701in}{2.145722in}%
\pgfsys@useobject{currentmarker}{}%
\end{pgfscope}%
\begin{pgfscope}%
\pgfsys@transformshift{3.195550in}{2.143814in}%
\pgfsys@useobject{currentmarker}{}%
\end{pgfscope}%
\begin{pgfscope}%
\pgfsys@transformshift{3.275398in}{2.208943in}%
\pgfsys@useobject{currentmarker}{}%
\end{pgfscope}%
\begin{pgfscope}%
\pgfsys@transformshift{3.355247in}{2.169896in}%
\pgfsys@useobject{currentmarker}{}%
\end{pgfscope}%
\begin{pgfscope}%
\pgfsys@transformshift{3.435095in}{2.283760in}%
\pgfsys@useobject{currentmarker}{}%
\end{pgfscope}%
\begin{pgfscope}%
\pgfsys@transformshift{3.514944in}{2.280045in}%
\pgfsys@useobject{currentmarker}{}%
\end{pgfscope}%
\begin{pgfscope}%
\pgfsys@transformshift{3.594792in}{2.396726in}%
\pgfsys@useobject{currentmarker}{}%
\end{pgfscope}%
\begin{pgfscope}%
\pgfsys@transformshift{3.674641in}{2.344475in}%
\pgfsys@useobject{currentmarker}{}%
\end{pgfscope}%
\begin{pgfscope}%
\pgfsys@transformshift{3.754489in}{2.564642in}%
\pgfsys@useobject{currentmarker}{}%
\end{pgfscope}%
\begin{pgfscope}%
\pgfsys@transformshift{3.834338in}{2.526462in}%
\pgfsys@useobject{currentmarker}{}%
\end{pgfscope}%
\begin{pgfscope}%
\pgfsys@transformshift{3.914186in}{2.535962in}%
\pgfsys@useobject{currentmarker}{}%
\end{pgfscope}%
\begin{pgfscope}%
\pgfsys@transformshift{3.994034in}{2.518730in}%
\pgfsys@useobject{currentmarker}{}%
\end{pgfscope}%
\begin{pgfscope}%
\pgfsys@transformshift{4.073883in}{2.518557in}%
\pgfsys@useobject{currentmarker}{}%
\end{pgfscope}%
\begin{pgfscope}%
\pgfsys@transformshift{4.153731in}{2.609485in}%
\pgfsys@useobject{currentmarker}{}%
\end{pgfscope}%
\begin{pgfscope}%
\pgfsys@transformshift{4.233580in}{2.564773in}%
\pgfsys@useobject{currentmarker}{}%
\end{pgfscope}%
\begin{pgfscope}%
\pgfsys@transformshift{4.313428in}{2.766016in}%
\pgfsys@useobject{currentmarker}{}%
\end{pgfscope}%
\end{pgfscope}%
\begin{pgfscope}%
\pgfpathrectangle{\pgfqpoint{0.540587in}{0.499691in}}{\pgfqpoint{3.952500in}{2.552550in}}%
\pgfusepath{clip}%
\pgfsetbuttcap%
\pgfsetroundjoin%
\pgfsetlinewidth{1.104125pt}%
\definecolor{currentstroke}{rgb}{0.000000,0.000000,0.000000}%
\pgfsetstrokecolor{currentstroke}%
\pgfsetdash{{7.040000pt}{1.760000pt}{1.100000pt}{1.760000pt}}{0.000000pt}%
\pgfpathmoveto{\pgfqpoint{0.720247in}{2.884621in}}%
\pgfpathlineto{\pgfqpoint{0.800095in}{2.878931in}}%
\pgfpathlineto{\pgfqpoint{0.879944in}{2.875944in}}%
\pgfpathlineto{\pgfqpoint{0.959792in}{2.872671in}}%
\pgfpathlineto{\pgfqpoint{1.039641in}{2.866033in}}%
\pgfpathlineto{\pgfqpoint{1.119489in}{2.850607in}}%
\pgfpathlineto{\pgfqpoint{1.199337in}{2.850427in}}%
\pgfpathlineto{\pgfqpoint{1.279186in}{2.836435in}}%
\pgfpathlineto{\pgfqpoint{1.359034in}{2.811454in}}%
\pgfpathlineto{\pgfqpoint{1.438883in}{2.809337in}}%
\pgfpathlineto{\pgfqpoint{1.518731in}{2.808120in}}%
\pgfpathlineto{\pgfqpoint{1.598580in}{2.782588in}}%
\pgfpathlineto{\pgfqpoint{1.678428in}{2.774715in}}%
\pgfpathlineto{\pgfqpoint{1.758277in}{2.720079in}}%
\pgfpathlineto{\pgfqpoint{1.838125in}{2.534402in}}%
\pgfpathlineto{\pgfqpoint{1.917974in}{2.534272in}}%
\pgfpathlineto{\pgfqpoint{1.997822in}{2.120529in}}%
\pgfpathlineto{\pgfqpoint{2.077671in}{2.152128in}}%
\pgfpathlineto{\pgfqpoint{2.157519in}{2.152094in}}%
\pgfpathlineto{\pgfqpoint{2.237368in}{2.064652in}}%
\pgfpathlineto{\pgfqpoint{2.317216in}{2.064545in}}%
\pgfpathlineto{\pgfqpoint{2.397065in}{1.974247in}}%
\pgfpathlineto{\pgfqpoint{2.476913in}{1.974213in}}%
\pgfpathlineto{\pgfqpoint{2.556762in}{1.884137in}}%
\pgfpathlineto{\pgfqpoint{2.636610in}{1.884030in}}%
\pgfpathlineto{\pgfqpoint{2.716459in}{1.793782in}}%
\pgfpathlineto{\pgfqpoint{2.796307in}{1.793748in}}%
\pgfpathlineto{\pgfqpoint{2.876156in}{1.703594in}}%
\pgfpathlineto{\pgfqpoint{2.956004in}{1.703487in}}%
\pgfpathlineto{\pgfqpoint{3.035853in}{1.614130in}}%
\pgfpathlineto{\pgfqpoint{3.115701in}{1.614096in}}%
\pgfpathlineto{\pgfqpoint{3.195550in}{1.523204in}}%
\pgfpathlineto{\pgfqpoint{3.275398in}{1.523085in}}%
\pgfpathlineto{\pgfqpoint{3.355247in}{1.417087in}}%
\pgfpathlineto{\pgfqpoint{3.435095in}{1.417055in}}%
\pgfpathlineto{\pgfqpoint{3.514944in}{1.269232in}}%
\pgfpathlineto{\pgfqpoint{3.594792in}{1.269201in}}%
\pgfpathlineto{\pgfqpoint{3.674641in}{1.093312in}}%
\pgfpathlineto{\pgfqpoint{3.754489in}{1.093263in}}%
\pgfpathlineto{\pgfqpoint{3.834338in}{1.041639in}}%
\pgfpathlineto{\pgfqpoint{3.914186in}{1.034910in}}%
\pgfpathlineto{\pgfqpoint{3.994034in}{1.014373in}}%
\pgfpathlineto{\pgfqpoint{4.073883in}{0.755235in}}%
\pgfpathlineto{\pgfqpoint{4.153731in}{0.686315in}}%
\pgfpathlineto{\pgfqpoint{4.226262in}{0.496358in}}%
\pgfusepath{stroke}%
\end{pgfscope}%
\begin{pgfscope}%
\pgfpathrectangle{\pgfqpoint{0.540587in}{0.499691in}}{\pgfqpoint{3.952500in}{2.552550in}}%
\pgfusepath{clip}%
\pgfsetbuttcap%
\pgfsetroundjoin%
\pgfsetlinewidth{1.104125pt}%
\definecolor{currentstroke}{rgb}{0.000000,0.000000,0.000000}%
\pgfsetstrokecolor{currentstroke}%
\pgfsetdash{{1.100000pt}{1.815000pt}}{0.000000pt}%
\pgfpathmoveto{\pgfqpoint{0.800095in}{2.771967in}}%
\pgfpathlineto{\pgfqpoint{0.879944in}{2.711994in}}%
\pgfpathlineto{\pgfqpoint{0.959792in}{2.704436in}}%
\pgfpathlineto{\pgfqpoint{1.039641in}{2.750986in}}%
\pgfpathlineto{\pgfqpoint{1.119489in}{2.789181in}}%
\pgfpathlineto{\pgfqpoint{1.199337in}{2.544749in}}%
\pgfpathlineto{\pgfqpoint{1.279186in}{2.767500in}}%
\pgfpathlineto{\pgfqpoint{1.359034in}{2.784937in}}%
\pgfpathlineto{\pgfqpoint{1.438883in}{2.612134in}}%
\pgfpathlineto{\pgfqpoint{1.518731in}{2.573983in}}%
\pgfpathlineto{\pgfqpoint{1.598580in}{2.757809in}}%
\pgfpathlineto{\pgfqpoint{1.678428in}{2.664678in}}%
\pgfpathlineto{\pgfqpoint{1.758277in}{2.760653in}}%
\pgfpathlineto{\pgfqpoint{1.838125in}{2.746187in}}%
\pgfpathlineto{\pgfqpoint{1.917974in}{2.369786in}}%
\pgfpathlineto{\pgfqpoint{1.997822in}{2.556873in}}%
\pgfpathlineto{\pgfqpoint{2.077671in}{2.130133in}}%
\pgfpathlineto{\pgfqpoint{2.157519in}{1.730311in}}%
\pgfpathlineto{\pgfqpoint{2.237368in}{2.155042in}}%
\pgfpathlineto{\pgfqpoint{2.317216in}{1.749072in}}%
\pgfpathlineto{\pgfqpoint{2.397065in}{2.068468in}}%
\pgfpathlineto{\pgfqpoint{2.476913in}{1.876646in}}%
\pgfpathlineto{\pgfqpoint{2.556762in}{1.978062in}}%
\pgfpathlineto{\pgfqpoint{2.636610in}{1.929607in}}%
\pgfpathlineto{\pgfqpoint{2.716459in}{1.887936in}}%
\pgfpathlineto{\pgfqpoint{2.796307in}{2.057103in}}%
\pgfpathlineto{\pgfqpoint{2.876156in}{1.797624in}}%
\pgfpathlineto{\pgfqpoint{2.956004in}{2.110517in}}%
\pgfpathlineto{\pgfqpoint{3.035853in}{1.722215in}}%
\pgfpathlineto{\pgfqpoint{3.115701in}{2.236691in}}%
\pgfpathlineto{\pgfqpoint{3.195550in}{1.811096in}}%
\pgfpathlineto{\pgfqpoint{3.275398in}{2.285117in}}%
\pgfpathlineto{\pgfqpoint{3.355247in}{1.897839in}}%
\pgfpathlineto{\pgfqpoint{3.435095in}{2.430162in}}%
\pgfpathlineto{\pgfqpoint{3.514944in}{1.996783in}}%
\pgfpathlineto{\pgfqpoint{3.594792in}{2.577128in}}%
\pgfpathlineto{\pgfqpoint{3.674641in}{2.142084in}}%
\pgfpathlineto{\pgfqpoint{3.754489in}{2.755747in}}%
\pgfpathlineto{\pgfqpoint{3.834338in}{2.352532in}}%
\pgfpathlineto{\pgfqpoint{3.914186in}{2.515255in}}%
\pgfpathlineto{\pgfqpoint{3.994034in}{2.456573in}}%
\pgfpathlineto{\pgfqpoint{4.073883in}{2.393649in}}%
\pgfpathlineto{\pgfqpoint{4.153731in}{2.679088in}}%
\pgfpathlineto{\pgfqpoint{4.233580in}{2.721708in}}%
\pgfpathlineto{\pgfqpoint{4.313428in}{2.929620in}}%
\pgfusepath{stroke}%
\end{pgfscope}%
\begin{pgfscope}%
\pgfpathrectangle{\pgfqpoint{0.540587in}{0.499691in}}{\pgfqpoint{3.952500in}{2.552550in}}%
\pgfusepath{clip}%
\pgfsetbuttcap%
\pgfsetroundjoin%
\pgfsetlinewidth{1.104125pt}%
\definecolor{currentstroke}{rgb}{0.000000,0.000000,0.000000}%
\pgfsetstrokecolor{currentstroke}%
\pgfsetdash{{4.070000pt}{1.760000pt}}{0.000000pt}%
\pgfpathmoveto{\pgfqpoint{0.720247in}{1.714655in}}%
\pgfpathlineto{\pgfqpoint{0.800095in}{1.714655in}}%
\pgfpathlineto{\pgfqpoint{0.879944in}{1.714655in}}%
\pgfpathlineto{\pgfqpoint{0.959792in}{1.714655in}}%
\pgfpathlineto{\pgfqpoint{1.039641in}{1.714655in}}%
\pgfpathlineto{\pgfqpoint{1.119489in}{1.714655in}}%
\pgfpathlineto{\pgfqpoint{1.199337in}{1.714655in}}%
\pgfpathlineto{\pgfqpoint{1.279186in}{1.714655in}}%
\pgfpathlineto{\pgfqpoint{1.359034in}{1.714655in}}%
\pgfpathlineto{\pgfqpoint{1.438883in}{1.714655in}}%
\pgfpathlineto{\pgfqpoint{1.518731in}{1.714655in}}%
\pgfpathlineto{\pgfqpoint{1.598580in}{1.714655in}}%
\pgfpathlineto{\pgfqpoint{1.678428in}{1.714655in}}%
\pgfpathlineto{\pgfqpoint{1.758277in}{1.714655in}}%
\pgfpathlineto{\pgfqpoint{1.838125in}{1.714655in}}%
\pgfpathlineto{\pgfqpoint{1.917974in}{1.714655in}}%
\pgfpathlineto{\pgfqpoint{1.997822in}{1.714655in}}%
\pgfpathlineto{\pgfqpoint{2.077671in}{1.714655in}}%
\pgfpathlineto{\pgfqpoint{2.157519in}{1.714655in}}%
\pgfpathlineto{\pgfqpoint{2.237368in}{1.714655in}}%
\pgfpathlineto{\pgfqpoint{2.317216in}{1.714655in}}%
\pgfpathlineto{\pgfqpoint{2.397065in}{1.714655in}}%
\pgfpathlineto{\pgfqpoint{2.476913in}{1.714655in}}%
\pgfpathlineto{\pgfqpoint{2.556762in}{1.714655in}}%
\pgfpathlineto{\pgfqpoint{2.636610in}{1.714655in}}%
\pgfpathlineto{\pgfqpoint{2.716459in}{1.714655in}}%
\pgfpathlineto{\pgfqpoint{2.796307in}{1.714655in}}%
\pgfpathlineto{\pgfqpoint{2.876156in}{1.714655in}}%
\pgfpathlineto{\pgfqpoint{2.956004in}{1.714655in}}%
\pgfpathlineto{\pgfqpoint{3.035853in}{1.714655in}}%
\pgfpathlineto{\pgfqpoint{3.115701in}{1.714655in}}%
\pgfpathlineto{\pgfqpoint{3.195550in}{1.714655in}}%
\pgfpathlineto{\pgfqpoint{3.275398in}{1.714655in}}%
\pgfpathlineto{\pgfqpoint{3.355247in}{1.714655in}}%
\pgfpathlineto{\pgfqpoint{3.435095in}{1.714655in}}%
\pgfpathlineto{\pgfqpoint{3.514944in}{1.714655in}}%
\pgfpathlineto{\pgfqpoint{3.594792in}{1.714655in}}%
\pgfpathlineto{\pgfqpoint{3.674641in}{1.714655in}}%
\pgfpathlineto{\pgfqpoint{3.754489in}{1.714655in}}%
\pgfpathlineto{\pgfqpoint{3.834338in}{1.714655in}}%
\pgfpathlineto{\pgfqpoint{3.914186in}{1.714655in}}%
\pgfpathlineto{\pgfqpoint{3.994034in}{1.714655in}}%
\pgfpathlineto{\pgfqpoint{4.073883in}{1.714655in}}%
\pgfpathlineto{\pgfqpoint{4.153731in}{1.714655in}}%
\pgfpathlineto{\pgfqpoint{4.233580in}{1.714655in}}%
\pgfpathlineto{\pgfqpoint{4.313428in}{1.714655in}}%
\pgfusepath{stroke}%
\end{pgfscope}%
\begin{pgfscope}%
\pgfsetrectcap%
\pgfsetmiterjoin%
\pgfsetlinewidth{0.803000pt}%
\definecolor{currentstroke}{rgb}{0.000000,0.000000,0.000000}%
\pgfsetstrokecolor{currentstroke}%
\pgfsetdash{}{0pt}%
\pgfpathmoveto{\pgfqpoint{0.540587in}{0.499691in}}%
\pgfpathlineto{\pgfqpoint{0.540587in}{3.052241in}}%
\pgfusepath{stroke}%
\end{pgfscope}%
\begin{pgfscope}%
\pgfsetrectcap%
\pgfsetmiterjoin%
\pgfsetlinewidth{0.803000pt}%
\definecolor{currentstroke}{rgb}{0.000000,0.000000,0.000000}%
\pgfsetstrokecolor{currentstroke}%
\pgfsetdash{}{0pt}%
\pgfpathmoveto{\pgfqpoint{4.493088in}{0.499691in}}%
\pgfpathlineto{\pgfqpoint{4.493088in}{3.052241in}}%
\pgfusepath{stroke}%
\end{pgfscope}%
\begin{pgfscope}%
\pgfsetrectcap%
\pgfsetmiterjoin%
\pgfsetlinewidth{0.803000pt}%
\definecolor{currentstroke}{rgb}{0.000000,0.000000,0.000000}%
\pgfsetstrokecolor{currentstroke}%
\pgfsetdash{}{0pt}%
\pgfpathmoveto{\pgfqpoint{0.540587in}{0.499691in}}%
\pgfpathlineto{\pgfqpoint{4.493088in}{0.499691in}}%
\pgfusepath{stroke}%
\end{pgfscope}%
\begin{pgfscope}%
\pgfsetrectcap%
\pgfsetmiterjoin%
\pgfsetlinewidth{0.803000pt}%
\definecolor{currentstroke}{rgb}{0.000000,0.000000,0.000000}%
\pgfsetstrokecolor{currentstroke}%
\pgfsetdash{}{0pt}%
\pgfpathmoveto{\pgfqpoint{0.540587in}{3.052241in}}%
\pgfpathlineto{\pgfqpoint{4.493088in}{3.052241in}}%
\pgfusepath{stroke}%
\end{pgfscope}%
\begin{pgfscope}%
\pgfsetbuttcap%
\pgfsetmiterjoin%
\definecolor{currentfill}{rgb}{1.000000,1.000000,1.000000}%
\pgfsetfillcolor{currentfill}%
\pgfsetfillopacity{0.800000}%
\pgfsetlinewidth{1.003750pt}%
\definecolor{currentstroke}{rgb}{0.800000,0.800000,0.800000}%
\pgfsetstrokecolor{currentstroke}%
\pgfsetstrokeopacity{0.800000}%
\pgfsetdash{}{0pt}%
\pgfpathmoveto{\pgfqpoint{0.637810in}{0.569136in}}%
\pgfpathlineto{\pgfqpoint{3.013827in}{0.569136in}}%
\pgfpathquadraticcurveto{\pgfqpoint{3.041605in}{0.569136in}}{\pgfqpoint{3.041605in}{0.596913in}}%
\pgfpathlineto{\pgfqpoint{3.041605in}{1.437747in}}%
\pgfpathquadraticcurveto{\pgfqpoint{3.041605in}{1.465524in}}{\pgfqpoint{3.013827in}{1.465524in}}%
\pgfpathlineto{\pgfqpoint{0.637810in}{1.465524in}}%
\pgfpathquadraticcurveto{\pgfqpoint{0.610032in}{1.465524in}}{\pgfqpoint{0.610032in}{1.437747in}}%
\pgfpathlineto{\pgfqpoint{0.610032in}{0.596913in}}%
\pgfpathquadraticcurveto{\pgfqpoint{0.610032in}{0.569136in}}{\pgfqpoint{0.637810in}{0.569136in}}%
\pgfpathlineto{\pgfqpoint{0.637810in}{0.569136in}}%
\pgfpathclose%
\pgfusepath{stroke,fill}%
\end{pgfscope}%
\begin{pgfscope}%
\pgfsetrectcap%
\pgfsetroundjoin%
\pgfsetlinewidth{1.104125pt}%
\definecolor{currentstroke}{rgb}{0.000000,0.000000,0.000000}%
\pgfsetstrokecolor{currentstroke}%
\pgfsetdash{}{0pt}%
\pgfpathmoveto{\pgfqpoint{0.665587in}{1.338579in}}%
\pgfpathlineto{\pgfqpoint{0.804476in}{1.338579in}}%
\pgfpathlineto{\pgfqpoint{0.943365in}{1.338579in}}%
\pgfusepath{stroke}%
\end{pgfscope}%
\begin{pgfscope}%
\pgfsetbuttcap%
\pgfsetroundjoin%
\definecolor{currentfill}{rgb}{0.000000,0.000000,0.000000}%
\pgfsetfillcolor{currentfill}%
\pgfsetfillopacity{0.000000}%
\pgfsetlinewidth{1.003750pt}%
\definecolor{currentstroke}{rgb}{0.000000,0.000000,0.000000}%
\pgfsetstrokecolor{currentstroke}%
\pgfsetdash{}{0pt}%
\pgfsys@defobject{currentmarker}{\pgfqpoint{-0.041667in}{-0.041667in}}{\pgfqpoint{0.041667in}{0.041667in}}{%
\pgfpathmoveto{\pgfqpoint{-0.041667in}{-0.041667in}}%
\pgfpathlineto{\pgfqpoint{0.041667in}{0.041667in}}%
\pgfpathmoveto{\pgfqpoint{-0.041667in}{0.041667in}}%
\pgfpathlineto{\pgfqpoint{0.041667in}{-0.041667in}}%
\pgfusepath{stroke,fill}%
}%
\begin{pgfscope}%
\pgfsys@transformshift{0.804476in}{1.338579in}%
\pgfsys@useobject{currentmarker}{}%
\end{pgfscope}%
\end{pgfscope}%
\begin{pgfscope}%
\definecolor{textcolor}{rgb}{0.000000,0.000000,0.000000}%
\pgfsetstrokecolor{textcolor}%
\pgfsetfillcolor{textcolor}%
\pgftext[x=1.054476in,y=1.289968in,left,base]{\color{textcolor}\rmfamily\fontsize{10.000000}{12.000000}\selectfont \(\displaystyle \norm{\tns{C}_k - D^3 f(\vek{x}_k)}_F / \norm{D^3 f(\vek{x}_k)}_F\)}%
\end{pgfscope}%
\begin{pgfscope}%
\pgfsetbuttcap%
\pgfsetroundjoin%
\pgfsetlinewidth{1.104125pt}%
\definecolor{currentstroke}{rgb}{0.000000,0.000000,0.000000}%
\pgfsetstrokecolor{currentstroke}%
\pgfsetdash{{7.040000pt}{1.760000pt}{1.100000pt}{1.760000pt}}{0.000000pt}%
\pgfpathmoveto{\pgfqpoint{0.665587in}{1.130246in}}%
\pgfpathlineto{\pgfqpoint{0.804476in}{1.130246in}}%
\pgfpathlineto{\pgfqpoint{0.943365in}{1.130246in}}%
\pgfusepath{stroke}%
\end{pgfscope}%
\begin{pgfscope}%
\definecolor{textcolor}{rgb}{0.000000,0.000000,0.000000}%
\pgfsetstrokecolor{textcolor}%
\pgfsetfillcolor{textcolor}%
\pgftext[x=1.054476in,y=1.081635in,left,base]{\color{textcolor}\rmfamily\fontsize{10.000000}{12.000000}\selectfont \(\displaystyle \norm{\vek{x}_k - \vek{x}_*}_2 / \norm{\vek{x}_*}_2\)}%
\end{pgfscope}%
\begin{pgfscope}%
\pgfsetbuttcap%
\pgfsetroundjoin%
\pgfsetlinewidth{1.104125pt}%
\definecolor{currentstroke}{rgb}{0.000000,0.000000,0.000000}%
\pgfsetstrokecolor{currentstroke}%
\pgfsetdash{{1.100000pt}{1.815000pt}}{0.000000pt}%
\pgfpathmoveto{\pgfqpoint{0.665587in}{0.921913in}}%
\pgfpathlineto{\pgfqpoint{0.804476in}{0.921913in}}%
\pgfpathlineto{\pgfqpoint{0.943365in}{0.921913in}}%
\pgfusepath{stroke}%
\end{pgfscope}%
\begin{pgfscope}%
\definecolor{textcolor}{rgb}{0.000000,0.000000,0.000000}%
\pgfsetstrokecolor{textcolor}%
\pgfsetfillcolor{textcolor}%
\pgftext[x=1.054476in,y=0.873302in,left,base]{\color{textcolor}\rmfamily\fontsize{10.000000}{12.000000}\selectfont \(\displaystyle \max \{\norm{\vek{s}_{k-1}}_2, \e_{\text{mach}} / \norm{\vek{s}_{k-1}}_2\}\)}%
\end{pgfscope}%
\begin{pgfscope}%
\pgfsetbuttcap%
\pgfsetroundjoin%
\pgfsetlinewidth{1.104125pt}%
\definecolor{currentstroke}{rgb}{0.000000,0.000000,0.000000}%
\pgfsetstrokecolor{currentstroke}%
\pgfsetdash{{4.070000pt}{1.760000pt}}{0.000000pt}%
\pgfpathmoveto{\pgfqpoint{0.665587in}{0.710233in}}%
\pgfpathlineto{\pgfqpoint{0.804476in}{0.710233in}}%
\pgfpathlineto{\pgfqpoint{0.943365in}{0.710233in}}%
\pgfusepath{stroke}%
\end{pgfscope}%
\begin{pgfscope}%
\definecolor{textcolor}{rgb}{0.000000,0.000000,0.000000}%
\pgfsetstrokecolor{textcolor}%
\pgfsetfillcolor{textcolor}%
\pgftext[x=1.054476in,y=0.661622in,left,base]{\color{textcolor}\rmfamily\fontsize{10.000000}{12.000000}\selectfont \(\displaystyle \sqrt{\e_{\text{mach}}}\)}%
\end{pgfscope}%
\end{pgfpicture}%
\makeatother%
\endgroup%

%% file: plots/convergence_quadratically_converging_min.pgf
\begingroup%
\makeatletter%
\begin{pgfpicture}%
\pgfpathrectangle{\pgfpointorigin}{\pgfqpoint{4.593087in}{3.152241in}}%
\pgfusepath{use as bounding box, clip}%
\begin{pgfscope}%
\pgfsetbuttcap%
\pgfsetmiterjoin%
\definecolor{currentfill}{rgb}{1.000000,1.000000,1.000000}%
\pgfsetfillcolor{currentfill}%
\pgfsetlinewidth{0.000000pt}%
\definecolor{currentstroke}{rgb}{1.000000,1.000000,1.000000}%
\pgfsetstrokecolor{currentstroke}%
\pgfsetdash{}{0pt}%
\pgfpathmoveto{\pgfqpoint{0.000000in}{0.000000in}}%
\pgfpathlineto{\pgfqpoint{4.593087in}{0.000000in}}%
\pgfpathlineto{\pgfqpoint{4.593087in}{3.152241in}}%
\pgfpathlineto{\pgfqpoint{0.000000in}{3.152241in}}%
\pgfpathlineto{\pgfqpoint{0.000000in}{0.000000in}}%
\pgfpathclose%
\pgfusepath{fill}%
\end{pgfscope}%
\begin{pgfscope}%
\pgfsetbuttcap%
\pgfsetmiterjoin%
\definecolor{currentfill}{rgb}{1.000000,1.000000,1.000000}%
\pgfsetfillcolor{currentfill}%
\pgfsetlinewidth{0.000000pt}%
\definecolor{currentstroke}{rgb}{0.000000,0.000000,0.000000}%
\pgfsetstrokecolor{currentstroke}%
\pgfsetstrokeopacity{0.000000}%
\pgfsetdash{}{0pt}%
\pgfpathmoveto{\pgfqpoint{0.540587in}{0.499691in}}%
\pgfpathlineto{\pgfqpoint{4.493088in}{0.499691in}}%
\pgfpathlineto{\pgfqpoint{4.493088in}{3.052241in}}%
\pgfpathlineto{\pgfqpoint{0.540587in}{3.052241in}}%
\pgfpathlineto{\pgfqpoint{0.540587in}{0.499691in}}%
\pgfpathclose%
\pgfusepath{fill}%
\end{pgfscope}%
\begin{pgfscope}%
\pgfsetbuttcap%
\pgfsetroundjoin%
\definecolor{currentfill}{rgb}{0.000000,0.000000,0.000000}%
\pgfsetfillcolor{currentfill}%
\pgfsetfillopacity{0.000000}%
\pgfsetlinewidth{0.803000pt}%
\definecolor{currentstroke}{rgb}{0.000000,0.000000,0.000000}%
\pgfsetstrokecolor{currentstroke}%
\pgfsetdash{}{0pt}%
\pgfsys@defobject{currentmarker}{\pgfqpoint{0.000000in}{-0.048611in}}{\pgfqpoint{0.000000in}{0.000000in}}{%
\pgfpathmoveto{\pgfqpoint{0.000000in}{0.000000in}}%
\pgfpathlineto{\pgfqpoint{0.000000in}{-0.048611in}}%
\pgfusepath{stroke,fill}%
}%
\begin{pgfscope}%
\pgfsys@transformshift{0.720247in}{0.499691in}%
\pgfsys@useobject{currentmarker}{}%
\end{pgfscope}%
\end{pgfscope}%
\begin{pgfscope}%
\definecolor{textcolor}{rgb}{0.000000,0.000000,0.000000}%
\pgfsetstrokecolor{textcolor}%
\pgfsetfillcolor{textcolor}%
\pgftext[x=0.720247in,y=0.402469in,,top]{\color{textcolor}\rmfamily\fontsize{10.000000}{12.000000}\selectfont \(\displaystyle {0}\)}%
\end{pgfscope}%
\begin{pgfscope}%
\pgfsetbuttcap%
\pgfsetroundjoin%
\definecolor{currentfill}{rgb}{0.000000,0.000000,0.000000}%
\pgfsetfillcolor{currentfill}%
\pgfsetfillopacity{0.000000}%
\pgfsetlinewidth{0.803000pt}%
\definecolor{currentstroke}{rgb}{0.000000,0.000000,0.000000}%
\pgfsetstrokecolor{currentstroke}%
\pgfsetdash{}{0pt}%
\pgfsys@defobject{currentmarker}{\pgfqpoint{0.000000in}{-0.048611in}}{\pgfqpoint{0.000000in}{0.000000in}}{%
\pgfpathmoveto{\pgfqpoint{0.000000in}{0.000000in}}%
\pgfpathlineto{\pgfqpoint{0.000000in}{-0.048611in}}%
\pgfusepath{stroke,fill}%
}%
\begin{pgfscope}%
\pgfsys@transformshift{1.169394in}{0.499691in}%
\pgfsys@useobject{currentmarker}{}%
\end{pgfscope}%
\end{pgfscope}%
\begin{pgfscope}%
\definecolor{textcolor}{rgb}{0.000000,0.000000,0.000000}%
\pgfsetstrokecolor{textcolor}%
\pgfsetfillcolor{textcolor}%
\pgftext[x=1.169394in,y=0.402469in,,top]{\color{textcolor}\rmfamily\fontsize{10.000000}{12.000000}\selectfont \(\displaystyle {2}\)}%
\end{pgfscope}%
\begin{pgfscope}%
\pgfsetbuttcap%
\pgfsetroundjoin%
\definecolor{currentfill}{rgb}{0.000000,0.000000,0.000000}%
\pgfsetfillcolor{currentfill}%
\pgfsetfillopacity{0.000000}%
\pgfsetlinewidth{0.803000pt}%
\definecolor{currentstroke}{rgb}{0.000000,0.000000,0.000000}%
\pgfsetstrokecolor{currentstroke}%
\pgfsetdash{}{0pt}%
\pgfsys@defobject{currentmarker}{\pgfqpoint{0.000000in}{-0.048611in}}{\pgfqpoint{0.000000in}{0.000000in}}{%
\pgfpathmoveto{\pgfqpoint{0.000000in}{0.000000in}}%
\pgfpathlineto{\pgfqpoint{0.000000in}{-0.048611in}}%
\pgfusepath{stroke,fill}%
}%
\begin{pgfscope}%
\pgfsys@transformshift{1.618542in}{0.499691in}%
\pgfsys@useobject{currentmarker}{}%
\end{pgfscope}%
\end{pgfscope}%
\begin{pgfscope}%
\definecolor{textcolor}{rgb}{0.000000,0.000000,0.000000}%
\pgfsetstrokecolor{textcolor}%
\pgfsetfillcolor{textcolor}%
\pgftext[x=1.618542in,y=0.402469in,,top]{\color{textcolor}\rmfamily\fontsize{10.000000}{12.000000}\selectfont \(\displaystyle {4}\)}%
\end{pgfscope}%
\begin{pgfscope}%
\pgfsetbuttcap%
\pgfsetroundjoin%
\definecolor{currentfill}{rgb}{0.000000,0.000000,0.000000}%
\pgfsetfillcolor{currentfill}%
\pgfsetfillopacity{0.000000}%
\pgfsetlinewidth{0.803000pt}%
\definecolor{currentstroke}{rgb}{0.000000,0.000000,0.000000}%
\pgfsetstrokecolor{currentstroke}%
\pgfsetdash{}{0pt}%
\pgfsys@defobject{currentmarker}{\pgfqpoint{0.000000in}{-0.048611in}}{\pgfqpoint{0.000000in}{0.000000in}}{%
\pgfpathmoveto{\pgfqpoint{0.000000in}{0.000000in}}%
\pgfpathlineto{\pgfqpoint{0.000000in}{-0.048611in}}%
\pgfusepath{stroke,fill}%
}%
\begin{pgfscope}%
\pgfsys@transformshift{2.067690in}{0.499691in}%
\pgfsys@useobject{currentmarker}{}%
\end{pgfscope}%
\end{pgfscope}%
\begin{pgfscope}%
\definecolor{textcolor}{rgb}{0.000000,0.000000,0.000000}%
\pgfsetstrokecolor{textcolor}%
\pgfsetfillcolor{textcolor}%
\pgftext[x=2.067690in,y=0.402469in,,top]{\color{textcolor}\rmfamily\fontsize{10.000000}{12.000000}\selectfont \(\displaystyle {6}\)}%
\end{pgfscope}%
\begin{pgfscope}%
\pgfsetbuttcap%
\pgfsetroundjoin%
\definecolor{currentfill}{rgb}{0.000000,0.000000,0.000000}%
\pgfsetfillcolor{currentfill}%
\pgfsetfillopacity{0.000000}%
\pgfsetlinewidth{0.803000pt}%
\definecolor{currentstroke}{rgb}{0.000000,0.000000,0.000000}%
\pgfsetstrokecolor{currentstroke}%
\pgfsetdash{}{0pt}%
\pgfsys@defobject{currentmarker}{\pgfqpoint{0.000000in}{-0.048611in}}{\pgfqpoint{0.000000in}{0.000000in}}{%
\pgfpathmoveto{\pgfqpoint{0.000000in}{0.000000in}}%
\pgfpathlineto{\pgfqpoint{0.000000in}{-0.048611in}}%
\pgfusepath{stroke,fill}%
}%
\begin{pgfscope}%
\pgfsys@transformshift{2.516838in}{0.499691in}%
\pgfsys@useobject{currentmarker}{}%
\end{pgfscope}%
\end{pgfscope}%
\begin{pgfscope}%
\definecolor{textcolor}{rgb}{0.000000,0.000000,0.000000}%
\pgfsetstrokecolor{textcolor}%
\pgfsetfillcolor{textcolor}%
\pgftext[x=2.516837in,y=0.402469in,,top]{\color{textcolor}\rmfamily\fontsize{10.000000}{12.000000}\selectfont \(\displaystyle {8}\)}%
\end{pgfscope}%
\begin{pgfscope}%
\pgfsetbuttcap%
\pgfsetroundjoin%
\definecolor{currentfill}{rgb}{0.000000,0.000000,0.000000}%
\pgfsetfillcolor{currentfill}%
\pgfsetfillopacity{0.000000}%
\pgfsetlinewidth{0.803000pt}%
\definecolor{currentstroke}{rgb}{0.000000,0.000000,0.000000}%
\pgfsetstrokecolor{currentstroke}%
\pgfsetdash{}{0pt}%
\pgfsys@defobject{currentmarker}{\pgfqpoint{0.000000in}{-0.048611in}}{\pgfqpoint{0.000000in}{0.000000in}}{%
\pgfpathmoveto{\pgfqpoint{0.000000in}{0.000000in}}%
\pgfpathlineto{\pgfqpoint{0.000000in}{-0.048611in}}%
\pgfusepath{stroke,fill}%
}%
\begin{pgfscope}%
\pgfsys@transformshift{2.965985in}{0.499691in}%
\pgfsys@useobject{currentmarker}{}%
\end{pgfscope}%
\end{pgfscope}%
\begin{pgfscope}%
\definecolor{textcolor}{rgb}{0.000000,0.000000,0.000000}%
\pgfsetstrokecolor{textcolor}%
\pgfsetfillcolor{textcolor}%
\pgftext[x=2.965985in,y=0.402469in,,top]{\color{textcolor}\rmfamily\fontsize{10.000000}{12.000000}\selectfont \(\displaystyle {10}\)}%
\end{pgfscope}%
\begin{pgfscope}%
\pgfsetbuttcap%
\pgfsetroundjoin%
\definecolor{currentfill}{rgb}{0.000000,0.000000,0.000000}%
\pgfsetfillcolor{currentfill}%
\pgfsetfillopacity{0.000000}%
\pgfsetlinewidth{0.803000pt}%
\definecolor{currentstroke}{rgb}{0.000000,0.000000,0.000000}%
\pgfsetstrokecolor{currentstroke}%
\pgfsetdash{}{0pt}%
\pgfsys@defobject{currentmarker}{\pgfqpoint{0.000000in}{-0.048611in}}{\pgfqpoint{0.000000in}{0.000000in}}{%
\pgfpathmoveto{\pgfqpoint{0.000000in}{0.000000in}}%
\pgfpathlineto{\pgfqpoint{0.000000in}{-0.048611in}}%
\pgfusepath{stroke,fill}%
}%
\begin{pgfscope}%
\pgfsys@transformshift{3.415133in}{0.499691in}%
\pgfsys@useobject{currentmarker}{}%
\end{pgfscope}%
\end{pgfscope}%
\begin{pgfscope}%
\definecolor{textcolor}{rgb}{0.000000,0.000000,0.000000}%
\pgfsetstrokecolor{textcolor}%
\pgfsetfillcolor{textcolor}%
\pgftext[x=3.415133in,y=0.402469in,,top]{\color{textcolor}\rmfamily\fontsize{10.000000}{12.000000}\selectfont \(\displaystyle {12}\)}%
\end{pgfscope}%
\begin{pgfscope}%
\pgfsetbuttcap%
\pgfsetroundjoin%
\definecolor{currentfill}{rgb}{0.000000,0.000000,0.000000}%
\pgfsetfillcolor{currentfill}%
\pgfsetfillopacity{0.000000}%
\pgfsetlinewidth{0.803000pt}%
\definecolor{currentstroke}{rgb}{0.000000,0.000000,0.000000}%
\pgfsetstrokecolor{currentstroke}%
\pgfsetdash{}{0pt}%
\pgfsys@defobject{currentmarker}{\pgfqpoint{0.000000in}{-0.048611in}}{\pgfqpoint{0.000000in}{0.000000in}}{%
\pgfpathmoveto{\pgfqpoint{0.000000in}{0.000000in}}%
\pgfpathlineto{\pgfqpoint{0.000000in}{-0.048611in}}%
\pgfusepath{stroke,fill}%
}%
\begin{pgfscope}%
\pgfsys@transformshift{3.864281in}{0.499691in}%
\pgfsys@useobject{currentmarker}{}%
\end{pgfscope}%
\end{pgfscope}%
\begin{pgfscope}%
\definecolor{textcolor}{rgb}{0.000000,0.000000,0.000000}%
\pgfsetstrokecolor{textcolor}%
\pgfsetfillcolor{textcolor}%
\pgftext[x=3.864281in,y=0.402469in,,top]{\color{textcolor}\rmfamily\fontsize{10.000000}{12.000000}\selectfont \(\displaystyle {14}\)}%
\end{pgfscope}%
\begin{pgfscope}%
\pgfsetbuttcap%
\pgfsetroundjoin%
\definecolor{currentfill}{rgb}{0.000000,0.000000,0.000000}%
\pgfsetfillcolor{currentfill}%
\pgfsetfillopacity{0.000000}%
\pgfsetlinewidth{0.803000pt}%
\definecolor{currentstroke}{rgb}{0.000000,0.000000,0.000000}%
\pgfsetstrokecolor{currentstroke}%
\pgfsetdash{}{0pt}%
\pgfsys@defobject{currentmarker}{\pgfqpoint{0.000000in}{-0.048611in}}{\pgfqpoint{0.000000in}{0.000000in}}{%
\pgfpathmoveto{\pgfqpoint{0.000000in}{0.000000in}}%
\pgfpathlineto{\pgfqpoint{0.000000in}{-0.048611in}}%
\pgfusepath{stroke,fill}%
}%
\begin{pgfscope}%
\pgfsys@transformshift{4.313428in}{0.499691in}%
\pgfsys@useobject{currentmarker}{}%
\end{pgfscope}%
\end{pgfscope}%
\begin{pgfscope}%
\definecolor{textcolor}{rgb}{0.000000,0.000000,0.000000}%
\pgfsetstrokecolor{textcolor}%
\pgfsetfillcolor{textcolor}%
\pgftext[x=4.313428in,y=0.402469in,,top]{\color{textcolor}\rmfamily\fontsize{10.000000}{12.000000}\selectfont \(\displaystyle {16}\)}%
\end{pgfscope}%
\begin{pgfscope}%
\definecolor{textcolor}{rgb}{0.000000,0.000000,0.000000}%
\pgfsetstrokecolor{textcolor}%
\pgfsetfillcolor{textcolor}%
\pgftext[x=2.516837in,y=0.223457in,,top]{\color{textcolor}\rmfamily\fontsize{10.000000}{12.000000}\selectfont Iteration \(\displaystyle k\)}%
\end{pgfscope}%
\begin{pgfscope}%
\pgfsetbuttcap%
\pgfsetroundjoin%
\definecolor{currentfill}{rgb}{0.000000,0.000000,0.000000}%
\pgfsetfillcolor{currentfill}%
\pgfsetfillopacity{0.000000}%
\pgfsetlinewidth{0.803000pt}%
\definecolor{currentstroke}{rgb}{0.000000,0.000000,0.000000}%
\pgfsetstrokecolor{currentstroke}%
\pgfsetdash{}{0pt}%
\pgfsys@defobject{currentmarker}{\pgfqpoint{-0.048611in}{0.000000in}}{\pgfqpoint{-0.000000in}{0.000000in}}{%
\pgfpathmoveto{\pgfqpoint{-0.000000in}{0.000000in}}%
\pgfpathlineto{\pgfqpoint{-0.048611in}{0.000000in}}%
\pgfusepath{stroke,fill}%
}%
\begin{pgfscope}%
\pgfsys@transformshift{0.540587in}{0.801300in}%
\pgfsys@useobject{currentmarker}{}%
\end{pgfscope}%
\end{pgfscope}%
\begin{pgfscope}%
\definecolor{textcolor}{rgb}{0.000000,0.000000,0.000000}%
\pgfsetstrokecolor{textcolor}%
\pgfsetfillcolor{textcolor}%
\pgftext[x=0.100000in, y=0.753074in, left, base]{\color{textcolor}\rmfamily\fontsize{10.000000}{12.000000}\selectfont \(\displaystyle {10^{-14}}\)}%
\end{pgfscope}%
\begin{pgfscope}%
\pgfsetbuttcap%
\pgfsetroundjoin%
\definecolor{currentfill}{rgb}{0.000000,0.000000,0.000000}%
\pgfsetfillcolor{currentfill}%
\pgfsetfillopacity{0.000000}%
\pgfsetlinewidth{0.803000pt}%
\definecolor{currentstroke}{rgb}{0.000000,0.000000,0.000000}%
\pgfsetstrokecolor{currentstroke}%
\pgfsetdash{}{0pt}%
\pgfsys@defobject{currentmarker}{\pgfqpoint{-0.048611in}{0.000000in}}{\pgfqpoint{-0.000000in}{0.000000in}}{%
\pgfpathmoveto{\pgfqpoint{-0.000000in}{0.000000in}}%
\pgfpathlineto{\pgfqpoint{-0.048611in}{0.000000in}}%
\pgfusepath{stroke,fill}%
}%
\begin{pgfscope}%
\pgfsys@transformshift{0.540587in}{1.109915in}%
\pgfsys@useobject{currentmarker}{}%
\end{pgfscope}%
\end{pgfscope}%
\begin{pgfscope}%
\definecolor{textcolor}{rgb}{0.000000,0.000000,0.000000}%
\pgfsetstrokecolor{textcolor}%
\pgfsetfillcolor{textcolor}%
\pgftext[x=0.100000in, y=1.061690in, left, base]{\color{textcolor}\rmfamily\fontsize{10.000000}{12.000000}\selectfont \(\displaystyle {10^{-12}}\)}%
\end{pgfscope}%
\begin{pgfscope}%
\pgfsetbuttcap%
\pgfsetroundjoin%
\definecolor{currentfill}{rgb}{0.000000,0.000000,0.000000}%
\pgfsetfillcolor{currentfill}%
\pgfsetfillopacity{0.000000}%
\pgfsetlinewidth{0.803000pt}%
\definecolor{currentstroke}{rgb}{0.000000,0.000000,0.000000}%
\pgfsetstrokecolor{currentstroke}%
\pgfsetdash{}{0pt}%
\pgfsys@defobject{currentmarker}{\pgfqpoint{-0.048611in}{0.000000in}}{\pgfqpoint{-0.000000in}{0.000000in}}{%
\pgfpathmoveto{\pgfqpoint{-0.000000in}{0.000000in}}%
\pgfpathlineto{\pgfqpoint{-0.048611in}{0.000000in}}%
\pgfusepath{stroke,fill}%
}%
\begin{pgfscope}%
\pgfsys@transformshift{0.540587in}{1.418531in}%
\pgfsys@useobject{currentmarker}{}%
\end{pgfscope}%
\end{pgfscope}%
\begin{pgfscope}%
\definecolor{textcolor}{rgb}{0.000000,0.000000,0.000000}%
\pgfsetstrokecolor{textcolor}%
\pgfsetfillcolor{textcolor}%
\pgftext[x=0.100000in, y=1.370306in, left, base]{\color{textcolor}\rmfamily\fontsize{10.000000}{12.000000}\selectfont \(\displaystyle {10^{-10}}\)}%
\end{pgfscope}%
\begin{pgfscope}%
\pgfsetbuttcap%
\pgfsetroundjoin%
\definecolor{currentfill}{rgb}{0.000000,0.000000,0.000000}%
\pgfsetfillcolor{currentfill}%
\pgfsetfillopacity{0.000000}%
\pgfsetlinewidth{0.803000pt}%
\definecolor{currentstroke}{rgb}{0.000000,0.000000,0.000000}%
\pgfsetstrokecolor{currentstroke}%
\pgfsetdash{}{0pt}%
\pgfsys@defobject{currentmarker}{\pgfqpoint{-0.048611in}{0.000000in}}{\pgfqpoint{-0.000000in}{0.000000in}}{%
\pgfpathmoveto{\pgfqpoint{-0.000000in}{0.000000in}}%
\pgfpathlineto{\pgfqpoint{-0.048611in}{0.000000in}}%
\pgfusepath{stroke,fill}%
}%
\begin{pgfscope}%
\pgfsys@transformshift{0.540587in}{1.727147in}%
\pgfsys@useobject{currentmarker}{}%
\end{pgfscope}%
\end{pgfscope}%
\begin{pgfscope}%
\definecolor{textcolor}{rgb}{0.000000,0.000000,0.000000}%
\pgfsetstrokecolor{textcolor}%
\pgfsetfillcolor{textcolor}%
\pgftext[x=0.155363in, y=1.678921in, left, base]{\color{textcolor}\rmfamily\fontsize{10.000000}{12.000000}\selectfont \(\displaystyle {10^{-8}}\)}%
\end{pgfscope}%
\begin{pgfscope}%
\pgfsetbuttcap%
\pgfsetroundjoin%
\definecolor{currentfill}{rgb}{0.000000,0.000000,0.000000}%
\pgfsetfillcolor{currentfill}%
\pgfsetfillopacity{0.000000}%
\pgfsetlinewidth{0.803000pt}%
\definecolor{currentstroke}{rgb}{0.000000,0.000000,0.000000}%
\pgfsetstrokecolor{currentstroke}%
\pgfsetdash{}{0pt}%
\pgfsys@defobject{currentmarker}{\pgfqpoint{-0.048611in}{0.000000in}}{\pgfqpoint{-0.000000in}{0.000000in}}{%
\pgfpathmoveto{\pgfqpoint{-0.000000in}{0.000000in}}%
\pgfpathlineto{\pgfqpoint{-0.048611in}{0.000000in}}%
\pgfusepath{stroke,fill}%
}%
\begin{pgfscope}%
\pgfsys@transformshift{0.540587in}{2.035762in}%
\pgfsys@useobject{currentmarker}{}%
\end{pgfscope}%
\end{pgfscope}%
\begin{pgfscope}%
\definecolor{textcolor}{rgb}{0.000000,0.000000,0.000000}%
\pgfsetstrokecolor{textcolor}%
\pgfsetfillcolor{textcolor}%
\pgftext[x=0.155363in, y=1.987537in, left, base]{\color{textcolor}\rmfamily\fontsize{10.000000}{12.000000}\selectfont \(\displaystyle {10^{-6}}\)}%
\end{pgfscope}%
\begin{pgfscope}%
\pgfsetbuttcap%
\pgfsetroundjoin%
\definecolor{currentfill}{rgb}{0.000000,0.000000,0.000000}%
\pgfsetfillcolor{currentfill}%
\pgfsetfillopacity{0.000000}%
\pgfsetlinewidth{0.803000pt}%
\definecolor{currentstroke}{rgb}{0.000000,0.000000,0.000000}%
\pgfsetstrokecolor{currentstroke}%
\pgfsetdash{}{0pt}%
\pgfsys@defobject{currentmarker}{\pgfqpoint{-0.048611in}{0.000000in}}{\pgfqpoint{-0.000000in}{0.000000in}}{%
\pgfpathmoveto{\pgfqpoint{-0.000000in}{0.000000in}}%
\pgfpathlineto{\pgfqpoint{-0.048611in}{0.000000in}}%
\pgfusepath{stroke,fill}%
}%
\begin{pgfscope}%
\pgfsys@transformshift{0.540587in}{2.344378in}%
\pgfsys@useobject{currentmarker}{}%
\end{pgfscope}%
\end{pgfscope}%
\begin{pgfscope}%
\definecolor{textcolor}{rgb}{0.000000,0.000000,0.000000}%
\pgfsetstrokecolor{textcolor}%
\pgfsetfillcolor{textcolor}%
\pgftext[x=0.155363in, y=2.296153in, left, base]{\color{textcolor}\rmfamily\fontsize{10.000000}{12.000000}\selectfont \(\displaystyle {10^{-4}}\)}%
\end{pgfscope}%
\begin{pgfscope}%
\pgfsetbuttcap%
\pgfsetroundjoin%
\definecolor{currentfill}{rgb}{0.000000,0.000000,0.000000}%
\pgfsetfillcolor{currentfill}%
\pgfsetfillopacity{0.000000}%
\pgfsetlinewidth{0.803000pt}%
\definecolor{currentstroke}{rgb}{0.000000,0.000000,0.000000}%
\pgfsetstrokecolor{currentstroke}%
\pgfsetdash{}{0pt}%
\pgfsys@defobject{currentmarker}{\pgfqpoint{-0.048611in}{0.000000in}}{\pgfqpoint{-0.000000in}{0.000000in}}{%
\pgfpathmoveto{\pgfqpoint{-0.000000in}{0.000000in}}%
\pgfpathlineto{\pgfqpoint{-0.048611in}{0.000000in}}%
\pgfusepath{stroke,fill}%
}%
\begin{pgfscope}%
\pgfsys@transformshift{0.540587in}{2.652993in}%
\pgfsys@useobject{currentmarker}{}%
\end{pgfscope}%
\end{pgfscope}%
\begin{pgfscope}%
\definecolor{textcolor}{rgb}{0.000000,0.000000,0.000000}%
\pgfsetstrokecolor{textcolor}%
\pgfsetfillcolor{textcolor}%
\pgftext[x=0.155363in, y=2.604768in, left, base]{\color{textcolor}\rmfamily\fontsize{10.000000}{12.000000}\selectfont \(\displaystyle {10^{-2}}\)}%
\end{pgfscope}%
\begin{pgfscope}%
\pgfsetbuttcap%
\pgfsetroundjoin%
\definecolor{currentfill}{rgb}{0.000000,0.000000,0.000000}%
\pgfsetfillcolor{currentfill}%
\pgfsetfillopacity{0.000000}%
\pgfsetlinewidth{0.803000pt}%
\definecolor{currentstroke}{rgb}{0.000000,0.000000,0.000000}%
\pgfsetstrokecolor{currentstroke}%
\pgfsetdash{}{0pt}%
\pgfsys@defobject{currentmarker}{\pgfqpoint{-0.048611in}{0.000000in}}{\pgfqpoint{-0.000000in}{0.000000in}}{%
\pgfpathmoveto{\pgfqpoint{-0.000000in}{0.000000in}}%
\pgfpathlineto{\pgfqpoint{-0.048611in}{0.000000in}}%
\pgfusepath{stroke,fill}%
}%
\begin{pgfscope}%
\pgfsys@transformshift{0.540587in}{2.961609in}%
\pgfsys@useobject{currentmarker}{}%
\end{pgfscope}%
\end{pgfscope}%
\begin{pgfscope}%
\definecolor{textcolor}{rgb}{0.000000,0.000000,0.000000}%
\pgfsetstrokecolor{textcolor}%
\pgfsetfillcolor{textcolor}%
\pgftext[x=0.242169in, y=2.913384in, left, base]{\color{textcolor}\rmfamily\fontsize{10.000000}{12.000000}\selectfont \(\displaystyle {10^{0}}\)}%
\end{pgfscope}%
\begin{pgfscope}%
\pgfpathrectangle{\pgfqpoint{0.540587in}{0.499691in}}{\pgfqpoint{3.952500in}{2.552550in}}%
\pgfusepath{clip}%
\pgfsetrectcap%
\pgfsetroundjoin%
\pgfsetlinewidth{1.104125pt}%
\definecolor{currentstroke}{rgb}{0.000000,0.000000,0.000000}%
\pgfsetstrokecolor{currentstroke}%
\pgfsetdash{}{0pt}%
\pgfpathmoveto{\pgfqpoint{0.720247in}{2.961609in}}%
\pgfpathlineto{\pgfqpoint{0.944820in}{2.811081in}}%
\pgfpathlineto{\pgfqpoint{1.169394in}{2.846700in}}%
\pgfpathlineto{\pgfqpoint{1.393968in}{2.855106in}}%
\pgfpathlineto{\pgfqpoint{1.618542in}{2.845315in}}%
\pgfpathlineto{\pgfqpoint{1.843116in}{2.864380in}}%
\pgfpathlineto{\pgfqpoint{2.067690in}{2.870333in}}%
\pgfpathlineto{\pgfqpoint{2.292264in}{2.880789in}}%
\pgfpathlineto{\pgfqpoint{2.516838in}{2.883323in}}%
\pgfpathlineto{\pgfqpoint{2.741411in}{2.884061in}}%
\pgfpathlineto{\pgfqpoint{2.965985in}{2.888828in}}%
\pgfpathlineto{\pgfqpoint{3.190559in}{2.890888in}}%
\pgfpathlineto{\pgfqpoint{3.415133in}{2.892110in}}%
\pgfpathlineto{\pgfqpoint{3.639707in}{2.892269in}}%
\pgfpathlineto{\pgfqpoint{3.864281in}{2.892250in}}%
\pgfpathlineto{\pgfqpoint{4.088855in}{2.892222in}}%
\pgfpathlineto{\pgfqpoint{4.313428in}{2.892193in}}%
\pgfusepath{stroke}%
\end{pgfscope}%
\begin{pgfscope}%
\pgfpathrectangle{\pgfqpoint{0.540587in}{0.499691in}}{\pgfqpoint{3.952500in}{2.552550in}}%
\pgfusepath{clip}%
\pgfsetbuttcap%
\pgfsetroundjoin%
\definecolor{currentfill}{rgb}{0.000000,0.000000,0.000000}%
\pgfsetfillcolor{currentfill}%
\pgfsetfillopacity{0.000000}%
\pgfsetlinewidth{1.003750pt}%
\definecolor{currentstroke}{rgb}{0.000000,0.000000,0.000000}%
\pgfsetstrokecolor{currentstroke}%
\pgfsetdash{}{0pt}%
\pgfsys@defobject{currentmarker}{\pgfqpoint{-0.041667in}{-0.041667in}}{\pgfqpoint{0.041667in}{0.041667in}}{%
\pgfpathmoveto{\pgfqpoint{-0.041667in}{-0.041667in}}%
\pgfpathlineto{\pgfqpoint{0.041667in}{0.041667in}}%
\pgfpathmoveto{\pgfqpoint{-0.041667in}{0.041667in}}%
\pgfpathlineto{\pgfqpoint{0.041667in}{-0.041667in}}%
\pgfusepath{stroke,fill}%
}%
\begin{pgfscope}%
\pgfsys@transformshift{0.720247in}{2.961609in}%
\pgfsys@useobject{currentmarker}{}%
\end{pgfscope}%
\begin{pgfscope}%
\pgfsys@transformshift{0.944820in}{2.811081in}%
\pgfsys@useobject{currentmarker}{}%
\end{pgfscope}%
\begin{pgfscope}%
\pgfsys@transformshift{1.169394in}{2.846700in}%
\pgfsys@useobject{currentmarker}{}%
\end{pgfscope}%
\begin{pgfscope}%
\pgfsys@transformshift{1.393968in}{2.855106in}%
\pgfsys@useobject{currentmarker}{}%
\end{pgfscope}%
\begin{pgfscope}%
\pgfsys@transformshift{1.618542in}{2.845315in}%
\pgfsys@useobject{currentmarker}{}%
\end{pgfscope}%
\begin{pgfscope}%
\pgfsys@transformshift{1.843116in}{2.864380in}%
\pgfsys@useobject{currentmarker}{}%
\end{pgfscope}%
\begin{pgfscope}%
\pgfsys@transformshift{2.067690in}{2.870333in}%
\pgfsys@useobject{currentmarker}{}%
\end{pgfscope}%
\begin{pgfscope}%
\pgfsys@transformshift{2.292264in}{2.880789in}%
\pgfsys@useobject{currentmarker}{}%
\end{pgfscope}%
\begin{pgfscope}%
\pgfsys@transformshift{2.516838in}{2.883323in}%
\pgfsys@useobject{currentmarker}{}%
\end{pgfscope}%
\begin{pgfscope}%
\pgfsys@transformshift{2.741411in}{2.884061in}%
\pgfsys@useobject{currentmarker}{}%
\end{pgfscope}%
\begin{pgfscope}%
\pgfsys@transformshift{2.965985in}{2.888828in}%
\pgfsys@useobject{currentmarker}{}%
\end{pgfscope}%
\begin{pgfscope}%
\pgfsys@transformshift{3.190559in}{2.890888in}%
\pgfsys@useobject{currentmarker}{}%
\end{pgfscope}%
\begin{pgfscope}%
\pgfsys@transformshift{3.415133in}{2.892110in}%
\pgfsys@useobject{currentmarker}{}%
\end{pgfscope}%
\begin{pgfscope}%
\pgfsys@transformshift{3.639707in}{2.892269in}%
\pgfsys@useobject{currentmarker}{}%
\end{pgfscope}%
\begin{pgfscope}%
\pgfsys@transformshift{3.864281in}{2.892250in}%
\pgfsys@useobject{currentmarker}{}%
\end{pgfscope}%
\begin{pgfscope}%
\pgfsys@transformshift{4.088855in}{2.892222in}%
\pgfsys@useobject{currentmarker}{}%
\end{pgfscope}%
\begin{pgfscope}%
\pgfsys@transformshift{4.313428in}{2.892193in}%
\pgfsys@useobject{currentmarker}{}%
\end{pgfscope}%
\end{pgfscope}%
\begin{pgfscope}%
\pgfpathrectangle{\pgfqpoint{0.540587in}{0.499691in}}{\pgfqpoint{3.952500in}{2.552550in}}%
\pgfusepath{clip}%
\pgfsetrectcap%
\pgfsetroundjoin%
\pgfsetlinewidth{1.104125pt}%
\definecolor{currentstroke}{rgb}{0.000000,0.000000,0.000000}%
\pgfsetstrokecolor{currentstroke}%
\pgfsetdash{}{0pt}%
\pgfpathmoveto{\pgfqpoint{0.720247in}{2.948023in}}%
\pgfpathlineto{\pgfqpoint{0.944820in}{2.810353in}}%
\pgfpathlineto{\pgfqpoint{1.169394in}{2.840315in}}%
\pgfpathlineto{\pgfqpoint{1.393968in}{2.822186in}}%
\pgfpathlineto{\pgfqpoint{1.618542in}{2.777739in}}%
\pgfpathlineto{\pgfqpoint{1.843116in}{2.771821in}}%
\pgfpathlineto{\pgfqpoint{2.067690in}{2.766750in}}%
\pgfpathlineto{\pgfqpoint{2.292264in}{2.749723in}}%
\pgfpathlineto{\pgfqpoint{2.516838in}{2.751560in}}%
\pgfpathlineto{\pgfqpoint{2.741411in}{2.612445in}}%
\pgfpathlineto{\pgfqpoint{2.965985in}{2.655978in}}%
\pgfpathlineto{\pgfqpoint{3.190559in}{2.642148in}}%
\pgfpathlineto{\pgfqpoint{3.415133in}{2.615803in}}%
\pgfpathlineto{\pgfqpoint{3.639707in}{2.620383in}}%
\pgfpathlineto{\pgfqpoint{3.864281in}{2.619014in}}%
\pgfpathlineto{\pgfqpoint{4.088855in}{2.619081in}}%
\pgfusepath{stroke}%
\end{pgfscope}%
\begin{pgfscope}%
\pgfpathrectangle{\pgfqpoint{0.540587in}{0.499691in}}{\pgfqpoint{3.952500in}{2.552550in}}%
\pgfusepath{clip}%
\pgfsetbuttcap%
\pgfsetroundjoin%
\definecolor{currentfill}{rgb}{0.000000,0.000000,0.000000}%
\pgfsetfillcolor{currentfill}%
\pgfsetfillopacity{0.000000}%
\pgfsetlinewidth{1.003750pt}%
\definecolor{currentstroke}{rgb}{0.000000,0.000000,0.000000}%
\pgfsetstrokecolor{currentstroke}%
\pgfsetdash{}{0pt}%
\pgfsys@defobject{currentmarker}{\pgfqpoint{-0.041667in}{-0.041667in}}{\pgfqpoint{0.041667in}{0.041667in}}{%
\pgfpathmoveto{\pgfqpoint{0.000000in}{-0.041667in}}%
\pgfpathcurveto{\pgfqpoint{0.011050in}{-0.041667in}}{\pgfqpoint{0.021649in}{-0.037276in}}{\pgfqpoint{0.029463in}{-0.029463in}}%
\pgfpathcurveto{\pgfqpoint{0.037276in}{-0.021649in}}{\pgfqpoint{0.041667in}{-0.011050in}}{\pgfqpoint{0.041667in}{0.000000in}}%
\pgfpathcurveto{\pgfqpoint{0.041667in}{0.011050in}}{\pgfqpoint{0.037276in}{0.021649in}}{\pgfqpoint{0.029463in}{0.029463in}}%
\pgfpathcurveto{\pgfqpoint{0.021649in}{0.037276in}}{\pgfqpoint{0.011050in}{0.041667in}}{\pgfqpoint{0.000000in}{0.041667in}}%
\pgfpathcurveto{\pgfqpoint{-0.011050in}{0.041667in}}{\pgfqpoint{-0.021649in}{0.037276in}}{\pgfqpoint{-0.029463in}{0.029463in}}%
\pgfpathcurveto{\pgfqpoint{-0.037276in}{0.021649in}}{\pgfqpoint{-0.041667in}{0.011050in}}{\pgfqpoint{-0.041667in}{0.000000in}}%
\pgfpathcurveto{\pgfqpoint{-0.041667in}{-0.011050in}}{\pgfqpoint{-0.037276in}{-0.021649in}}{\pgfqpoint{-0.029463in}{-0.029463in}}%
\pgfpathcurveto{\pgfqpoint{-0.021649in}{-0.037276in}}{\pgfqpoint{-0.011050in}{-0.041667in}}{\pgfqpoint{0.000000in}{-0.041667in}}%
\pgfpathlineto{\pgfqpoint{0.000000in}{-0.041667in}}%
\pgfpathclose%
\pgfusepath{stroke,fill}%
}%
\begin{pgfscope}%
\pgfsys@transformshift{0.720247in}{2.948023in}%
\pgfsys@useobject{currentmarker}{}%
\end{pgfscope}%
\begin{pgfscope}%
\pgfsys@transformshift{0.944820in}{2.810353in}%
\pgfsys@useobject{currentmarker}{}%
\end{pgfscope}%
\begin{pgfscope}%
\pgfsys@transformshift{1.169394in}{2.840315in}%
\pgfsys@useobject{currentmarker}{}%
\end{pgfscope}%
\begin{pgfscope}%
\pgfsys@transformshift{1.393968in}{2.822186in}%
\pgfsys@useobject{currentmarker}{}%
\end{pgfscope}%
\begin{pgfscope}%
\pgfsys@transformshift{1.618542in}{2.777739in}%
\pgfsys@useobject{currentmarker}{}%
\end{pgfscope}%
\begin{pgfscope}%
\pgfsys@transformshift{1.843116in}{2.771821in}%
\pgfsys@useobject{currentmarker}{}%
\end{pgfscope}%
\begin{pgfscope}%
\pgfsys@transformshift{2.067690in}{2.766750in}%
\pgfsys@useobject{currentmarker}{}%
\end{pgfscope}%
\begin{pgfscope}%
\pgfsys@transformshift{2.292264in}{2.749723in}%
\pgfsys@useobject{currentmarker}{}%
\end{pgfscope}%
\begin{pgfscope}%
\pgfsys@transformshift{2.516838in}{2.751560in}%
\pgfsys@useobject{currentmarker}{}%
\end{pgfscope}%
\begin{pgfscope}%
\pgfsys@transformshift{2.741411in}{2.612445in}%
\pgfsys@useobject{currentmarker}{}%
\end{pgfscope}%
\begin{pgfscope}%
\pgfsys@transformshift{2.965985in}{2.655978in}%
\pgfsys@useobject{currentmarker}{}%
\end{pgfscope}%
\begin{pgfscope}%
\pgfsys@transformshift{3.190559in}{2.642148in}%
\pgfsys@useobject{currentmarker}{}%
\end{pgfscope}%
\begin{pgfscope}%
\pgfsys@transformshift{3.415133in}{2.615803in}%
\pgfsys@useobject{currentmarker}{}%
\end{pgfscope}%
\begin{pgfscope}%
\pgfsys@transformshift{3.639707in}{2.620383in}%
\pgfsys@useobject{currentmarker}{}%
\end{pgfscope}%
\begin{pgfscope}%
\pgfsys@transformshift{3.864281in}{2.619014in}%
\pgfsys@useobject{currentmarker}{}%
\end{pgfscope}%
\begin{pgfscope}%
\pgfsys@transformshift{4.088855in}{2.619081in}%
\pgfsys@useobject{currentmarker}{}%
\end{pgfscope}%
\end{pgfscope}%
\begin{pgfscope}%
\pgfpathrectangle{\pgfqpoint{0.540587in}{0.499691in}}{\pgfqpoint{3.952500in}{2.552550in}}%
\pgfusepath{clip}%
\pgfsetbuttcap%
\pgfsetroundjoin%
\pgfsetlinewidth{1.104125pt}%
\definecolor{currentstroke}{rgb}{0.000000,0.000000,0.000000}%
\pgfsetstrokecolor{currentstroke}%
\pgfsetdash{{7.040000pt}{1.760000pt}{1.100000pt}{1.760000pt}}{0.000000pt}%
\pgfpathmoveto{\pgfqpoint{0.720247in}{2.961609in}}%
\pgfpathlineto{\pgfqpoint{0.944820in}{2.959516in}}%
\pgfpathlineto{\pgfqpoint{1.169394in}{2.954672in}}%
\pgfpathlineto{\pgfqpoint{1.393968in}{2.944159in}}%
\pgfpathlineto{\pgfqpoint{1.618542in}{2.936382in}}%
\pgfpathlineto{\pgfqpoint{1.843116in}{2.917178in}}%
\pgfpathlineto{\pgfqpoint{2.067690in}{2.906862in}}%
\pgfpathlineto{\pgfqpoint{2.292264in}{2.872943in}}%
\pgfpathlineto{\pgfqpoint{2.516838in}{2.856808in}}%
\pgfpathlineto{\pgfqpoint{2.741411in}{2.840927in}}%
\pgfpathlineto{\pgfqpoint{2.965985in}{2.789632in}}%
\pgfpathlineto{\pgfqpoint{3.190559in}{2.732968in}}%
\pgfpathlineto{\pgfqpoint{3.415133in}{2.598616in}}%
\pgfpathlineto{\pgfqpoint{3.639707in}{2.409291in}}%
\pgfpathlineto{\pgfqpoint{3.864281in}{1.964532in}}%
\pgfpathlineto{\pgfqpoint{4.088855in}{1.148968in}}%
\pgfpathlineto{\pgfqpoint{4.089816in}{0.496358in}}%
\pgfusepath{stroke}%
\end{pgfscope}%
\begin{pgfscope}%
\pgfpathrectangle{\pgfqpoint{0.540587in}{0.499691in}}{\pgfqpoint{3.952500in}{2.552550in}}%
\pgfusepath{clip}%
\pgfsetbuttcap%
\pgfsetroundjoin%
\pgfsetlinewidth{1.104125pt}%
\definecolor{currentstroke}{rgb}{0.000000,0.000000,0.000000}%
\pgfsetstrokecolor{currentstroke}%
\pgfsetdash{{1.100000pt}{1.815000pt}}{0.000000pt}%
\pgfpathmoveto{\pgfqpoint{0.944820in}{2.775804in}}%
\pgfpathlineto{\pgfqpoint{1.169394in}{2.822267in}}%
\pgfpathlineto{\pgfqpoint{1.393968in}{2.856942in}}%
\pgfpathlineto{\pgfqpoint{1.618542in}{2.819397in}}%
\pgfpathlineto{\pgfqpoint{1.843116in}{2.868707in}}%
\pgfpathlineto{\pgfqpoint{2.067690in}{2.810060in}}%
\pgfpathlineto{\pgfqpoint{2.292264in}{2.868707in}}%
\pgfpathlineto{\pgfqpoint{2.516838in}{2.793052in}}%
\pgfpathlineto{\pgfqpoint{2.741411in}{2.775804in}}%
\pgfpathlineto{\pgfqpoint{2.965985in}{2.822255in}}%
\pgfpathlineto{\pgfqpoint{3.190559in}{2.775267in}}%
\pgfpathlineto{\pgfqpoint{3.415133in}{2.746501in}}%
\pgfpathlineto{\pgfqpoint{3.639707in}{2.617745in}}%
\pgfpathlineto{\pgfqpoint{3.864281in}{2.432428in}}%
\pgfpathlineto{\pgfqpoint{4.088855in}{1.987757in}}%
\pgfpathlineto{\pgfqpoint{4.313428in}{2.335557in}}%
\pgfusepath{stroke}%
\end{pgfscope}%
\begin{pgfscope}%
\pgfpathrectangle{\pgfqpoint{0.540587in}{0.499691in}}{\pgfqpoint{3.952500in}{2.552550in}}%
\pgfusepath{clip}%
\pgfsetbuttcap%
\pgfsetroundjoin%
\pgfsetlinewidth{1.104125pt}%
\definecolor{currentstroke}{rgb}{0.000000,0.000000,0.000000}%
\pgfsetstrokecolor{currentstroke}%
\pgfsetdash{{4.070000pt}{1.760000pt}}{0.000000pt}%
\pgfpathmoveto{\pgfqpoint{0.720247in}{1.753876in}}%
\pgfpathlineto{\pgfqpoint{0.944820in}{1.753876in}}%
\pgfpathlineto{\pgfqpoint{1.169394in}{1.753876in}}%
\pgfpathlineto{\pgfqpoint{1.393968in}{1.753876in}}%
\pgfpathlineto{\pgfqpoint{1.618542in}{1.753876in}}%
\pgfpathlineto{\pgfqpoint{1.843116in}{1.753876in}}%
\pgfpathlineto{\pgfqpoint{2.067690in}{1.753876in}}%
\pgfpathlineto{\pgfqpoint{2.292264in}{1.753876in}}%
\pgfpathlineto{\pgfqpoint{2.516838in}{1.753876in}}%
\pgfpathlineto{\pgfqpoint{2.741411in}{1.753876in}}%
\pgfpathlineto{\pgfqpoint{2.965985in}{1.753876in}}%
\pgfpathlineto{\pgfqpoint{3.190559in}{1.753876in}}%
\pgfpathlineto{\pgfqpoint{3.415133in}{1.753876in}}%
\pgfpathlineto{\pgfqpoint{3.639707in}{1.753876in}}%
\pgfpathlineto{\pgfqpoint{3.864281in}{1.753876in}}%
\pgfpathlineto{\pgfqpoint{4.088855in}{1.753876in}}%
\pgfpathlineto{\pgfqpoint{4.313428in}{1.753876in}}%
\pgfusepath{stroke}%
\end{pgfscope}%
\begin{pgfscope}%
\pgfsetrectcap%
\pgfsetmiterjoin%
\pgfsetlinewidth{0.803000pt}%
\definecolor{currentstroke}{rgb}{0.000000,0.000000,0.000000}%
\pgfsetstrokecolor{currentstroke}%
\pgfsetdash{}{0pt}%
\pgfpathmoveto{\pgfqpoint{0.540587in}{0.499691in}}%
\pgfpathlineto{\pgfqpoint{0.540587in}{3.052241in}}%
\pgfusepath{stroke}%
\end{pgfscope}%
\begin{pgfscope}%
\pgfsetrectcap%
\pgfsetmiterjoin%
\pgfsetlinewidth{0.803000pt}%
\definecolor{currentstroke}{rgb}{0.000000,0.000000,0.000000}%
\pgfsetstrokecolor{currentstroke}%
\pgfsetdash{}{0pt}%
\pgfpathmoveto{\pgfqpoint{4.493088in}{0.499691in}}%
\pgfpathlineto{\pgfqpoint{4.493088in}{3.052241in}}%
\pgfusepath{stroke}%
\end{pgfscope}%
\begin{pgfscope}%
\pgfsetrectcap%
\pgfsetmiterjoin%
\pgfsetlinewidth{0.803000pt}%
\definecolor{currentstroke}{rgb}{0.000000,0.000000,0.000000}%
\pgfsetstrokecolor{currentstroke}%
\pgfsetdash{}{0pt}%
\pgfpathmoveto{\pgfqpoint{0.540587in}{0.499691in}}%
\pgfpathlineto{\pgfqpoint{4.493088in}{0.499691in}}%
\pgfusepath{stroke}%
\end{pgfscope}%
\begin{pgfscope}%
\pgfsetrectcap%
\pgfsetmiterjoin%
\pgfsetlinewidth{0.803000pt}%
\definecolor{currentstroke}{rgb}{0.000000,0.000000,0.000000}%
\pgfsetstrokecolor{currentstroke}%
\pgfsetdash{}{0pt}%
\pgfpathmoveto{\pgfqpoint{0.540587in}{3.052241in}}%
\pgfpathlineto{\pgfqpoint{4.493088in}{3.052241in}}%
\pgfusepath{stroke}%
\end{pgfscope}%
\begin{pgfscope}%
\pgfsetbuttcap%
\pgfsetmiterjoin%
\definecolor{currentfill}{rgb}{1.000000,1.000000,1.000000}%
\pgfsetfillcolor{currentfill}%
\pgfsetfillopacity{0.800000}%
\pgfsetlinewidth{1.003750pt}%
\definecolor{currentstroke}{rgb}{0.800000,0.800000,0.800000}%
\pgfsetstrokecolor{currentstroke}%
\pgfsetstrokeopacity{0.800000}%
\pgfsetdash{}{0pt}%
\pgfpathmoveto{\pgfqpoint{0.637810in}{0.569136in}}%
\pgfpathlineto{\pgfqpoint{3.390468in}{0.569136in}}%
\pgfpathquadraticcurveto{\pgfqpoint{3.418245in}{0.569136in}}{\pgfqpoint{3.418245in}{0.596913in}}%
\pgfpathlineto{\pgfqpoint{3.418245in}{1.661914in}}%
\pgfpathquadraticcurveto{\pgfqpoint{3.418245in}{1.689691in}}{\pgfqpoint{3.390468in}{1.689691in}}%
\pgfpathlineto{\pgfqpoint{0.637810in}{1.689691in}}%
\pgfpathquadraticcurveto{\pgfqpoint{0.610032in}{1.689691in}}{\pgfqpoint{0.610032in}{1.661914in}}%
\pgfpathlineto{\pgfqpoint{0.610032in}{0.596913in}}%
\pgfpathquadraticcurveto{\pgfqpoint{0.610032in}{0.569136in}}{\pgfqpoint{0.637810in}{0.569136in}}%
\pgfpathlineto{\pgfqpoint{0.637810in}{0.569136in}}%
\pgfpathclose%
\pgfusepath{stroke,fill}%
\end{pgfscope}%
\begin{pgfscope}%
\pgfsetrectcap%
\pgfsetroundjoin%
\pgfsetlinewidth{1.104125pt}%
\definecolor{currentstroke}{rgb}{0.000000,0.000000,0.000000}%
\pgfsetstrokecolor{currentstroke}%
\pgfsetdash{}{0pt}%
\pgfpathmoveto{\pgfqpoint{0.665587in}{1.562747in}}%
\pgfpathlineto{\pgfqpoint{0.804476in}{1.562747in}}%
\pgfpathlineto{\pgfqpoint{0.943365in}{1.562747in}}%
\pgfusepath{stroke}%
\end{pgfscope}%
\begin{pgfscope}%
\pgfsetbuttcap%
\pgfsetroundjoin%
\definecolor{currentfill}{rgb}{0.000000,0.000000,0.000000}%
\pgfsetfillcolor{currentfill}%
\pgfsetfillopacity{0.000000}%
\pgfsetlinewidth{1.003750pt}%
\definecolor{currentstroke}{rgb}{0.000000,0.000000,0.000000}%
\pgfsetstrokecolor{currentstroke}%
\pgfsetdash{}{0pt}%
\pgfsys@defobject{currentmarker}{\pgfqpoint{-0.041667in}{-0.041667in}}{\pgfqpoint{0.041667in}{0.041667in}}{%
\pgfpathmoveto{\pgfqpoint{-0.041667in}{-0.041667in}}%
\pgfpathlineto{\pgfqpoint{0.041667in}{0.041667in}}%
\pgfpathmoveto{\pgfqpoint{-0.041667in}{0.041667in}}%
\pgfpathlineto{\pgfqpoint{0.041667in}{-0.041667in}}%
\pgfusepath{stroke,fill}%
}%
\begin{pgfscope}%
\pgfsys@transformshift{0.804476in}{1.562747in}%
\pgfsys@useobject{currentmarker}{}%
\end{pgfscope}%
\end{pgfscope}%
\begin{pgfscope}%
\definecolor{textcolor}{rgb}{0.000000,0.000000,0.000000}%
\pgfsetstrokecolor{textcolor}%
\pgfsetfillcolor{textcolor}%
\pgftext[x=1.054476in,y=1.514135in,left,base]{\color{textcolor}\rmfamily\fontsize{10.000000}{12.000000}\selectfont \(\displaystyle \norm{\tns{C}_k - D^3 f(\vek{x}_k)}_F / \norm{D^3 f(\vek{x}_k)}_F\)}%
\end{pgfscope}%
\begin{pgfscope}%
\pgfsetrectcap%
\pgfsetroundjoin%
\pgfsetlinewidth{1.104125pt}%
\definecolor{currentstroke}{rgb}{0.000000,0.000000,0.000000}%
\pgfsetstrokecolor{currentstroke}%
\pgfsetdash{}{0pt}%
\pgfpathmoveto{\pgfqpoint{0.665587in}{1.338579in}}%
\pgfpathlineto{\pgfqpoint{0.804476in}{1.338579in}}%
\pgfpathlineto{\pgfqpoint{0.943365in}{1.338579in}}%
\pgfusepath{stroke}%
\end{pgfscope}%
\begin{pgfscope}%
\pgfsetbuttcap%
\pgfsetroundjoin%
\definecolor{currentfill}{rgb}{0.000000,0.000000,0.000000}%
\pgfsetfillcolor{currentfill}%
\pgfsetfillopacity{0.000000}%
\pgfsetlinewidth{1.003750pt}%
\definecolor{currentstroke}{rgb}{0.000000,0.000000,0.000000}%
\pgfsetstrokecolor{currentstroke}%
\pgfsetdash{}{0pt}%
\pgfsys@defobject{currentmarker}{\pgfqpoint{-0.041667in}{-0.041667in}}{\pgfqpoint{0.041667in}{0.041667in}}{%
\pgfpathmoveto{\pgfqpoint{0.000000in}{-0.041667in}}%
\pgfpathcurveto{\pgfqpoint{0.011050in}{-0.041667in}}{\pgfqpoint{0.021649in}{-0.037276in}}{\pgfqpoint{0.029463in}{-0.029463in}}%
\pgfpathcurveto{\pgfqpoint{0.037276in}{-0.021649in}}{\pgfqpoint{0.041667in}{-0.011050in}}{\pgfqpoint{0.041667in}{0.000000in}}%
\pgfpathcurveto{\pgfqpoint{0.041667in}{0.011050in}}{\pgfqpoint{0.037276in}{0.021649in}}{\pgfqpoint{0.029463in}{0.029463in}}%
\pgfpathcurveto{\pgfqpoint{0.021649in}{0.037276in}}{\pgfqpoint{0.011050in}{0.041667in}}{\pgfqpoint{0.000000in}{0.041667in}}%
\pgfpathcurveto{\pgfqpoint{-0.011050in}{0.041667in}}{\pgfqpoint{-0.021649in}{0.037276in}}{\pgfqpoint{-0.029463in}{0.029463in}}%
\pgfpathcurveto{\pgfqpoint{-0.037276in}{0.021649in}}{\pgfqpoint{-0.041667in}{0.011050in}}{\pgfqpoint{-0.041667in}{0.000000in}}%
\pgfpathcurveto{\pgfqpoint{-0.041667in}{-0.011050in}}{\pgfqpoint{-0.037276in}{-0.021649in}}{\pgfqpoint{-0.029463in}{-0.029463in}}%
\pgfpathcurveto{\pgfqpoint{-0.021649in}{-0.037276in}}{\pgfqpoint{-0.011050in}{-0.041667in}}{\pgfqpoint{0.000000in}{-0.041667in}}%
\pgfpathlineto{\pgfqpoint{0.000000in}{-0.041667in}}%
\pgfpathclose%
\pgfusepath{stroke,fill}%
}%
\begin{pgfscope}%
\pgfsys@transformshift{0.804476in}{1.338579in}%
\pgfsys@useobject{currentmarker}{}%
\end{pgfscope}%
\end{pgfscope}%
\begin{pgfscope}%
\definecolor{textcolor}{rgb}{0.000000,0.000000,0.000000}%
\pgfsetstrokecolor{textcolor}%
\pgfsetfillcolor{textcolor}%
\pgftext[x=1.054476in,y=1.289968in,left,base]{\color{textcolor}\rmfamily\fontsize{10.000000}{12.000000}\selectfont \(\displaystyle \norm{(\tns{C}_k - D^3 f(\vek{x}_k))\brack{\vek{s}_k^{\rightarrow}}}_F / \norm{D^3 f(\vek{x}_k)}_F\)}%
\end{pgfscope}%
\begin{pgfscope}%
\pgfsetbuttcap%
\pgfsetroundjoin%
\pgfsetlinewidth{1.104125pt}%
\definecolor{currentstroke}{rgb}{0.000000,0.000000,0.000000}%
\pgfsetstrokecolor{currentstroke}%
\pgfsetdash{{7.040000pt}{1.760000pt}{1.100000pt}{1.760000pt}}{0.000000pt}%
\pgfpathmoveto{\pgfqpoint{0.665587in}{1.130246in}}%
\pgfpathlineto{\pgfqpoint{0.804476in}{1.130246in}}%
\pgfpathlineto{\pgfqpoint{0.943365in}{1.130246in}}%
\pgfusepath{stroke}%
\end{pgfscope}%
\begin{pgfscope}%
\definecolor{textcolor}{rgb}{0.000000,0.000000,0.000000}%
\pgfsetstrokecolor{textcolor}%
\pgfsetfillcolor{textcolor}%
\pgftext[x=1.054476in,y=1.081635in,left,base]{\color{textcolor}\rmfamily\fontsize{10.000000}{12.000000}\selectfont \(\displaystyle \norm{\vek{x}_k - \vek{x}_*}_2 / \norm{\vek{x}_*}_2\)}%
\end{pgfscope}%
\begin{pgfscope}%
\pgfsetbuttcap%
\pgfsetroundjoin%
\pgfsetlinewidth{1.104125pt}%
\definecolor{currentstroke}{rgb}{0.000000,0.000000,0.000000}%
\pgfsetstrokecolor{currentstroke}%
\pgfsetdash{{1.100000pt}{1.815000pt}}{0.000000pt}%
\pgfpathmoveto{\pgfqpoint{0.665587in}{0.921913in}}%
\pgfpathlineto{\pgfqpoint{0.804476in}{0.921913in}}%
\pgfpathlineto{\pgfqpoint{0.943365in}{0.921913in}}%
\pgfusepath{stroke}%
\end{pgfscope}%
\begin{pgfscope}%
\definecolor{textcolor}{rgb}{0.000000,0.000000,0.000000}%
\pgfsetstrokecolor{textcolor}%
\pgfsetfillcolor{textcolor}%
\pgftext[x=1.054476in,y=0.873302in,left,base]{\color{textcolor}\rmfamily\fontsize{10.000000}{12.000000}\selectfont \(\displaystyle \max \{\norm{\vek{s}_{k-1}}_2, \e_{\text{mach}} / \norm{\vek{s}_{k-1}}_2\}\)}%
\end{pgfscope}%
\begin{pgfscope}%
\pgfsetbuttcap%
\pgfsetroundjoin%
\pgfsetlinewidth{1.104125pt}%
\definecolor{currentstroke}{rgb}{0.000000,0.000000,0.000000}%
\pgfsetstrokecolor{currentstroke}%
\pgfsetdash{{4.070000pt}{1.760000pt}}{0.000000pt}%
\pgfpathmoveto{\pgfqpoint{0.665587in}{0.710233in}}%
\pgfpathlineto{\pgfqpoint{0.804476in}{0.710233in}}%
\pgfpathlineto{\pgfqpoint{0.943365in}{0.710233in}}%
\pgfusepath{stroke}%
\end{pgfscope}%
\begin{pgfscope}%
\definecolor{textcolor}{rgb}{0.000000,0.000000,0.000000}%
\pgfsetstrokecolor{textcolor}%
\pgfsetfillcolor{textcolor}%
\pgftext[x=1.054476in,y=0.661622in,left,base]{\color{textcolor}\rmfamily\fontsize{10.000000}{12.000000}\selectfont \(\displaystyle \sqrt{\e_{\text{mach}}}\)}%
\end{pgfscope}%
\end{pgfpicture}%
\makeatother%
\endgroup%

%% file: plots/angles_comparison.pgf
\begingroup%
\makeatletter%
\begin{pgfpicture}%
\pgfpathrectangle{\pgfpointorigin}{\pgfqpoint{4.683353in}{3.200466in}}%
\pgfusepath{use as bounding box, clip}%
\begin{pgfscope}%
\pgfsetbuttcap%
\pgfsetmiterjoin%
\definecolor{currentfill}{rgb}{1.000000,1.000000,1.000000}%
\pgfsetfillcolor{currentfill}%
\pgfsetlinewidth{0.000000pt}%
\definecolor{currentstroke}{rgb}{1.000000,1.000000,1.000000}%
\pgfsetstrokecolor{currentstroke}%
\pgfsetdash{}{0pt}%
\pgfpathmoveto{\pgfqpoint{0.000000in}{0.000000in}}%
\pgfpathlineto{\pgfqpoint{4.683353in}{0.000000in}}%
\pgfpathlineto{\pgfqpoint{4.683353in}{3.200466in}}%
\pgfpathlineto{\pgfqpoint{0.000000in}{3.200466in}}%
\pgfpathlineto{\pgfqpoint{0.000000in}{0.000000in}}%
\pgfpathclose%
\pgfusepath{fill}%
\end{pgfscope}%
\begin{pgfscope}%
\pgfsetbuttcap%
\pgfsetmiterjoin%
\definecolor{currentfill}{rgb}{1.000000,1.000000,1.000000}%
\pgfsetfillcolor{currentfill}%
\pgfsetlinewidth{0.000000pt}%
\definecolor{currentstroke}{rgb}{0.000000,0.000000,0.000000}%
\pgfsetstrokecolor{currentstroke}%
\pgfsetstrokeopacity{0.000000}%
\pgfsetdash{}{0pt}%
\pgfpathmoveto{\pgfqpoint{0.630853in}{0.499691in}}%
\pgfpathlineto{\pgfqpoint{2.427444in}{0.499691in}}%
\pgfpathlineto{\pgfqpoint{2.427444in}{3.052241in}}%
\pgfpathlineto{\pgfqpoint{0.630853in}{3.052241in}}%
\pgfpathlineto{\pgfqpoint{0.630853in}{0.499691in}}%
\pgfpathclose%
\pgfusepath{fill}%
\end{pgfscope}%
\begin{pgfscope}%
\pgfsetbuttcap%
\pgfsetroundjoin%
\definecolor{currentfill}{rgb}{0.000000,0.000000,0.000000}%
\pgfsetfillcolor{currentfill}%
\pgfsetfillopacity{0.000000}%
\pgfsetlinewidth{0.803000pt}%
\definecolor{currentstroke}{rgb}{0.000000,0.000000,0.000000}%
\pgfsetstrokecolor{currentstroke}%
\pgfsetdash{}{0pt}%
\pgfsys@defobject{currentmarker}{\pgfqpoint{0.000000in}{-0.048611in}}{\pgfqpoint{0.000000in}{0.000000in}}{%
\pgfpathmoveto{\pgfqpoint{0.000000in}{0.000000in}}%
\pgfpathlineto{\pgfqpoint{0.000000in}{-0.048611in}}%
\pgfusepath{stroke,fill}%
}%
\begin{pgfscope}%
\pgfsys@transformshift{0.712517in}{0.499691in}%
\pgfsys@useobject{currentmarker}{}%
\end{pgfscope}%
\end{pgfscope}%
\begin{pgfscope}%
\definecolor{textcolor}{rgb}{0.000000,0.000000,0.000000}%
\pgfsetstrokecolor{textcolor}%
\pgfsetfillcolor{textcolor}%
\pgftext[x=0.712517in,y=0.402469in,,top]{\color{textcolor}\rmfamily\fontsize{10.000000}{12.000000}\selectfont \(\displaystyle {0}\)}%
\end{pgfscope}%
\begin{pgfscope}%
\pgfsetbuttcap%
\pgfsetroundjoin%
\definecolor{currentfill}{rgb}{0.000000,0.000000,0.000000}%
\pgfsetfillcolor{currentfill}%
\pgfsetfillopacity{0.000000}%
\pgfsetlinewidth{0.803000pt}%
\definecolor{currentstroke}{rgb}{0.000000,0.000000,0.000000}%
\pgfsetstrokecolor{currentstroke}%
\pgfsetdash{}{0pt}%
\pgfsys@defobject{currentmarker}{\pgfqpoint{0.000000in}{-0.048611in}}{\pgfqpoint{0.000000in}{0.000000in}}{%
\pgfpathmoveto{\pgfqpoint{0.000000in}{0.000000in}}%
\pgfpathlineto{\pgfqpoint{0.000000in}{-0.048611in}}%
\pgfusepath{stroke,fill}%
}%
\begin{pgfscope}%
\pgfsys@transformshift{1.454910in}{0.499691in}%
\pgfsys@useobject{currentmarker}{}%
\end{pgfscope}%
\end{pgfscope}%
\begin{pgfscope}%
\definecolor{textcolor}{rgb}{0.000000,0.000000,0.000000}%
\pgfsetstrokecolor{textcolor}%
\pgfsetfillcolor{textcolor}%
\pgftext[x=1.454910in,y=0.402469in,,top]{\color{textcolor}\rmfamily\fontsize{10.000000}{12.000000}\selectfont \(\displaystyle {20}\)}%
\end{pgfscope}%
\begin{pgfscope}%
\pgfsetbuttcap%
\pgfsetroundjoin%
\definecolor{currentfill}{rgb}{0.000000,0.000000,0.000000}%
\pgfsetfillcolor{currentfill}%
\pgfsetfillopacity{0.000000}%
\pgfsetlinewidth{0.803000pt}%
\definecolor{currentstroke}{rgb}{0.000000,0.000000,0.000000}%
\pgfsetstrokecolor{currentstroke}%
\pgfsetdash{}{0pt}%
\pgfsys@defobject{currentmarker}{\pgfqpoint{0.000000in}{-0.048611in}}{\pgfqpoint{0.000000in}{0.000000in}}{%
\pgfpathmoveto{\pgfqpoint{0.000000in}{0.000000in}}%
\pgfpathlineto{\pgfqpoint{0.000000in}{-0.048611in}}%
\pgfusepath{stroke,fill}%
}%
\begin{pgfscope}%
\pgfsys@transformshift{2.197303in}{0.499691in}%
\pgfsys@useobject{currentmarker}{}%
\end{pgfscope}%
\end{pgfscope}%
\begin{pgfscope}%
\definecolor{textcolor}{rgb}{0.000000,0.000000,0.000000}%
\pgfsetstrokecolor{textcolor}%
\pgfsetfillcolor{textcolor}%
\pgftext[x=2.197303in,y=0.402469in,,top]{\color{textcolor}\rmfamily\fontsize{10.000000}{12.000000}\selectfont \(\displaystyle {40}\)}%
\end{pgfscope}%
\begin{pgfscope}%
\definecolor{textcolor}{rgb}{0.000000,0.000000,0.000000}%
\pgfsetstrokecolor{textcolor}%
\pgfsetfillcolor{textcolor}%
\pgftext[x=1.529149in,y=0.223457in,,top]{\color{textcolor}\rmfamily\fontsize{10.000000}{12.000000}\selectfont Iteration \(\displaystyle k\)}%
\end{pgfscope}%
\begin{pgfscope}%
\pgfsetbuttcap%
\pgfsetroundjoin%
\definecolor{currentfill}{rgb}{0.000000,0.000000,0.000000}%
\pgfsetfillcolor{currentfill}%
\pgfsetfillopacity{0.000000}%
\pgfsetlinewidth{0.803000pt}%
\definecolor{currentstroke}{rgb}{0.000000,0.000000,0.000000}%
\pgfsetstrokecolor{currentstroke}%
\pgfsetdash{}{0pt}%
\pgfsys@defobject{currentmarker}{\pgfqpoint{-0.048611in}{0.000000in}}{\pgfqpoint{-0.000000in}{0.000000in}}{%
\pgfpathmoveto{\pgfqpoint{-0.000000in}{0.000000in}}%
\pgfpathlineto{\pgfqpoint{-0.048611in}{0.000000in}}%
\pgfusepath{stroke,fill}%
}%
\begin{pgfscope}%
\pgfsys@transformshift{0.630853in}{0.499691in}%
\pgfsys@useobject{currentmarker}{}%
\end{pgfscope}%
\end{pgfscope}%
\begin{pgfscope}%
\definecolor{textcolor}{rgb}{0.000000,0.000000,0.000000}%
\pgfsetstrokecolor{textcolor}%
\pgfsetfillcolor{textcolor}%
\pgftext[x=0.417902in, y=0.451466in, left, base]{\color{textcolor}\rmfamily\fontsize{10.000000}{12.000000}\selectfont 0\textdegree{}}%
\end{pgfscope}%
\begin{pgfscope}%
\pgfsetbuttcap%
\pgfsetroundjoin%
\definecolor{currentfill}{rgb}{0.000000,0.000000,0.000000}%
\pgfsetfillcolor{currentfill}%
\pgfsetfillopacity{0.000000}%
\pgfsetlinewidth{0.803000pt}%
\definecolor{currentstroke}{rgb}{0.000000,0.000000,0.000000}%
\pgfsetstrokecolor{currentstroke}%
\pgfsetdash{}{0pt}%
\pgfsys@defobject{currentmarker}{\pgfqpoint{-0.048611in}{0.000000in}}{\pgfqpoint{-0.000000in}{0.000000in}}{%
\pgfpathmoveto{\pgfqpoint{-0.000000in}{0.000000in}}%
\pgfpathlineto{\pgfqpoint{-0.048611in}{0.000000in}}%
\pgfusepath{stroke,fill}%
}%
\begin{pgfscope}%
\pgfsys@transformshift{0.630853in}{0.783308in}%
\pgfsys@useobject{currentmarker}{}%
\end{pgfscope}%
\end{pgfscope}%
\begin{pgfscope}%
\definecolor{textcolor}{rgb}{0.000000,0.000000,0.000000}%
\pgfsetstrokecolor{textcolor}%
\pgfsetfillcolor{textcolor}%
\pgftext[x=0.348457in, y=0.735082in, left, base]{\color{textcolor}\rmfamily\fontsize{10.000000}{12.000000}\selectfont 20\textdegree{}}%
\end{pgfscope}%
\begin{pgfscope}%
\pgfsetbuttcap%
\pgfsetroundjoin%
\definecolor{currentfill}{rgb}{0.000000,0.000000,0.000000}%
\pgfsetfillcolor{currentfill}%
\pgfsetfillopacity{0.000000}%
\pgfsetlinewidth{0.803000pt}%
\definecolor{currentstroke}{rgb}{0.000000,0.000000,0.000000}%
\pgfsetstrokecolor{currentstroke}%
\pgfsetdash{}{0pt}%
\pgfsys@defobject{currentmarker}{\pgfqpoint{-0.048611in}{0.000000in}}{\pgfqpoint{-0.000000in}{0.000000in}}{%
\pgfpathmoveto{\pgfqpoint{-0.000000in}{0.000000in}}%
\pgfpathlineto{\pgfqpoint{-0.048611in}{0.000000in}}%
\pgfusepath{stroke,fill}%
}%
\begin{pgfscope}%
\pgfsys@transformshift{0.630853in}{1.066924in}%
\pgfsys@useobject{currentmarker}{}%
\end{pgfscope}%
\end{pgfscope}%
\begin{pgfscope}%
\definecolor{textcolor}{rgb}{0.000000,0.000000,0.000000}%
\pgfsetstrokecolor{textcolor}%
\pgfsetfillcolor{textcolor}%
\pgftext[x=0.348457in, y=1.018699in, left, base]{\color{textcolor}\rmfamily\fontsize{10.000000}{12.000000}\selectfont 40\textdegree{}}%
\end{pgfscope}%
\begin{pgfscope}%
\pgfsetbuttcap%
\pgfsetroundjoin%
\definecolor{currentfill}{rgb}{0.000000,0.000000,0.000000}%
\pgfsetfillcolor{currentfill}%
\pgfsetfillopacity{0.000000}%
\pgfsetlinewidth{0.803000pt}%
\definecolor{currentstroke}{rgb}{0.000000,0.000000,0.000000}%
\pgfsetstrokecolor{currentstroke}%
\pgfsetdash{}{0pt}%
\pgfsys@defobject{currentmarker}{\pgfqpoint{-0.048611in}{0.000000in}}{\pgfqpoint{-0.000000in}{0.000000in}}{%
\pgfpathmoveto{\pgfqpoint{-0.000000in}{0.000000in}}%
\pgfpathlineto{\pgfqpoint{-0.048611in}{0.000000in}}%
\pgfusepath{stroke,fill}%
}%
\begin{pgfscope}%
\pgfsys@transformshift{0.630853in}{1.350541in}%
\pgfsys@useobject{currentmarker}{}%
\end{pgfscope}%
\end{pgfscope}%
\begin{pgfscope}%
\definecolor{textcolor}{rgb}{0.000000,0.000000,0.000000}%
\pgfsetstrokecolor{textcolor}%
\pgfsetfillcolor{textcolor}%
\pgftext[x=0.348457in, y=1.302316in, left, base]{\color{textcolor}\rmfamily\fontsize{10.000000}{12.000000}\selectfont 60\textdegree{}}%
\end{pgfscope}%
\begin{pgfscope}%
\pgfsetbuttcap%
\pgfsetroundjoin%
\definecolor{currentfill}{rgb}{0.000000,0.000000,0.000000}%
\pgfsetfillcolor{currentfill}%
\pgfsetfillopacity{0.000000}%
\pgfsetlinewidth{0.803000pt}%
\definecolor{currentstroke}{rgb}{0.000000,0.000000,0.000000}%
\pgfsetstrokecolor{currentstroke}%
\pgfsetdash{}{0pt}%
\pgfsys@defobject{currentmarker}{\pgfqpoint{-0.048611in}{0.000000in}}{\pgfqpoint{-0.000000in}{0.000000in}}{%
\pgfpathmoveto{\pgfqpoint{-0.000000in}{0.000000in}}%
\pgfpathlineto{\pgfqpoint{-0.048611in}{0.000000in}}%
\pgfusepath{stroke,fill}%
}%
\begin{pgfscope}%
\pgfsys@transformshift{0.630853in}{1.634158in}%
\pgfsys@useobject{currentmarker}{}%
\end{pgfscope}%
\end{pgfscope}%
\begin{pgfscope}%
\definecolor{textcolor}{rgb}{0.000000,0.000000,0.000000}%
\pgfsetstrokecolor{textcolor}%
\pgfsetfillcolor{textcolor}%
\pgftext[x=0.348457in, y=1.585932in, left, base]{\color{textcolor}\rmfamily\fontsize{10.000000}{12.000000}\selectfont 80\textdegree{}}%
\end{pgfscope}%
\begin{pgfscope}%
\pgfsetbuttcap%
\pgfsetroundjoin%
\definecolor{currentfill}{rgb}{0.000000,0.000000,0.000000}%
\pgfsetfillcolor{currentfill}%
\pgfsetfillopacity{0.000000}%
\pgfsetlinewidth{0.803000pt}%
\definecolor{currentstroke}{rgb}{0.000000,0.000000,0.000000}%
\pgfsetstrokecolor{currentstroke}%
\pgfsetdash{}{0pt}%
\pgfsys@defobject{currentmarker}{\pgfqpoint{-0.048611in}{0.000000in}}{\pgfqpoint{-0.000000in}{0.000000in}}{%
\pgfpathmoveto{\pgfqpoint{-0.000000in}{0.000000in}}%
\pgfpathlineto{\pgfqpoint{-0.048611in}{0.000000in}}%
\pgfusepath{stroke,fill}%
}%
\begin{pgfscope}%
\pgfsys@transformshift{0.630853in}{1.917774in}%
\pgfsys@useobject{currentmarker}{}%
\end{pgfscope}%
\end{pgfscope}%
\begin{pgfscope}%
\definecolor{textcolor}{rgb}{0.000000,0.000000,0.000000}%
\pgfsetstrokecolor{textcolor}%
\pgfsetfillcolor{textcolor}%
\pgftext[x=0.279012in, y=1.869549in, left, base]{\color{textcolor}\rmfamily\fontsize{10.000000}{12.000000}\selectfont 100\textdegree{}}%
\end{pgfscope}%
\begin{pgfscope}%
\pgfsetbuttcap%
\pgfsetroundjoin%
\definecolor{currentfill}{rgb}{0.000000,0.000000,0.000000}%
\pgfsetfillcolor{currentfill}%
\pgfsetfillopacity{0.000000}%
\pgfsetlinewidth{0.803000pt}%
\definecolor{currentstroke}{rgb}{0.000000,0.000000,0.000000}%
\pgfsetstrokecolor{currentstroke}%
\pgfsetdash{}{0pt}%
\pgfsys@defobject{currentmarker}{\pgfqpoint{-0.048611in}{0.000000in}}{\pgfqpoint{-0.000000in}{0.000000in}}{%
\pgfpathmoveto{\pgfqpoint{-0.000000in}{0.000000in}}%
\pgfpathlineto{\pgfqpoint{-0.048611in}{0.000000in}}%
\pgfusepath{stroke,fill}%
}%
\begin{pgfscope}%
\pgfsys@transformshift{0.630853in}{2.201391in}%
\pgfsys@useobject{currentmarker}{}%
\end{pgfscope}%
\end{pgfscope}%
\begin{pgfscope}%
\definecolor{textcolor}{rgb}{0.000000,0.000000,0.000000}%
\pgfsetstrokecolor{textcolor}%
\pgfsetfillcolor{textcolor}%
\pgftext[x=0.279012in, y=2.153166in, left, base]{\color{textcolor}\rmfamily\fontsize{10.000000}{12.000000}\selectfont 120\textdegree{}}%
\end{pgfscope}%
\begin{pgfscope}%
\pgfsetbuttcap%
\pgfsetroundjoin%
\definecolor{currentfill}{rgb}{0.000000,0.000000,0.000000}%
\pgfsetfillcolor{currentfill}%
\pgfsetfillopacity{0.000000}%
\pgfsetlinewidth{0.803000pt}%
\definecolor{currentstroke}{rgb}{0.000000,0.000000,0.000000}%
\pgfsetstrokecolor{currentstroke}%
\pgfsetdash{}{0pt}%
\pgfsys@defobject{currentmarker}{\pgfqpoint{-0.048611in}{0.000000in}}{\pgfqpoint{-0.000000in}{0.000000in}}{%
\pgfpathmoveto{\pgfqpoint{-0.000000in}{0.000000in}}%
\pgfpathlineto{\pgfqpoint{-0.048611in}{0.000000in}}%
\pgfusepath{stroke,fill}%
}%
\begin{pgfscope}%
\pgfsys@transformshift{0.630853in}{2.485008in}%
\pgfsys@useobject{currentmarker}{}%
\end{pgfscope}%
\end{pgfscope}%
\begin{pgfscope}%
\definecolor{textcolor}{rgb}{0.000000,0.000000,0.000000}%
\pgfsetstrokecolor{textcolor}%
\pgfsetfillcolor{textcolor}%
\pgftext[x=0.279012in, y=2.436783in, left, base]{\color{textcolor}\rmfamily\fontsize{10.000000}{12.000000}\selectfont 140\textdegree{}}%
\end{pgfscope}%
\begin{pgfscope}%
\pgfsetbuttcap%
\pgfsetroundjoin%
\definecolor{currentfill}{rgb}{0.000000,0.000000,0.000000}%
\pgfsetfillcolor{currentfill}%
\pgfsetfillopacity{0.000000}%
\pgfsetlinewidth{0.803000pt}%
\definecolor{currentstroke}{rgb}{0.000000,0.000000,0.000000}%
\pgfsetstrokecolor{currentstroke}%
\pgfsetdash{}{0pt}%
\pgfsys@defobject{currentmarker}{\pgfqpoint{-0.048611in}{0.000000in}}{\pgfqpoint{-0.000000in}{0.000000in}}{%
\pgfpathmoveto{\pgfqpoint{-0.000000in}{0.000000in}}%
\pgfpathlineto{\pgfqpoint{-0.048611in}{0.000000in}}%
\pgfusepath{stroke,fill}%
}%
\begin{pgfscope}%
\pgfsys@transformshift{0.630853in}{2.768624in}%
\pgfsys@useobject{currentmarker}{}%
\end{pgfscope}%
\end{pgfscope}%
\begin{pgfscope}%
\definecolor{textcolor}{rgb}{0.000000,0.000000,0.000000}%
\pgfsetstrokecolor{textcolor}%
\pgfsetfillcolor{textcolor}%
\pgftext[x=0.279012in, y=2.720399in, left, base]{\color{textcolor}\rmfamily\fontsize{10.000000}{12.000000}\selectfont 160\textdegree{}}%
\end{pgfscope}%
\begin{pgfscope}%
\pgfsetbuttcap%
\pgfsetroundjoin%
\definecolor{currentfill}{rgb}{0.000000,0.000000,0.000000}%
\pgfsetfillcolor{currentfill}%
\pgfsetfillopacity{0.000000}%
\pgfsetlinewidth{0.803000pt}%
\definecolor{currentstroke}{rgb}{0.000000,0.000000,0.000000}%
\pgfsetstrokecolor{currentstroke}%
\pgfsetdash{}{0pt}%
\pgfsys@defobject{currentmarker}{\pgfqpoint{-0.048611in}{0.000000in}}{\pgfqpoint{-0.000000in}{0.000000in}}{%
\pgfpathmoveto{\pgfqpoint{-0.000000in}{0.000000in}}%
\pgfpathlineto{\pgfqpoint{-0.048611in}{0.000000in}}%
\pgfusepath{stroke,fill}%
}%
\begin{pgfscope}%
\pgfsys@transformshift{0.630853in}{3.052241in}%
\pgfsys@useobject{currentmarker}{}%
\end{pgfscope}%
\end{pgfscope}%
\begin{pgfscope}%
\definecolor{textcolor}{rgb}{0.000000,0.000000,0.000000}%
\pgfsetstrokecolor{textcolor}%
\pgfsetfillcolor{textcolor}%
\pgftext[x=0.279012in, y=3.004016in, left, base]{\color{textcolor}\rmfamily\fontsize{10.000000}{12.000000}\selectfont 180\textdegree{}}%
\end{pgfscope}%
\begin{pgfscope}%
\definecolor{textcolor}{rgb}{0.000000,0.000000,0.000000}%
\pgfsetstrokecolor{textcolor}%
\pgfsetfillcolor{textcolor}%
\pgftext[x=0.223457in,y=1.775966in,,bottom,rotate=90.000000]{\color{textcolor}\rmfamily\fontsize{10.000000}{12.000000}\selectfont Angle of \(\displaystyle \vek{s}_k\) with x-axis}%
\end{pgfscope}%
\begin{pgfscope}%
\pgfpathrectangle{\pgfqpoint{0.630853in}{0.499691in}}{\pgfqpoint{1.796591in}{2.552550in}}%
\pgfusepath{clip}%
\pgfsetbuttcap%
\pgfsetmiterjoin%
\definecolor{currentfill}{rgb}{0.000000,0.000000,0.000000}%
\pgfsetfillcolor{currentfill}%
\pgfsetfillopacity{0.000000}%
\pgfsetlinewidth{1.003750pt}%
\definecolor{currentstroke}{rgb}{0.000000,0.000000,0.000000}%
\pgfsetstrokecolor{currentstroke}%
\pgfsetdash{}{0pt}%
\pgfsys@defobject{currentmarker}{\pgfqpoint{-0.041667in}{-0.041667in}}{\pgfqpoint{0.041667in}{0.041667in}}{%
\pgfpathmoveto{\pgfqpoint{0.000000in}{0.041667in}}%
\pgfpathlineto{\pgfqpoint{-0.041667in}{-0.041667in}}%
\pgfpathlineto{\pgfqpoint{0.041667in}{-0.041667in}}%
\pgfpathlineto{\pgfqpoint{0.000000in}{0.041667in}}%
\pgfpathclose%
\pgfusepath{stroke,fill}%
}%
\begin{pgfscope}%
\pgfsys@transformshift{0.712517in}{0.499691in}%
\pgfsys@useobject{currentmarker}{}%
\end{pgfscope}%
\begin{pgfscope}%
\pgfsys@transformshift{0.749636in}{0.750658in}%
\pgfsys@useobject{currentmarker}{}%
\end{pgfscope}%
\begin{pgfscope}%
\pgfsys@transformshift{0.786756in}{1.414430in}%
\pgfsys@useobject{currentmarker}{}%
\end{pgfscope}%
\begin{pgfscope}%
\pgfsys@transformshift{0.823876in}{0.877993in}%
\pgfsys@useobject{currentmarker}{}%
\end{pgfscope}%
\begin{pgfscope}%
\pgfsys@transformshift{0.860995in}{1.118739in}%
\pgfsys@useobject{currentmarker}{}%
\end{pgfscope}%
\begin{pgfscope}%
\pgfsys@transformshift{0.898115in}{2.223628in}%
\pgfsys@useobject{currentmarker}{}%
\end{pgfscope}%
\begin{pgfscope}%
\pgfsys@transformshift{0.935235in}{1.178691in}%
\pgfsys@useobject{currentmarker}{}%
\end{pgfscope}%
\begin{pgfscope}%
\pgfsys@transformshift{0.972354in}{1.306920in}%
\pgfsys@useobject{currentmarker}{}%
\end{pgfscope}%
\begin{pgfscope}%
\pgfsys@transformshift{1.009474in}{1.190276in}%
\pgfsys@useobject{currentmarker}{}%
\end{pgfscope}%
\begin{pgfscope}%
\pgfsys@transformshift{1.046594in}{1.231478in}%
\pgfsys@useobject{currentmarker}{}%
\end{pgfscope}%
\begin{pgfscope}%
\pgfsys@transformshift{1.083713in}{1.327010in}%
\pgfsys@useobject{currentmarker}{}%
\end{pgfscope}%
\begin{pgfscope}%
\pgfsys@transformshift{1.120833in}{1.477701in}%
\pgfsys@useobject{currentmarker}{}%
\end{pgfscope}%
\begin{pgfscope}%
\pgfsys@transformshift{1.157952in}{1.356028in}%
\pgfsys@useobject{currentmarker}{}%
\end{pgfscope}%
\begin{pgfscope}%
\pgfsys@transformshift{1.195072in}{1.400061in}%
\pgfsys@useobject{currentmarker}{}%
\end{pgfscope}%
\begin{pgfscope}%
\pgfsys@transformshift{1.232192in}{2.682929in}%
\pgfsys@useobject{currentmarker}{}%
\end{pgfscope}%
\begin{pgfscope}%
\pgfsys@transformshift{1.269311in}{1.400928in}%
\pgfsys@useobject{currentmarker}{}%
\end{pgfscope}%
\begin{pgfscope}%
\pgfsys@transformshift{1.306431in}{0.980504in}%
\pgfsys@useobject{currentmarker}{}%
\end{pgfscope}%
\begin{pgfscope}%
\pgfsys@transformshift{1.343551in}{2.254046in}%
\pgfsys@useobject{currentmarker}{}%
\end{pgfscope}%
\begin{pgfscope}%
\pgfsys@transformshift{1.380670in}{1.399713in}%
\pgfsys@useobject{currentmarker}{}%
\end{pgfscope}%
\begin{pgfscope}%
\pgfsys@transformshift{1.417790in}{0.982888in}%
\pgfsys@useobject{currentmarker}{}%
\end{pgfscope}%
\begin{pgfscope}%
\pgfsys@transformshift{1.454910in}{1.400475in}%
\pgfsys@useobject{currentmarker}{}%
\end{pgfscope}%
\begin{pgfscope}%
\pgfsys@transformshift{1.492029in}{2.258783in}%
\pgfsys@useobject{currentmarker}{}%
\end{pgfscope}%
\begin{pgfscope}%
\pgfsys@transformshift{1.529149in}{1.399718in}%
\pgfsys@useobject{currentmarker}{}%
\end{pgfscope}%
\begin{pgfscope}%
\pgfsys@transformshift{1.566269in}{0.982804in}%
\pgfsys@useobject{currentmarker}{}%
\end{pgfscope}%
\begin{pgfscope}%
\pgfsys@transformshift{1.603388in}{1.400475in}%
\pgfsys@useobject{currentmarker}{}%
\end{pgfscope}%
\begin{pgfscope}%
\pgfsys@transformshift{1.640508in}{2.258920in}%
\pgfsys@useobject{currentmarker}{}%
\end{pgfscope}%
\begin{pgfscope}%
\pgfsys@transformshift{1.677628in}{1.399718in}%
\pgfsys@useobject{currentmarker}{}%
\end{pgfscope}%
\begin{pgfscope}%
\pgfsys@transformshift{1.714747in}{0.981256in}%
\pgfsys@useobject{currentmarker}{}%
\end{pgfscope}%
\begin{pgfscope}%
\pgfsys@transformshift{1.751867in}{1.400471in}%
\pgfsys@useobject{currentmarker}{}%
\end{pgfscope}%
\begin{pgfscope}%
\pgfsys@transformshift{1.788986in}{2.260177in}%
\pgfsys@useobject{currentmarker}{}%
\end{pgfscope}%
\begin{pgfscope}%
\pgfsys@transformshift{1.826106in}{1.399715in}%
\pgfsys@useobject{currentmarker}{}%
\end{pgfscope}%
\begin{pgfscope}%
\pgfsys@transformshift{1.863226in}{1.004613in}%
\pgfsys@useobject{currentmarker}{}%
\end{pgfscope}%
\begin{pgfscope}%
\pgfsys@transformshift{1.900345in}{1.400499in}%
\pgfsys@useobject{currentmarker}{}%
\end{pgfscope}%
\begin{pgfscope}%
\pgfsys@transformshift{1.937465in}{2.311181in}%
\pgfsys@useobject{currentmarker}{}%
\end{pgfscope}%
\begin{pgfscope}%
\pgfsys@transformshift{1.974585in}{1.399806in}%
\pgfsys@useobject{currentmarker}{}%
\end{pgfscope}%
\begin{pgfscope}%
\pgfsys@transformshift{2.011704in}{3.030864in}%
\pgfsys@useobject{currentmarker}{}%
\end{pgfscope}%
\begin{pgfscope}%
\pgfsys@transformshift{2.048824in}{1.400088in}%
\pgfsys@useobject{currentmarker}{}%
\end{pgfscope}%
\begin{pgfscope}%
\pgfsys@transformshift{2.085944in}{0.725810in}%
\pgfsys@useobject{currentmarker}{}%
\end{pgfscope}%
\begin{pgfscope}%
\pgfsys@transformshift{2.123063in}{1.394994in}%
\pgfsys@useobject{currentmarker}{}%
\end{pgfscope}%
\begin{pgfscope}%
\pgfsys@transformshift{2.160183in}{1.492824in}%
\pgfsys@useobject{currentmarker}{}%
\end{pgfscope}%
\begin{pgfscope}%
\pgfsys@transformshift{2.197303in}{1.421643in}%
\pgfsys@useobject{currentmarker}{}%
\end{pgfscope}%
\begin{pgfscope}%
\pgfsys@transformshift{2.234422in}{1.387644in}%
\pgfsys@useobject{currentmarker}{}%
\end{pgfscope}%
\begin{pgfscope}%
\pgfsys@transformshift{2.271542in}{1.422241in}%
\pgfsys@useobject{currentmarker}{}%
\end{pgfscope}%
\begin{pgfscope}%
\pgfsys@transformshift{2.308662in}{1.399252in}%
\pgfsys@useobject{currentmarker}{}%
\end{pgfscope}%
\begin{pgfscope}%
\pgfsys@transformshift{2.345781in}{0.499691in}%
\pgfsys@useobject{currentmarker}{}%
\end{pgfscope}%
\end{pgfscope}%
\begin{pgfscope}%
\pgfpathrectangle{\pgfqpoint{0.630853in}{0.499691in}}{\pgfqpoint{1.796591in}{2.552550in}}%
\pgfusepath{clip}%
\pgfsetbuttcap%
\pgfsetroundjoin%
\pgfsetlinewidth{1.104125pt}%
\definecolor{currentstroke}{rgb}{0.000000,0.000000,0.000000}%
\pgfsetstrokecolor{currentstroke}%
\pgfsetdash{{1.100000pt}{1.815000pt}}{0.000000pt}%
\pgfpathmoveto{\pgfqpoint{0.712517in}{1.399252in}}%
\pgfpathlineto{\pgfqpoint{0.749636in}{1.399252in}}%
\pgfpathlineto{\pgfqpoint{0.786756in}{1.399252in}}%
\pgfpathlineto{\pgfqpoint{0.823876in}{1.399252in}}%
\pgfpathlineto{\pgfqpoint{0.860995in}{1.399252in}}%
\pgfpathlineto{\pgfqpoint{0.898115in}{1.399252in}}%
\pgfpathlineto{\pgfqpoint{0.935235in}{1.399252in}}%
\pgfpathlineto{\pgfqpoint{0.972354in}{1.399252in}}%
\pgfpathlineto{\pgfqpoint{1.009474in}{1.399252in}}%
\pgfpathlineto{\pgfqpoint{1.046594in}{1.399252in}}%
\pgfpathlineto{\pgfqpoint{1.083713in}{1.399252in}}%
\pgfpathlineto{\pgfqpoint{1.120833in}{1.399252in}}%
\pgfpathlineto{\pgfqpoint{1.157952in}{1.399252in}}%
\pgfpathlineto{\pgfqpoint{1.195072in}{1.399252in}}%
\pgfpathlineto{\pgfqpoint{1.232192in}{1.399252in}}%
\pgfpathlineto{\pgfqpoint{1.269311in}{1.399252in}}%
\pgfpathlineto{\pgfqpoint{1.306431in}{1.399252in}}%
\pgfpathlineto{\pgfqpoint{1.343551in}{1.399252in}}%
\pgfpathlineto{\pgfqpoint{1.380670in}{1.399252in}}%
\pgfpathlineto{\pgfqpoint{1.417790in}{1.399252in}}%
\pgfpathlineto{\pgfqpoint{1.454910in}{1.399252in}}%
\pgfpathlineto{\pgfqpoint{1.492029in}{1.399252in}}%
\pgfpathlineto{\pgfqpoint{1.529149in}{1.399252in}}%
\pgfpathlineto{\pgfqpoint{1.566269in}{1.399252in}}%
\pgfpathlineto{\pgfqpoint{1.603388in}{1.399252in}}%
\pgfpathlineto{\pgfqpoint{1.640508in}{1.399252in}}%
\pgfpathlineto{\pgfqpoint{1.677628in}{1.399252in}}%
\pgfpathlineto{\pgfqpoint{1.714747in}{1.399252in}}%
\pgfpathlineto{\pgfqpoint{1.751867in}{1.399252in}}%
\pgfpathlineto{\pgfqpoint{1.788986in}{1.399252in}}%
\pgfpathlineto{\pgfqpoint{1.826106in}{1.399252in}}%
\pgfpathlineto{\pgfqpoint{1.863226in}{1.399252in}}%
\pgfpathlineto{\pgfqpoint{1.900345in}{1.399252in}}%
\pgfpathlineto{\pgfqpoint{1.937465in}{1.399252in}}%
\pgfpathlineto{\pgfqpoint{1.974585in}{1.399252in}}%
\pgfpathlineto{\pgfqpoint{2.011704in}{1.399252in}}%
\pgfpathlineto{\pgfqpoint{2.048824in}{1.399252in}}%
\pgfpathlineto{\pgfqpoint{2.085944in}{1.399252in}}%
\pgfpathlineto{\pgfqpoint{2.123063in}{1.399252in}}%
\pgfpathlineto{\pgfqpoint{2.160183in}{1.399252in}}%
\pgfpathlineto{\pgfqpoint{2.197303in}{1.399252in}}%
\pgfpathlineto{\pgfqpoint{2.234422in}{1.399252in}}%
\pgfpathlineto{\pgfqpoint{2.271542in}{1.399252in}}%
\pgfpathlineto{\pgfqpoint{2.308662in}{1.399252in}}%
\pgfpathlineto{\pgfqpoint{2.345781in}{1.399252in}}%
\pgfusepath{stroke}%
\end{pgfscope}%
\begin{pgfscope}%
\pgfsetrectcap%
\pgfsetmiterjoin%
\pgfsetlinewidth{0.803000pt}%
\definecolor{currentstroke}{rgb}{0.000000,0.000000,0.000000}%
\pgfsetstrokecolor{currentstroke}%
\pgfsetdash{}{0pt}%
\pgfpathmoveto{\pgfqpoint{0.630853in}{0.499691in}}%
\pgfpathlineto{\pgfqpoint{0.630853in}{3.052241in}}%
\pgfusepath{stroke}%
\end{pgfscope}%
\begin{pgfscope}%
\pgfsetrectcap%
\pgfsetmiterjoin%
\pgfsetlinewidth{0.803000pt}%
\definecolor{currentstroke}{rgb}{0.000000,0.000000,0.000000}%
\pgfsetstrokecolor{currentstroke}%
\pgfsetdash{}{0pt}%
\pgfpathmoveto{\pgfqpoint{2.427444in}{0.499691in}}%
\pgfpathlineto{\pgfqpoint{2.427444in}{3.052241in}}%
\pgfusepath{stroke}%
\end{pgfscope}%
\begin{pgfscope}%
\pgfsetrectcap%
\pgfsetmiterjoin%
\pgfsetlinewidth{0.803000pt}%
\definecolor{currentstroke}{rgb}{0.000000,0.000000,0.000000}%
\pgfsetstrokecolor{currentstroke}%
\pgfsetdash{}{0pt}%
\pgfpathmoveto{\pgfqpoint{0.630853in}{0.499691in}}%
\pgfpathlineto{\pgfqpoint{2.427444in}{0.499691in}}%
\pgfusepath{stroke}%
\end{pgfscope}%
\begin{pgfscope}%
\pgfsetrectcap%
\pgfsetmiterjoin%
\pgfsetlinewidth{0.803000pt}%
\definecolor{currentstroke}{rgb}{0.000000,0.000000,0.000000}%
\pgfsetstrokecolor{currentstroke}%
\pgfsetdash{}{0pt}%
\pgfpathmoveto{\pgfqpoint{0.630853in}{3.052241in}}%
\pgfpathlineto{\pgfqpoint{2.427444in}{3.052241in}}%
\pgfusepath{stroke}%
\end{pgfscope}%
\begin{pgfscope}%
\pgfsetbuttcap%
\pgfsetmiterjoin%
\definecolor{currentfill}{rgb}{1.000000,1.000000,1.000000}%
\pgfsetfillcolor{currentfill}%
\pgfsetlinewidth{0.000000pt}%
\definecolor{currentstroke}{rgb}{0.000000,0.000000,0.000000}%
\pgfsetstrokecolor{currentstroke}%
\pgfsetstrokeopacity{0.000000}%
\pgfsetdash{}{0pt}%
\pgfpathmoveto{\pgfqpoint{2.786763in}{0.499691in}}%
\pgfpathlineto{\pgfqpoint{4.583353in}{0.499691in}}%
\pgfpathlineto{\pgfqpoint{4.583353in}{3.052241in}}%
\pgfpathlineto{\pgfqpoint{2.786763in}{3.052241in}}%
\pgfpathlineto{\pgfqpoint{2.786763in}{0.499691in}}%
\pgfpathclose%
\pgfusepath{fill}%
\end{pgfscope}%
\begin{pgfscope}%
\pgfsetbuttcap%
\pgfsetroundjoin%
\definecolor{currentfill}{rgb}{0.000000,0.000000,0.000000}%
\pgfsetfillcolor{currentfill}%
\pgfsetfillopacity{0.000000}%
\pgfsetlinewidth{0.803000pt}%
\definecolor{currentstroke}{rgb}{0.000000,0.000000,0.000000}%
\pgfsetstrokecolor{currentstroke}%
\pgfsetdash{}{0pt}%
\pgfsys@defobject{currentmarker}{\pgfqpoint{0.000000in}{-0.048611in}}{\pgfqpoint{0.000000in}{0.000000in}}{%
\pgfpathmoveto{\pgfqpoint{0.000000in}{0.000000in}}%
\pgfpathlineto{\pgfqpoint{0.000000in}{-0.048611in}}%
\pgfusepath{stroke,fill}%
}%
\begin{pgfscope}%
\pgfsys@transformshift{2.868426in}{0.499691in}%
\pgfsys@useobject{currentmarker}{}%
\end{pgfscope}%
\end{pgfscope}%
\begin{pgfscope}%
\definecolor{textcolor}{rgb}{0.000000,0.000000,0.000000}%
\pgfsetstrokecolor{textcolor}%
\pgfsetfillcolor{textcolor}%
\pgftext[x=2.868426in,y=0.402469in,,top]{\color{textcolor}\rmfamily\fontsize{10.000000}{12.000000}\selectfont \(\displaystyle {0}\)}%
\end{pgfscope}%
\begin{pgfscope}%
\pgfsetbuttcap%
\pgfsetroundjoin%
\definecolor{currentfill}{rgb}{0.000000,0.000000,0.000000}%
\pgfsetfillcolor{currentfill}%
\pgfsetfillopacity{0.000000}%
\pgfsetlinewidth{0.803000pt}%
\definecolor{currentstroke}{rgb}{0.000000,0.000000,0.000000}%
\pgfsetstrokecolor{currentstroke}%
\pgfsetdash{}{0pt}%
\pgfsys@defobject{currentmarker}{\pgfqpoint{0.000000in}{-0.048611in}}{\pgfqpoint{0.000000in}{0.000000in}}{%
\pgfpathmoveto{\pgfqpoint{0.000000in}{0.000000in}}%
\pgfpathlineto{\pgfqpoint{0.000000in}{-0.048611in}}%
\pgfusepath{stroke,fill}%
}%
\begin{pgfscope}%
\pgfsys@transformshift{3.412847in}{0.499691in}%
\pgfsys@useobject{currentmarker}{}%
\end{pgfscope}%
\end{pgfscope}%
\begin{pgfscope}%
\definecolor{textcolor}{rgb}{0.000000,0.000000,0.000000}%
\pgfsetstrokecolor{textcolor}%
\pgfsetfillcolor{textcolor}%
\pgftext[x=3.412847in,y=0.402469in,,top]{\color{textcolor}\rmfamily\fontsize{10.000000}{12.000000}\selectfont \(\displaystyle {5}\)}%
\end{pgfscope}%
\begin{pgfscope}%
\pgfsetbuttcap%
\pgfsetroundjoin%
\definecolor{currentfill}{rgb}{0.000000,0.000000,0.000000}%
\pgfsetfillcolor{currentfill}%
\pgfsetfillopacity{0.000000}%
\pgfsetlinewidth{0.803000pt}%
\definecolor{currentstroke}{rgb}{0.000000,0.000000,0.000000}%
\pgfsetstrokecolor{currentstroke}%
\pgfsetdash{}{0pt}%
\pgfsys@defobject{currentmarker}{\pgfqpoint{0.000000in}{-0.048611in}}{\pgfqpoint{0.000000in}{0.000000in}}{%
\pgfpathmoveto{\pgfqpoint{0.000000in}{0.000000in}}%
\pgfpathlineto{\pgfqpoint{0.000000in}{-0.048611in}}%
\pgfusepath{stroke,fill}%
}%
\begin{pgfscope}%
\pgfsys@transformshift{3.957269in}{0.499691in}%
\pgfsys@useobject{currentmarker}{}%
\end{pgfscope}%
\end{pgfscope}%
\begin{pgfscope}%
\definecolor{textcolor}{rgb}{0.000000,0.000000,0.000000}%
\pgfsetstrokecolor{textcolor}%
\pgfsetfillcolor{textcolor}%
\pgftext[x=3.957269in,y=0.402469in,,top]{\color{textcolor}\rmfamily\fontsize{10.000000}{12.000000}\selectfont \(\displaystyle {10}\)}%
\end{pgfscope}%
\begin{pgfscope}%
\pgfsetbuttcap%
\pgfsetroundjoin%
\definecolor{currentfill}{rgb}{0.000000,0.000000,0.000000}%
\pgfsetfillcolor{currentfill}%
\pgfsetfillopacity{0.000000}%
\pgfsetlinewidth{0.803000pt}%
\definecolor{currentstroke}{rgb}{0.000000,0.000000,0.000000}%
\pgfsetstrokecolor{currentstroke}%
\pgfsetdash{}{0pt}%
\pgfsys@defobject{currentmarker}{\pgfqpoint{0.000000in}{-0.048611in}}{\pgfqpoint{0.000000in}{0.000000in}}{%
\pgfpathmoveto{\pgfqpoint{0.000000in}{0.000000in}}%
\pgfpathlineto{\pgfqpoint{0.000000in}{-0.048611in}}%
\pgfusepath{stroke,fill}%
}%
\begin{pgfscope}%
\pgfsys@transformshift{4.501690in}{0.499691in}%
\pgfsys@useobject{currentmarker}{}%
\end{pgfscope}%
\end{pgfscope}%
\begin{pgfscope}%
\definecolor{textcolor}{rgb}{0.000000,0.000000,0.000000}%
\pgfsetstrokecolor{textcolor}%
\pgfsetfillcolor{textcolor}%
\pgftext[x=4.501690in,y=0.402469in,,top]{\color{textcolor}\rmfamily\fontsize{10.000000}{12.000000}\selectfont \(\displaystyle {15}\)}%
\end{pgfscope}%
\begin{pgfscope}%
\definecolor{textcolor}{rgb}{0.000000,0.000000,0.000000}%
\pgfsetstrokecolor{textcolor}%
\pgfsetfillcolor{textcolor}%
\pgftext[x=3.685058in,y=0.223457in,,top]{\color{textcolor}\rmfamily\fontsize{10.000000}{12.000000}\selectfont Iteration \(\displaystyle k\)}%
\end{pgfscope}%
\begin{pgfscope}%
\pgfsetbuttcap%
\pgfsetroundjoin%
\definecolor{currentfill}{rgb}{0.000000,0.000000,0.000000}%
\pgfsetfillcolor{currentfill}%
\pgfsetfillopacity{0.000000}%
\pgfsetlinewidth{0.803000pt}%
\definecolor{currentstroke}{rgb}{0.000000,0.000000,0.000000}%
\pgfsetstrokecolor{currentstroke}%
\pgfsetdash{}{0pt}%
\pgfsys@defobject{currentmarker}{\pgfqpoint{-0.048611in}{0.000000in}}{\pgfqpoint{-0.000000in}{0.000000in}}{%
\pgfpathmoveto{\pgfqpoint{-0.000000in}{0.000000in}}%
\pgfpathlineto{\pgfqpoint{-0.048611in}{0.000000in}}%
\pgfusepath{stroke,fill}%
}%
\begin{pgfscope}%
\pgfsys@transformshift{2.786763in}{0.499691in}%
\pgfsys@useobject{currentmarker}{}%
\end{pgfscope}%
\end{pgfscope}%
\begin{pgfscope}%
\pgfsetbuttcap%
\pgfsetroundjoin%
\definecolor{currentfill}{rgb}{0.000000,0.000000,0.000000}%
\pgfsetfillcolor{currentfill}%
\pgfsetfillopacity{0.000000}%
\pgfsetlinewidth{0.803000pt}%
\definecolor{currentstroke}{rgb}{0.000000,0.000000,0.000000}%
\pgfsetstrokecolor{currentstroke}%
\pgfsetdash{}{0pt}%
\pgfsys@defobject{currentmarker}{\pgfqpoint{-0.048611in}{0.000000in}}{\pgfqpoint{-0.000000in}{0.000000in}}{%
\pgfpathmoveto{\pgfqpoint{-0.000000in}{0.000000in}}%
\pgfpathlineto{\pgfqpoint{-0.048611in}{0.000000in}}%
\pgfusepath{stroke,fill}%
}%
\begin{pgfscope}%
\pgfsys@transformshift{2.786763in}{0.783308in}%
\pgfsys@useobject{currentmarker}{}%
\end{pgfscope}%
\end{pgfscope}%
\begin{pgfscope}%
\pgfsetbuttcap%
\pgfsetroundjoin%
\definecolor{currentfill}{rgb}{0.000000,0.000000,0.000000}%
\pgfsetfillcolor{currentfill}%
\pgfsetfillopacity{0.000000}%
\pgfsetlinewidth{0.803000pt}%
\definecolor{currentstroke}{rgb}{0.000000,0.000000,0.000000}%
\pgfsetstrokecolor{currentstroke}%
\pgfsetdash{}{0pt}%
\pgfsys@defobject{currentmarker}{\pgfqpoint{-0.048611in}{0.000000in}}{\pgfqpoint{-0.000000in}{0.000000in}}{%
\pgfpathmoveto{\pgfqpoint{-0.000000in}{0.000000in}}%
\pgfpathlineto{\pgfqpoint{-0.048611in}{0.000000in}}%
\pgfusepath{stroke,fill}%
}%
\begin{pgfscope}%
\pgfsys@transformshift{2.786763in}{1.066924in}%
\pgfsys@useobject{currentmarker}{}%
\end{pgfscope}%
\end{pgfscope}%
\begin{pgfscope}%
\pgfsetbuttcap%
\pgfsetroundjoin%
\definecolor{currentfill}{rgb}{0.000000,0.000000,0.000000}%
\pgfsetfillcolor{currentfill}%
\pgfsetfillopacity{0.000000}%
\pgfsetlinewidth{0.803000pt}%
\definecolor{currentstroke}{rgb}{0.000000,0.000000,0.000000}%
\pgfsetstrokecolor{currentstroke}%
\pgfsetdash{}{0pt}%
\pgfsys@defobject{currentmarker}{\pgfqpoint{-0.048611in}{0.000000in}}{\pgfqpoint{-0.000000in}{0.000000in}}{%
\pgfpathmoveto{\pgfqpoint{-0.000000in}{0.000000in}}%
\pgfpathlineto{\pgfqpoint{-0.048611in}{0.000000in}}%
\pgfusepath{stroke,fill}%
}%
\begin{pgfscope}%
\pgfsys@transformshift{2.786763in}{1.350541in}%
\pgfsys@useobject{currentmarker}{}%
\end{pgfscope}%
\end{pgfscope}%
\begin{pgfscope}%
\pgfsetbuttcap%
\pgfsetroundjoin%
\definecolor{currentfill}{rgb}{0.000000,0.000000,0.000000}%
\pgfsetfillcolor{currentfill}%
\pgfsetfillopacity{0.000000}%
\pgfsetlinewidth{0.803000pt}%
\definecolor{currentstroke}{rgb}{0.000000,0.000000,0.000000}%
\pgfsetstrokecolor{currentstroke}%
\pgfsetdash{}{0pt}%
\pgfsys@defobject{currentmarker}{\pgfqpoint{-0.048611in}{0.000000in}}{\pgfqpoint{-0.000000in}{0.000000in}}{%
\pgfpathmoveto{\pgfqpoint{-0.000000in}{0.000000in}}%
\pgfpathlineto{\pgfqpoint{-0.048611in}{0.000000in}}%
\pgfusepath{stroke,fill}%
}%
\begin{pgfscope}%
\pgfsys@transformshift{2.786763in}{1.634158in}%
\pgfsys@useobject{currentmarker}{}%
\end{pgfscope}%
\end{pgfscope}%
\begin{pgfscope}%
\pgfsetbuttcap%
\pgfsetroundjoin%
\definecolor{currentfill}{rgb}{0.000000,0.000000,0.000000}%
\pgfsetfillcolor{currentfill}%
\pgfsetfillopacity{0.000000}%
\pgfsetlinewidth{0.803000pt}%
\definecolor{currentstroke}{rgb}{0.000000,0.000000,0.000000}%
\pgfsetstrokecolor{currentstroke}%
\pgfsetdash{}{0pt}%
\pgfsys@defobject{currentmarker}{\pgfqpoint{-0.048611in}{0.000000in}}{\pgfqpoint{-0.000000in}{0.000000in}}{%
\pgfpathmoveto{\pgfqpoint{-0.000000in}{0.000000in}}%
\pgfpathlineto{\pgfqpoint{-0.048611in}{0.000000in}}%
\pgfusepath{stroke,fill}%
}%
\begin{pgfscope}%
\pgfsys@transformshift{2.786763in}{1.917774in}%
\pgfsys@useobject{currentmarker}{}%
\end{pgfscope}%
\end{pgfscope}%
\begin{pgfscope}%
\pgfsetbuttcap%
\pgfsetroundjoin%
\definecolor{currentfill}{rgb}{0.000000,0.000000,0.000000}%
\pgfsetfillcolor{currentfill}%
\pgfsetfillopacity{0.000000}%
\pgfsetlinewidth{0.803000pt}%
\definecolor{currentstroke}{rgb}{0.000000,0.000000,0.000000}%
\pgfsetstrokecolor{currentstroke}%
\pgfsetdash{}{0pt}%
\pgfsys@defobject{currentmarker}{\pgfqpoint{-0.048611in}{0.000000in}}{\pgfqpoint{-0.000000in}{0.000000in}}{%
\pgfpathmoveto{\pgfqpoint{-0.000000in}{0.000000in}}%
\pgfpathlineto{\pgfqpoint{-0.048611in}{0.000000in}}%
\pgfusepath{stroke,fill}%
}%
\begin{pgfscope}%
\pgfsys@transformshift{2.786763in}{2.201391in}%
\pgfsys@useobject{currentmarker}{}%
\end{pgfscope}%
\end{pgfscope}%
\begin{pgfscope}%
\pgfsetbuttcap%
\pgfsetroundjoin%
\definecolor{currentfill}{rgb}{0.000000,0.000000,0.000000}%
\pgfsetfillcolor{currentfill}%
\pgfsetfillopacity{0.000000}%
\pgfsetlinewidth{0.803000pt}%
\definecolor{currentstroke}{rgb}{0.000000,0.000000,0.000000}%
\pgfsetstrokecolor{currentstroke}%
\pgfsetdash{}{0pt}%
\pgfsys@defobject{currentmarker}{\pgfqpoint{-0.048611in}{0.000000in}}{\pgfqpoint{-0.000000in}{0.000000in}}{%
\pgfpathmoveto{\pgfqpoint{-0.000000in}{0.000000in}}%
\pgfpathlineto{\pgfqpoint{-0.048611in}{0.000000in}}%
\pgfusepath{stroke,fill}%
}%
\begin{pgfscope}%
\pgfsys@transformshift{2.786763in}{2.485008in}%
\pgfsys@useobject{currentmarker}{}%
\end{pgfscope}%
\end{pgfscope}%
\begin{pgfscope}%
\pgfsetbuttcap%
\pgfsetroundjoin%
\definecolor{currentfill}{rgb}{0.000000,0.000000,0.000000}%
\pgfsetfillcolor{currentfill}%
\pgfsetfillopacity{0.000000}%
\pgfsetlinewidth{0.803000pt}%
\definecolor{currentstroke}{rgb}{0.000000,0.000000,0.000000}%
\pgfsetstrokecolor{currentstroke}%
\pgfsetdash{}{0pt}%
\pgfsys@defobject{currentmarker}{\pgfqpoint{-0.048611in}{0.000000in}}{\pgfqpoint{-0.000000in}{0.000000in}}{%
\pgfpathmoveto{\pgfqpoint{-0.000000in}{0.000000in}}%
\pgfpathlineto{\pgfqpoint{-0.048611in}{0.000000in}}%
\pgfusepath{stroke,fill}%
}%
\begin{pgfscope}%
\pgfsys@transformshift{2.786763in}{2.768624in}%
\pgfsys@useobject{currentmarker}{}%
\end{pgfscope}%
\end{pgfscope}%
\begin{pgfscope}%
\pgfsetbuttcap%
\pgfsetroundjoin%
\definecolor{currentfill}{rgb}{0.000000,0.000000,0.000000}%
\pgfsetfillcolor{currentfill}%
\pgfsetfillopacity{0.000000}%
\pgfsetlinewidth{0.803000pt}%
\definecolor{currentstroke}{rgb}{0.000000,0.000000,0.000000}%
\pgfsetstrokecolor{currentstroke}%
\pgfsetdash{}{0pt}%
\pgfsys@defobject{currentmarker}{\pgfqpoint{-0.048611in}{0.000000in}}{\pgfqpoint{-0.000000in}{0.000000in}}{%
\pgfpathmoveto{\pgfqpoint{-0.000000in}{0.000000in}}%
\pgfpathlineto{\pgfqpoint{-0.048611in}{0.000000in}}%
\pgfusepath{stroke,fill}%
}%
\begin{pgfscope}%
\pgfsys@transformshift{2.786763in}{3.052241in}%
\pgfsys@useobject{currentmarker}{}%
\end{pgfscope}%
\end{pgfscope}%
\begin{pgfscope}%
\pgfpathrectangle{\pgfqpoint{2.786763in}{0.499691in}}{\pgfqpoint{1.796591in}{2.552550in}}%
\pgfusepath{clip}%
\pgfsetbuttcap%
\pgfsetmiterjoin%
\definecolor{currentfill}{rgb}{0.000000,0.000000,0.000000}%
\pgfsetfillcolor{currentfill}%
\pgfsetfillopacity{0.000000}%
\pgfsetlinewidth{1.003750pt}%
\definecolor{currentstroke}{rgb}{0.000000,0.000000,0.000000}%
\pgfsetstrokecolor{currentstroke}%
\pgfsetdash{}{0pt}%
\pgfsys@defobject{currentmarker}{\pgfqpoint{-0.058926in}{-0.058926in}}{\pgfqpoint{0.058926in}{0.058926in}}{%
\pgfpathmoveto{\pgfqpoint{-0.000000in}{-0.058926in}}%
\pgfpathlineto{\pgfqpoint{0.058926in}{0.000000in}}%
\pgfpathlineto{\pgfqpoint{0.000000in}{0.058926in}}%
\pgfpathlineto{\pgfqpoint{-0.058926in}{0.000000in}}%
\pgfpathlineto{\pgfqpoint{-0.000000in}{-0.058926in}}%
\pgfpathclose%
\pgfusepath{stroke,fill}%
}%
\begin{pgfscope}%
\pgfsys@transformshift{2.868426in}{0.499691in}%
\pgfsys@useobject{currentmarker}{}%
\end{pgfscope}%
\begin{pgfscope}%
\pgfsys@transformshift{2.977310in}{0.619194in}%
\pgfsys@useobject{currentmarker}{}%
\end{pgfscope}%
\begin{pgfscope}%
\pgfsys@transformshift{3.086194in}{0.848663in}%
\pgfsys@useobject{currentmarker}{}%
\end{pgfscope}%
\begin{pgfscope}%
\pgfsys@transformshift{3.195079in}{1.223449in}%
\pgfsys@useobject{currentmarker}{}%
\end{pgfscope}%
\begin{pgfscope}%
\pgfsys@transformshift{3.303963in}{1.108826in}%
\pgfsys@useobject{currentmarker}{}%
\end{pgfscope}%
\begin{pgfscope}%
\pgfsys@transformshift{3.412847in}{1.400588in}%
\pgfsys@useobject{currentmarker}{}%
\end{pgfscope}%
\begin{pgfscope}%
\pgfsys@transformshift{3.521732in}{1.264895in}%
\pgfsys@useobject{currentmarker}{}%
\end{pgfscope}%
\begin{pgfscope}%
\pgfsys@transformshift{3.630616in}{1.450261in}%
\pgfsys@useobject{currentmarker}{}%
\end{pgfscope}%
\begin{pgfscope}%
\pgfsys@transformshift{3.739500in}{1.351769in}%
\pgfsys@useobject{currentmarker}{}%
\end{pgfscope}%
\begin{pgfscope}%
\pgfsys@transformshift{3.848384in}{1.364721in}%
\pgfsys@useobject{currentmarker}{}%
\end{pgfscope}%
\begin{pgfscope}%
\pgfsys@transformshift{3.957269in}{1.406220in}%
\pgfsys@useobject{currentmarker}{}%
\end{pgfscope}%
\begin{pgfscope}%
\pgfsys@transformshift{4.066153in}{1.399197in}%
\pgfsys@useobject{currentmarker}{}%
\end{pgfscope}%
\begin{pgfscope}%
\pgfsys@transformshift{4.175037in}{1.418462in}%
\pgfsys@useobject{currentmarker}{}%
\end{pgfscope}%
\begin{pgfscope}%
\pgfsys@transformshift{4.283922in}{1.405387in}%
\pgfsys@useobject{currentmarker}{}%
\end{pgfscope}%
\begin{pgfscope}%
\pgfsys@transformshift{4.392806in}{1.419227in}%
\pgfsys@useobject{currentmarker}{}%
\end{pgfscope}%
\begin{pgfscope}%
\pgfsys@transformshift{4.501690in}{1.405431in}%
\pgfsys@useobject{currentmarker}{}%
\end{pgfscope}%
\end{pgfscope}%
\begin{pgfscope}%
\pgfpathrectangle{\pgfqpoint{2.786763in}{0.499691in}}{\pgfqpoint{1.796591in}{2.552550in}}%
\pgfusepath{clip}%
\pgfsetbuttcap%
\pgfsetroundjoin%
\pgfsetlinewidth{1.104125pt}%
\definecolor{currentstroke}{rgb}{0.000000,0.000000,0.000000}%
\pgfsetstrokecolor{currentstroke}%
\pgfsetdash{{1.100000pt}{1.815000pt}}{0.000000pt}%
\pgfpathmoveto{\pgfqpoint{2.868426in}{1.399252in}}%
\pgfpathlineto{\pgfqpoint{2.977310in}{1.399252in}}%
\pgfpathlineto{\pgfqpoint{3.086194in}{1.399252in}}%
\pgfpathlineto{\pgfqpoint{3.195079in}{1.399252in}}%
\pgfpathlineto{\pgfqpoint{3.303963in}{1.399252in}}%
\pgfpathlineto{\pgfqpoint{3.412847in}{1.399252in}}%
\pgfpathlineto{\pgfqpoint{3.521732in}{1.399252in}}%
\pgfpathlineto{\pgfqpoint{3.630616in}{1.399252in}}%
\pgfpathlineto{\pgfqpoint{3.739500in}{1.399252in}}%
\pgfpathlineto{\pgfqpoint{3.848384in}{1.399252in}}%
\pgfpathlineto{\pgfqpoint{3.957269in}{1.399252in}}%
\pgfpathlineto{\pgfqpoint{4.066153in}{1.399252in}}%
\pgfpathlineto{\pgfqpoint{4.175037in}{1.399252in}}%
\pgfpathlineto{\pgfqpoint{4.283922in}{1.399252in}}%
\pgfpathlineto{\pgfqpoint{4.392806in}{1.399252in}}%
\pgfpathlineto{\pgfqpoint{4.501690in}{1.399252in}}%
\pgfusepath{stroke}%
\end{pgfscope}%
\begin{pgfscope}%
\pgfsetrectcap%
\pgfsetmiterjoin%
\pgfsetlinewidth{0.803000pt}%
\definecolor{currentstroke}{rgb}{0.000000,0.000000,0.000000}%
\pgfsetstrokecolor{currentstroke}%
\pgfsetdash{}{0pt}%
\pgfpathmoveto{\pgfqpoint{2.786763in}{0.499691in}}%
\pgfpathlineto{\pgfqpoint{2.786763in}{3.052241in}}%
\pgfusepath{stroke}%
\end{pgfscope}%
\begin{pgfscope}%
\pgfsetrectcap%
\pgfsetmiterjoin%
\pgfsetlinewidth{0.803000pt}%
\definecolor{currentstroke}{rgb}{0.000000,0.000000,0.000000}%
\pgfsetstrokecolor{currentstroke}%
\pgfsetdash{}{0pt}%
\pgfpathmoveto{\pgfqpoint{4.583353in}{0.499691in}}%
\pgfpathlineto{\pgfqpoint{4.583353in}{3.052241in}}%
\pgfusepath{stroke}%
\end{pgfscope}%
\begin{pgfscope}%
\pgfsetrectcap%
\pgfsetmiterjoin%
\pgfsetlinewidth{0.803000pt}%
\definecolor{currentstroke}{rgb}{0.000000,0.000000,0.000000}%
\pgfsetstrokecolor{currentstroke}%
\pgfsetdash{}{0pt}%
\pgfpathmoveto{\pgfqpoint{2.786763in}{0.499691in}}%
\pgfpathlineto{\pgfqpoint{4.583353in}{0.499691in}}%
\pgfusepath{stroke}%
\end{pgfscope}%
\begin{pgfscope}%
\pgfsetrectcap%
\pgfsetmiterjoin%
\pgfsetlinewidth{0.803000pt}%
\definecolor{currentstroke}{rgb}{0.000000,0.000000,0.000000}%
\pgfsetstrokecolor{currentstroke}%
\pgfsetdash{}{0pt}%
\pgfpathmoveto{\pgfqpoint{2.786763in}{3.052241in}}%
\pgfpathlineto{\pgfqpoint{4.583353in}{3.052241in}}%
\pgfusepath{stroke}%
\end{pgfscope}%
\begin{pgfscope}%
\pgfsetbuttcap%
\pgfsetmiterjoin%
\definecolor{currentfill}{rgb}{1.000000,1.000000,1.000000}%
\pgfsetfillcolor{currentfill}%
\pgfsetfillopacity{0.800000}%
\pgfsetlinewidth{1.003750pt}%
\definecolor{currentstroke}{rgb}{0.800000,0.800000,0.800000}%
\pgfsetstrokecolor{currentstroke}%
\pgfsetstrokeopacity{0.800000}%
\pgfsetdash{}{0pt}%
\pgfpathmoveto{\pgfqpoint{2.930959in}{2.201161in}}%
\pgfpathlineto{\pgfqpoint{4.486131in}{2.201161in}}%
\pgfpathquadraticcurveto{\pgfqpoint{4.513909in}{2.201161in}}{\pgfqpoint{4.513909in}{2.228939in}}%
\pgfpathlineto{\pgfqpoint{4.513909in}{2.955019in}}%
\pgfpathquadraticcurveto{\pgfqpoint{4.513909in}{2.982797in}}{\pgfqpoint{4.486131in}{2.982797in}}%
\pgfpathlineto{\pgfqpoint{2.930959in}{2.982797in}}%
\pgfpathquadraticcurveto{\pgfqpoint{2.903182in}{2.982797in}}{\pgfqpoint{2.903182in}{2.955019in}}%
\pgfpathlineto{\pgfqpoint{2.903182in}{2.228939in}}%
\pgfpathquadraticcurveto{\pgfqpoint{2.903182in}{2.201161in}}{\pgfqpoint{2.930959in}{2.201161in}}%
\pgfpathlineto{\pgfqpoint{2.930959in}{2.201161in}}%
\pgfpathclose%
\pgfusepath{stroke,fill}%
\end{pgfscope}%
\begin{pgfscope}%
\pgfsetbuttcap%
\pgfsetmiterjoin%
\definecolor{currentfill}{rgb}{0.000000,0.000000,0.000000}%
\pgfsetfillcolor{currentfill}%
\pgfsetfillopacity{0.000000}%
\pgfsetlinewidth{1.003750pt}%
\definecolor{currentstroke}{rgb}{0.000000,0.000000,0.000000}%
\pgfsetstrokecolor{currentstroke}%
\pgfsetdash{}{0pt}%
\pgfsys@defobject{currentmarker}{\pgfqpoint{-0.041667in}{-0.041667in}}{\pgfqpoint{0.041667in}{0.041667in}}{%
\pgfpathmoveto{\pgfqpoint{0.000000in}{0.041667in}}%
\pgfpathlineto{\pgfqpoint{-0.041667in}{-0.041667in}}%
\pgfpathlineto{\pgfqpoint{0.041667in}{-0.041667in}}%
\pgfpathlineto{\pgfqpoint{0.000000in}{0.041667in}}%
\pgfpathclose%
\pgfusepath{stroke,fill}%
}%
\begin{pgfscope}%
\pgfsys@transformshift{3.097626in}{2.878630in}%
\pgfsys@useobject{currentmarker}{}%
\end{pgfscope}%
\end{pgfscope}%
\begin{pgfscope}%
\definecolor{textcolor}{rgb}{0.000000,0.000000,0.000000}%
\pgfsetstrokecolor{textcolor}%
\pgfsetfillcolor{textcolor}%
\pgftext[x=3.347626in,y=2.830019in,left,base]{\color{textcolor}\rmfamily\fontsize{10.000000}{12.000000}\selectfont Nonlinear CG}%
\end{pgfscope}%
\begin{pgfscope}%
\pgfsetbuttcap%
\pgfsetmiterjoin%
\definecolor{currentfill}{rgb}{0.000000,0.000000,0.000000}%
\pgfsetfillcolor{currentfill}%
\pgfsetfillopacity{0.000000}%
\pgfsetlinewidth{1.003750pt}%
\definecolor{currentstroke}{rgb}{0.000000,0.000000,0.000000}%
\pgfsetstrokecolor{currentstroke}%
\pgfsetdash{}{0pt}%
\pgfsys@defobject{currentmarker}{\pgfqpoint{-0.058926in}{-0.058926in}}{\pgfqpoint{0.058926in}{0.058926in}}{%
\pgfpathmoveto{\pgfqpoint{-0.000000in}{-0.058926in}}%
\pgfpathlineto{\pgfqpoint{0.058926in}{0.000000in}}%
\pgfpathlineto{\pgfqpoint{0.000000in}{0.058926in}}%
\pgfpathlineto{\pgfqpoint{-0.058926in}{0.000000in}}%
\pgfpathlineto{\pgfqpoint{-0.000000in}{-0.058926in}}%
\pgfpathclose%
\pgfusepath{stroke,fill}%
}%
\begin{pgfscope}%
\pgfsys@transformshift{3.097626in}{2.684957in}%
\pgfsys@useobject{currentmarker}{}%
\end{pgfscope}%
\end{pgfscope}%
\begin{pgfscope}%
\definecolor{textcolor}{rgb}{0.000000,0.000000,0.000000}%
\pgfsetstrokecolor{textcolor}%
\pgfsetfillcolor{textcolor}%
\pgftext[x=3.347626in,y=2.636346in,left,base]{\color{textcolor}\rmfamily\fontsize{10.000000}{12.000000}\selectfont Exact trust region}%
\end{pgfscope}%
\begin{pgfscope}%
\pgfsetbuttcap%
\pgfsetroundjoin%
\pgfsetlinewidth{1.104125pt}%
\definecolor{currentstroke}{rgb}{0.000000,0.000000,0.000000}%
\pgfsetstrokecolor{currentstroke}%
\pgfsetdash{{1.100000pt}{1.815000pt}}{0.000000pt}%
\pgfpathmoveto{\pgfqpoint{2.958737in}{2.398692in}}%
\pgfpathlineto{\pgfqpoint{3.097626in}{2.398692in}}%
\pgfpathlineto{\pgfqpoint{3.236515in}{2.398692in}}%
\pgfusepath{stroke}%
\end{pgfscope}%
\begin{pgfscope}%
\definecolor{textcolor}{rgb}{0.000000,0.000000,0.000000}%
\pgfsetstrokecolor{textcolor}%
\pgfsetfillcolor{textcolor}%
\pgftext[x=3.347626in, y=2.443445in, left, base]{\color{textcolor}\rmfamily\fontsize{10.000000}{12.000000}\selectfont Rosenbrock valley}%
\end{pgfscope}%
\begin{pgfscope}%
\definecolor{textcolor}{rgb}{0.000000,0.000000,0.000000}%
\pgfsetstrokecolor{textcolor}%
\pgfsetfillcolor{textcolor}%
\pgftext[x=3.347626in, y=2.291439in, left, base]{\color{textcolor}\rmfamily\fontsize{10.000000}{12.000000}\selectfont subspace at \((1, 1)\)}%
\end{pgfscope}%
\end{pgfpicture}%
\makeatother%
\endgroup%